# A Simple Galerkin Meshless Method, the Fragile Points Method (FPM) Using Point Stiffness Matrices, for 2D Linear Elastic Problems in Complex Domains with Crack and Rupture Propagation

Tian Yang[1], Leiting Dong[2,*], Satya N. Atluri[3]

[1]*Graduate student, School of Aeronautic Science and Engineering, Beihang University, China*

[2]*Professor, School of Aeronautic Science and Engineering, Beihang University, China*

[3]*Presidential Chair & University Distinguished Professor of Texas Tech University, USA*

## Abstract

The Fragile Points Method (FPM) is an elementarily simple Galerkin meshless method, employing Point-based discontinuous trial and test functions only, *without using element-based trial and test functions*. In this study, the algorithmic formulations of FPM for linear elasticity are given in detail, by exploring the concepts of point stiffness matrices and numerical flux corrections. Advantages of FPM for simulating the deformations of complex structures, and for simulating complex crack propagations and rupture developments, are also thoroughly discussed. Numerical examples of deformation and stress analyses of benchmark problems, as well as of realistic structures with complex geometries, demonstrate the accuracy, efficiency and robustness of the proposed FPM. Simulations of crack initiation and propagations are also given in this study, demonstrating the advantages of the present FPM in modeling complex rupture and fracture phenomena. The crack and rupture propagation modeling in FPM is achieved without remeshing or augmenting the trial functions as in standard, extended or generalized FEM. The simulation of impact, penetration and other extreme problems by FPM will be discussed in our future papers.
KEY WORDS: Elasticity; Meshfree methods; Fragile Points Method; Numerical Flux Corrections; Fracture

* Corresponding author. Email address: ltdong@buaa.edu.cn (L. Dong).

# 1. Introduction

Structural stress analysis is crucial and necessary in diverse engineering fields, such as aeronautics, astronautics, automobile engineering, etc. From the design and manufacture to maintenance of products, structural stress analysis plays a crucial role. Because of its significance, numerous researchers have been focusing on improving the accuracy and efficiency of this procedure for decades. Moreover, under certain extreme conditions, crack initiation and propagation would result in serious deterioration to the integrity of the structure. Therefore, efficient and accurate simulations of the deformation, stress, as well as the crack initiation and propagation are of significant importance.

The Finite Element Method (FEM) is mature, reliable and widely used in structural stress analysis [1]. This method employs contiguous elements, and Element-based, local, polynomial, interelement-continuous trial and test functions. Because the trial and test functions are Element-based, the Galerkin weak form leads to Element Stiffness Matrices. Therefore, integrals in the Galerkin weak form underlying the FEM are easy to compute. The symmetry and sparsity of the global stiffness matrix make the FEM suitable and efficient in large-scale simulations. However, the accuracy of the FEM greatly depends on the quality of mesh. In order to obtain satisfactory solutions, many efforts are usually made on meshing. Especially, even if simulations are initialized with a high-quality mesh structure, mesh distortion will occur in the case of large deformations and the precision of solutions decreases dramatically. In order to study the formation of cracks, rupture and fragmentation, methods such as remeshing, and deleting elements, are often used.

Meshless methods, which eliminate the mesh structure partly or completely, have been invented and developed since the end of last century. Element Free Galerkin (EFG) [2] and Meshless Local Petrov-Galerkin (MLPG) [3] methods are two classical meshless weak-form methods based on the "Global Galerkin" and "Local Petrov-Galerkin" weak forms, respectively. While the EFG method uses the same Node-based trial and test functions, the MLPG method uses different local trial and test function spaces. These two meshless methods have utilized Moving Least Squares (MLS), Radial Basis Function (RBF), or other methods to derive Node-based trial functions. With MLS and RBF approximations, higher-order continuity can be easily achieved. Besides, since individual nodes have replaced element-based mesh structure, EFG and MLPG can conveniently bypass the influence of mesh distortion even in large deformation and fracture simulations (e.g. [4], [5]). However, on the other hand, the trial functions given by MLS or RBF are rational functions and grossly complex. Therefore, the computation of integrals in the weak forms in either EFG or MLPG is very tedious, less accurate and can influence the method's stability. To reduce the computational cost and improve the accuracy of integration, some special, new types of numerical integration methods, for example, the series of nodal integration methods [6], are often adopted.

Smoothed Particle Hydrodynamics (SPH) method [7], as a kind of meshless particle method,  needs less computational cost. Nodal smoothing together with the collocation of governing differential equations are used to derive discretized algebraic equations. This, on one hand, makes it very simple and easy for implementation. On the other hand,  proving the stability of a strong form method is not an easy task. In fact, tensile instability will occur in the SPH, if we use Smoothed Kernel functions to calculate derivatives.

From the above discussion, we can conclude that simple, local, polynomial, "Point-Based" shape functions are helpful in the calculation of integrals in the weak form. Besides, a weak-form method can have a better performance on stability. But with these requirements, it is difficult to keep the trial and test functions continuous over the entire domain. In our previous paper, we have developed the Fragile Points Method [8] for the first time, for Poisson's equations. The FPM approach employs Point-based and discontinuous trial and test functions (piece-wise polynomials) instead of continuous ones. Substituting the “Point-Based” functions in a Galerkin weak form, the method leads to “Point Stiffness Matrices” as opposed to the Element Stiffness Matrices in the FEM. Numerical Flux Corrections are introduced in the FPM to solve the inconsistency caused by the method’s discontinuity of trial and test functions. Integrals in the Galerkin weak form can be computed easily by using Gauss Integration or even just analytically. Like the FEM, since the FPM is based on a Galerkin weak form, a symmetric, sparse and positive definitive global matrix can be obtained in the FPM, which means that the FPM can be easily used in large-scale simulations. More importantly, this method is named as “Fragile Points Method” because we can easily cut off the interaction between two neighboring Points and introduce cracks, to simulate damage, rupture, fragmentation, material “Fragility” without much effort. The major differences between FPM and other mesh-based and meshless methods are listed in Table I.

Table I. The difference between FPM and other methods

| Methods | Trial Functions | Weak/Strong forms | Numerical Integration |
|---|---|---|---|
| FEM | Element-based interpolations | Galerkin weak form | Gauss integration for Isoparametric Elements, inaccurate when the element is highly distorted |
| EFG | Point-based continuous trial functions (MLS, RBF etc.) | Galerkin weak form | Numerical integration with many quadrature points, or stabilized nodal integration with cell-based smoothing |
| MLPG | Point-based continuous trial functions (MLS, RBF etc.) | Local Petrov-Galerkin weak form | Numerical integration in the local subdomain with many quadrature points |
| SPH | Kernel smoothing of nodal variables | Strong form | Collocation for spatial discretization, may cause |
| | | | instability |
| FPM | Point-based discontinuous trial functions (piece-wise polynomial) | Galerkin weak form with Numerical Flux corrections | Exact integration with simple Gauss integration (one-point-integral for linear trial functions), as only simple Cartesian strains are involved |

In this paper, we formulate and apply the FPM for solving linear elastic problems in complex shaped domains, and also for simulating crack and rupture initiation and propagation. The procedure of constructing Point-based trial and test functions is introduced in Section 2. The Interior Penalty Numerical Fluxes and the numerical implementation of the FPM for elasticity are discussed in Section 3. Detailed steps to deal with cracks in the FPM are also introduced in Section 3. Several benchmark problems as well as realistic structures with complex geometries, involving stress and deformation analyses, are studied in Section 4. Simulations involving crack initiation and propagation and their comparison with experimental results are also presented in Section 4. Lastly, conclusions and some discussions are given in Section 5.

## 2. Local, Polynomial, Point-Based, Discontinuous Trial and Test Functions

For linear elasticity, the governing equations are given in Eq. (2.1),

$$\begin{cases} \varepsilon_{ij}(\mathbf{u}) = \frac{1}{2}\left(u_{i,j} + u_{j,i}\right) \\ \sigma_{ij,j}(\mathbf{u}) + f_i = 0 \qquad \text{in } \Omega \\ \sigma_{ij}(\mathbf{u}) = D_{ijkl}\varepsilon_{kl}(\mathbf{u}) \end{cases} \tag{2.1}$$

where $\Omega$ is the problem domain; $\sigma_{ij}$, $\varepsilon_{ij}$ and $u_i$ stand for the stress tensor, strain tensor and displacement vector, respectively; $f_i$ is the body force and $D_{ijkl}$ is the fourth order linear elasticity tensor.

The corresponding boundary conditions are shown in Eq. (2.2), where $\Gamma_u$ and $\Gamma_t$ are displacement prescribed and traction prescribed boundaries, respectively; $\bar{u}_i$ and $\bar{t}_i$ denote the prescribed displacements and tractions on the corresponding boundaries, respectively; $n_j$ stands for the unit vector outward to the external boundary $\partial\Omega$.

$$\begin{cases} u_i = \bar{u}_i & \text{on } \Gamma_u \\ \sigma_{ij}(\mathbf{u}) n_j = \bar{t}_i & \text{on } \Gamma_t \end{cases} \tag{2.2}$$

Considering the problem domain $\Omega$, as shown in Figure 1(a), several Points are distributed randomly inside the domain or on its boundary. Utilizing these Points, the domain can be partitioned into contiguous and nonoverlapping subdomains of arbitrary shape, with only one Point being involved in each subdomain (shown in Figure 1(b)). It should be noted that such an initial domain partition is necessary for FPM, and numerous methods can be used for this partition. And in this paper, the Voronoi Diagram method is employed as a simple choice. One can also convert the contiguous elements used for FEM to be the subdomains which are needed for FPM (shown in Figure 2).

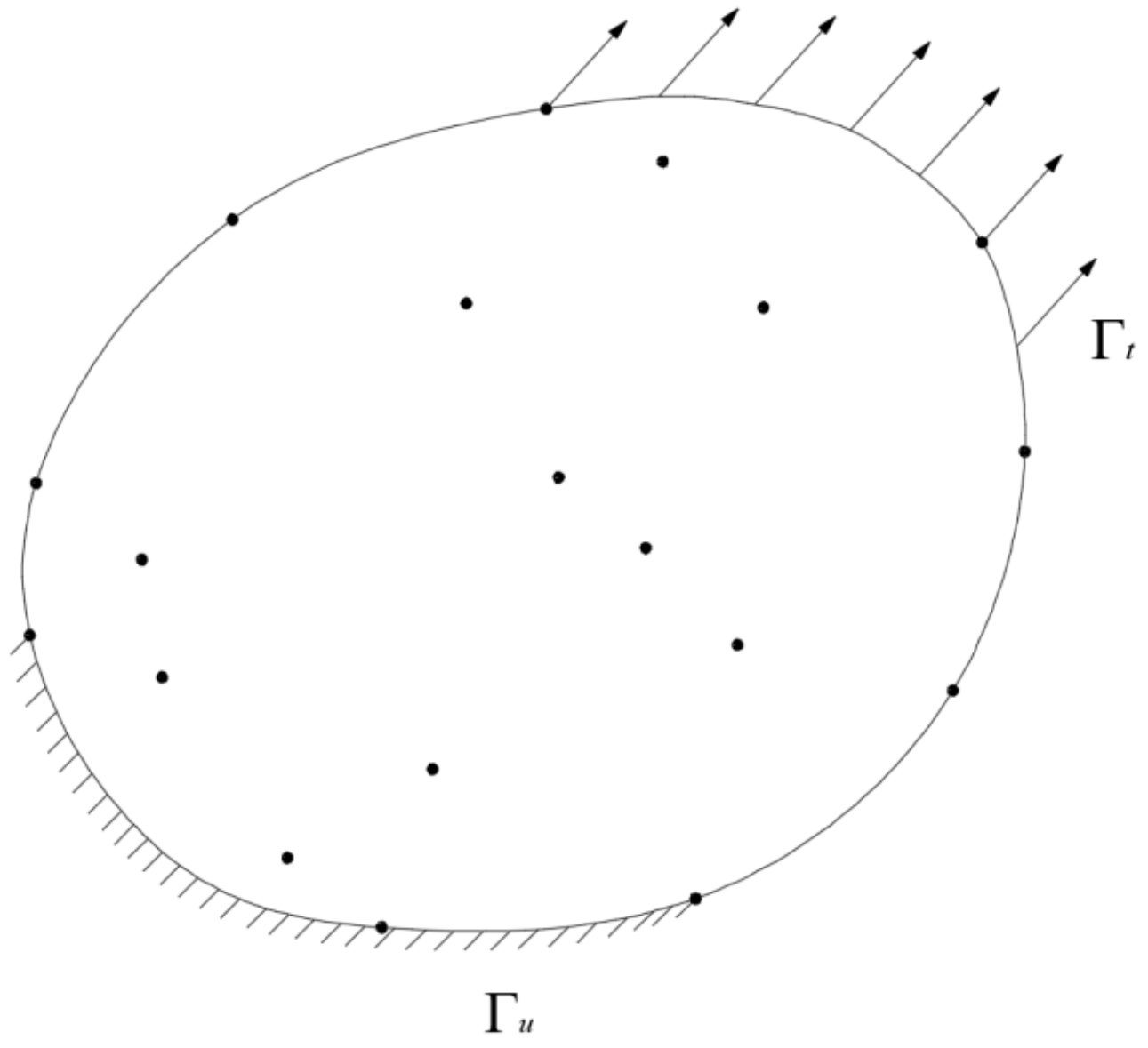

(a)

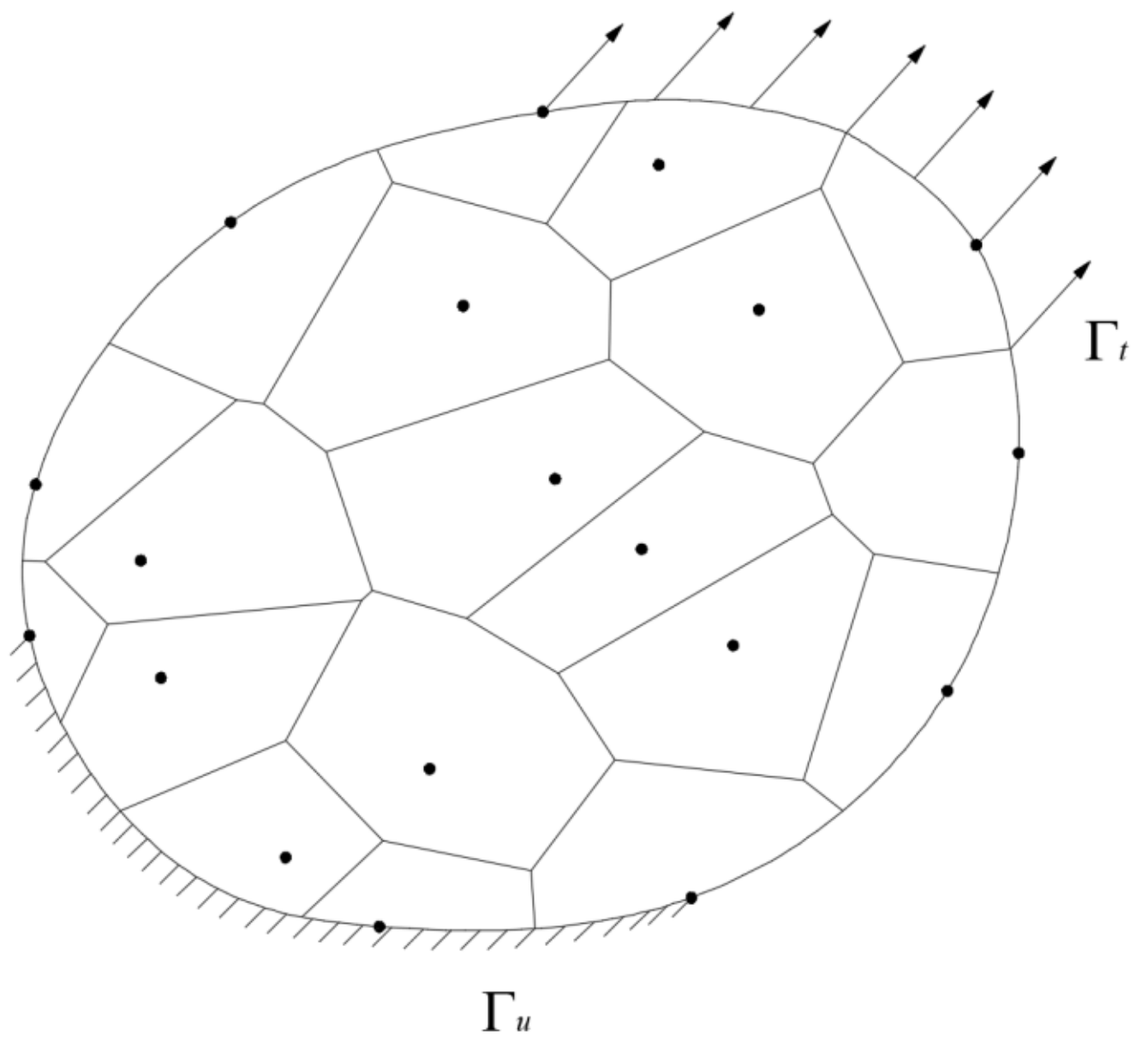

(b)

Figure 1 (a) (b). The problem domain and its partition

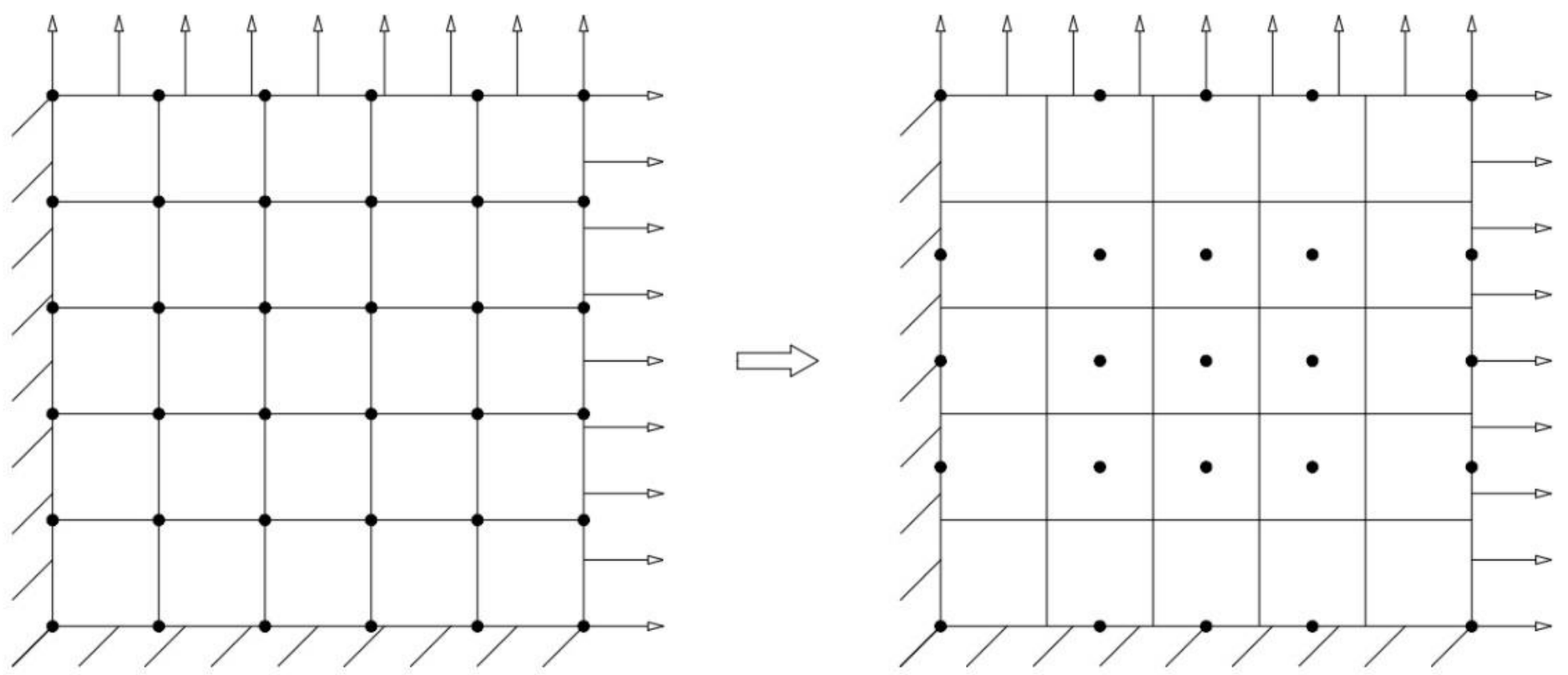

Figure 2. Subdomains based on a FEM mesh

In each subdomain, we define the simple, local, polynomial, discontinuous displacement vector or trial function $\mathbf{u}^h$ using Taylor's expansion, in terms of the values $u_1$ and $u_2$ (displacements in $x_1$ and $x_2$ directions, respectively) and their derivatives, at the internal Point. For instance, the approximation to the local displacement field within the subdomain $E_0$ which contains the Point $P_0$ is given in Eq. (2.3),

$$\mathbf{u}^h(x_1,x_2)=\begin{bmatrix}u_1^h\\u_2^h\end{bmatrix}=\begin{bmatrix}u_1^0\\u_2^0\end{bmatrix}+\begin{bmatrix}\mathbf{h}(x_1,x_2) & 0\\ 0 & \mathbf{h}(x_1,x_2)\end{bmatrix}\begin{bmatrix}\mathbf{a}_1\\\mathbf{a}_2\end{bmatrix},\quad (x_1,x_2)\in E_0. \tag{2.3}$$

For linear trial functions,

$$\mathbf{h}(x_1,x_2)=\begin{bmatrix}x_1-x_1^0 & x_2-x_2^0\end{bmatrix},\ \mathbf{a}_i=\begin{bmatrix}\dfrac{\partial u_i}{\partial x_1} & \dfrac{\partial u_i}{\partial x_2}\end{bmatrix}^{\mathrm{T}}\Bigg|_{P_0}\quad i=1,2.$$

For quadratic trial functions,

$$\mathbf{h}(x_1,x_2)=\begin{bmatrix}x_1-x_1^0 & x_2-x_2^0 & \dfrac{1}{2}(x_1-x_1^0)^2 & \dfrac{1}{2}(x_2-x_2^0)^2 & (x_1-x_1^0)(x_2-x_2^0)\end{bmatrix},$$

$$\mathbf{a}_i=\begin{bmatrix}\dfrac{\partial u_i}{\partial x_1} & \dfrac{\partial u_i}{\partial x_2} & \dfrac{\partial^2 u_i}{\partial x_1^{\ 2}} & \dfrac{\partial^2 u_i}{\partial x_2^{\ 2}} & \dfrac{\partial^2 u_i}{\partial x_1\partial x_2}\end{bmatrix}^{\mathrm{T}}\Bigg|_{P_0}\quad i=1,2,$$

where $(x_1^0, x_2^0)$ are the coordinates of the Point $P_0$; $\begin{bmatrix}u_1^0 & u_2^0\end{bmatrix}^{\mathrm{T}}$ are nodal degrees of freedom at $P_0$; Derivatives $\mathbf{a}^{\mathrm{T}}=\begin{bmatrix}\mathbf{a}_1^{\ \mathrm{T}} & \mathbf{a}_2^{\ \mathrm{T}}\end{bmatrix}$ are to be related to a finite number of nodal DoFs in the support domain of $P_0$.

Since the derivatives $\mathbf{a}^{\mathrm{T}}=\begin{bmatrix}\mathbf{a}_1^{\ \mathrm{T}} & \mathbf{a}_2^{\ \mathrm{T}}\end{bmatrix}$ are only determined at each point, numerous approaches may be employed. Specifically, we use the Generalized Finite Difference (GFD) method [9] and the Compactly-Supported Radial Basis Function (CSRBF) method [10], respectively in the following study.

## 2.1 The Generalized Finite Difference method

The first step for the GFD method is to define the local support of the Point $P_0$. Usually, we prefer to define the support by drawing a circle at $P_0$ and assume that all the Points included in that circle have interactions with $P_0$ (shown in Figure 3(a)). Alternatively, we can replace the circular support with a square one or other shapes. In this paper, the support of $P_0$ is defined so as to contain all of its neighboring points in the subdomain partition (shown in Figure 3(b)). These neighboring points are named as

$P_1,\ P_2, \ldots,\ P_m$.

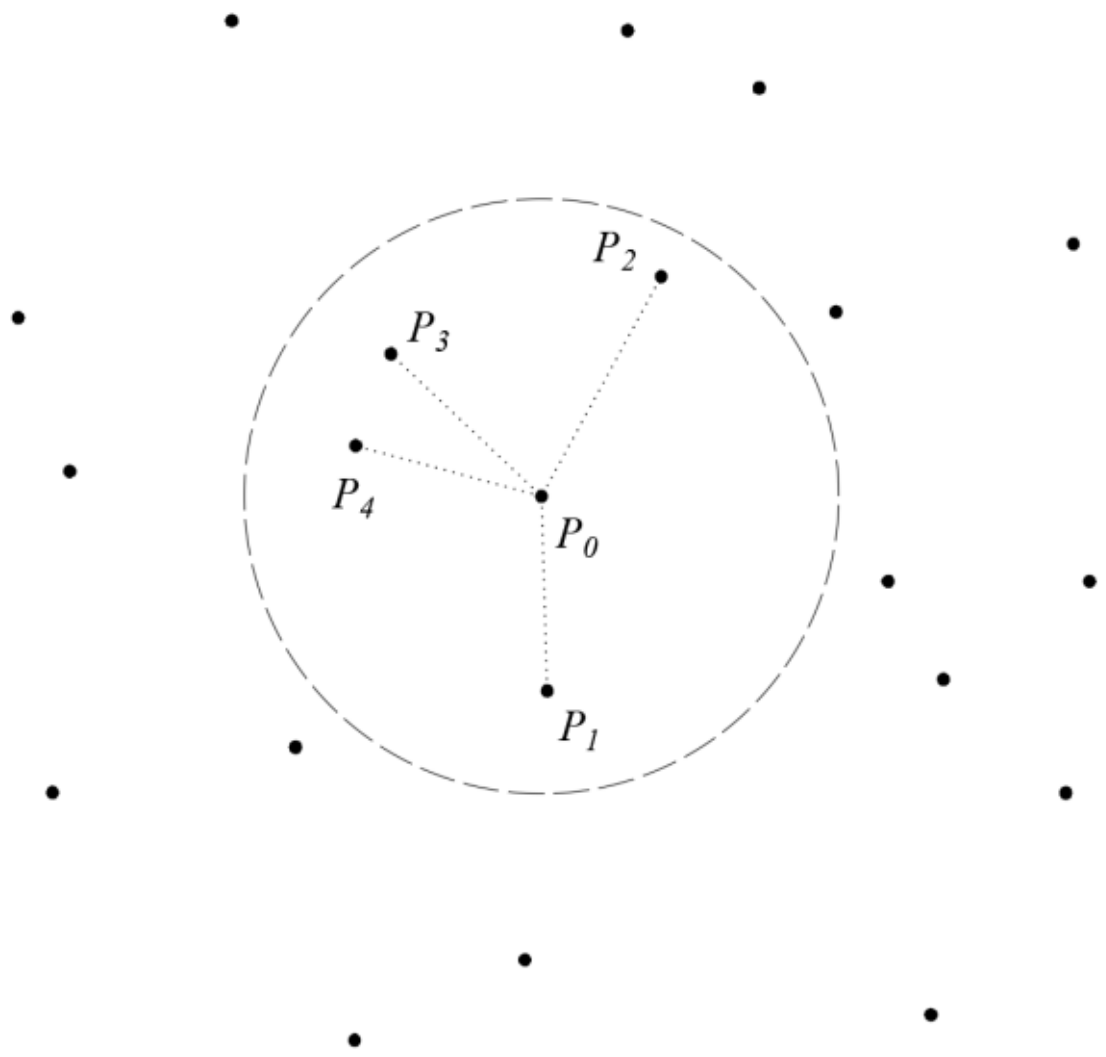

(a)

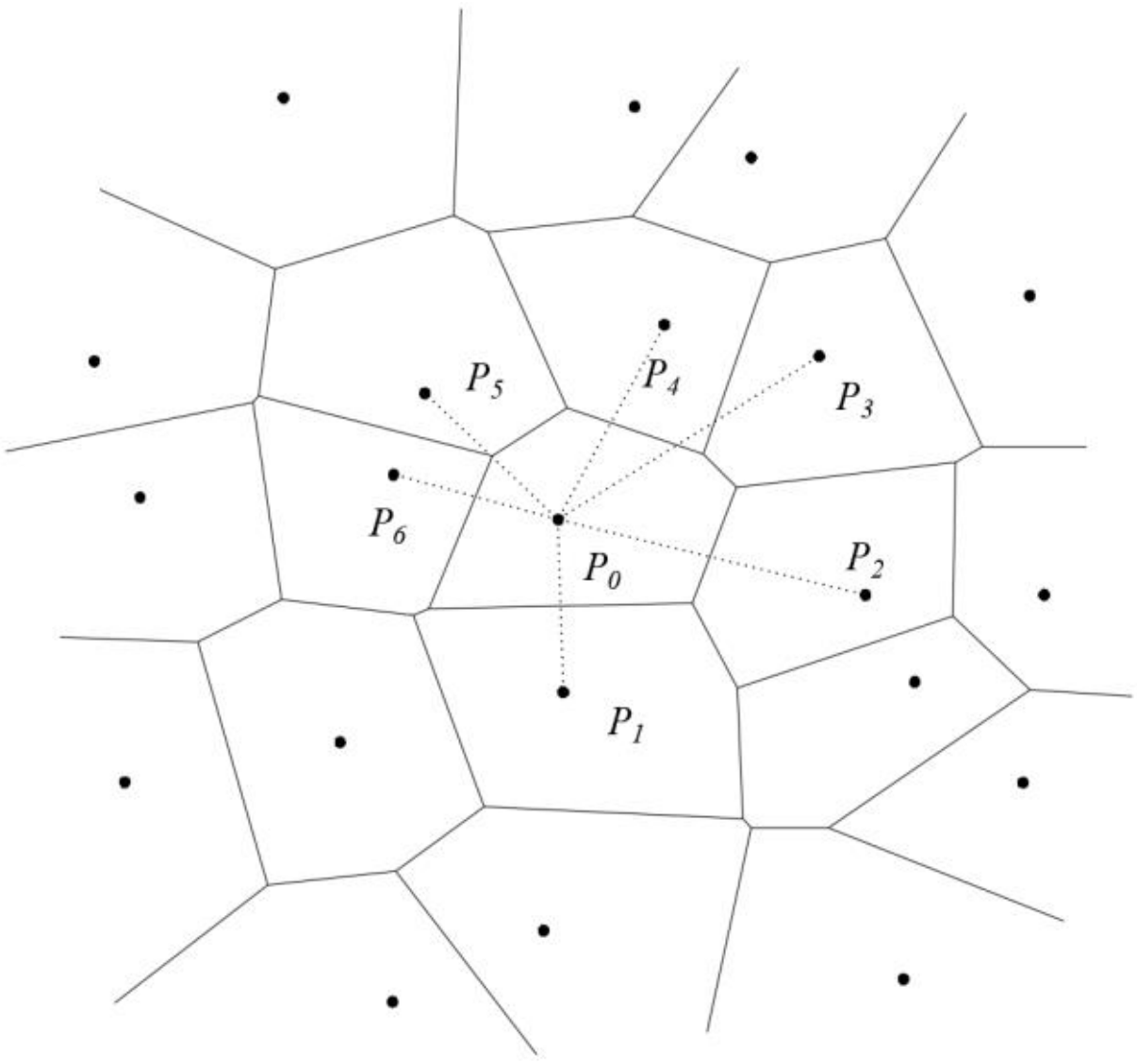

(b)

Figure 3 (a) (b). Two kinds of support of $P_0$

After defining the support of $P_0$, we define a weighted discrete $L^2$ norm $J_1$ in a matrix form,

$$J_1 = \left(\mathbf{A}\mathbf{a}_1 + \mathbf{u}_1^0 - \mathbf{u}_1^m\right)^{\mathrm{T}} \mathbf{W}_1 \left(\mathbf{A}\mathbf{a}_1 + \mathbf{u}_1^0 - \mathbf{u}_1^m\right) \tag{2.4}$$

where

$$\mathbf{u}_1^0 = \left(\begin{bmatrix} u_1^0 & u_1^0 & \ldots & u_1^0 \end{bmatrix}_{1\times m}\right)^{\mathrm{T}},\ \mathbf{u}_1^m = \begin{bmatrix} u_1^1 & u_1^2 & \ldots & u_1^m \end{bmatrix}^{\mathrm{T}}$$

$$\mathbf{A}=\begin{bmatrix}\mathbf{h}(\mathbf{x}_1)\\ \dots\\ \mathbf{h}(\mathbf{x}_m)\end{bmatrix}\quad \mathbf{W}_1=\begin{bmatrix} w_1 & 0 & 0 & \dots & 0\\ 0 & w_2 & 0 & \dots & 0\\ \dots & \dots & \dots & \dots & \dots\\ 0 & \dots & 0 & w_{m-1} & 0\\ 0 & \dots & \dots & 0 & w_m\end{bmatrix}$$

$\mathbf{x}_i=\left(x_1^i, x_2^i\right)$ are the coordinates of $P_i$; $u_1^i$ is the value of $u_1$ at $P_i$; $w_i$ is the value of the weight function at $P_i$ $(i=1,2,3,\dots,m)$. For convenience, constant weight function is used in this paper.

By considering the stationarity condition of $J_1$ in Eq. (2.4), we can relate $\mathbf{a}_1$ at $P_0$ to nodal DoFs:

$$\mathbf{a}_1=\left(\mathbf{A}^{\mathrm{T}}\mathbf{W}_1\mathbf{A}\right)^{-1}\mathbf{A}^{\mathrm{T}}\mathbf{W}_1\left(\mathbf{u}_1^m-\mathbf{u}_1^0\right) \tag{2.5}$$

Besides, $\mathbf{u}_1^m-\mathbf{u}_1^0$ can be transformed into the following form,

$$\left(\mathbf{u}_1^m-\mathbf{u}_1^0\right)=\left[\mathbf{I}_1 \quad \mathbf{I}_2\right]\mathbf{u}_1 \tag{2.6}$$

where $\mathbf{u}_1=\left[u_1^0 \quad u_1^1 \quad \dots \quad u_1^m\right]^{\mathrm{T}}$, $\mathbf{I}_1=\left(\left[-1 \quad -1 \quad \dots \quad -1\right]_{1\times m}\right)^{\mathrm{T}}$,

$$\mathbf{I}_2=\begin{bmatrix}1 & 0 & \dots & \dots & 0\\ 0 & 1 & 0 & \dots & 0\\ \dots & \dots & \dots & \dots & \dots\\ 0 & \dots & 0 & 1 & 0\\ 0 & \dots & \dots & 0 & 1\end{bmatrix}_{m\times m}$$

Substituting Eq. (2.6) into Eq. (2.5), we obtain the relation between $\mathbf{a}_1$ and $\mathbf{u}_1$.

$$\mathbf{a}_1=\mathbf{C}_1\mathbf{u}_1$$

where

$$\mathbf{C}_1=\left(\mathbf{A}^{\mathrm{T}}\mathbf{W}_1\mathbf{A}\right)^{-1}\mathbf{A}^{\mathrm{T}}\mathbf{W}_1\left[\mathbf{I}_1 \quad \mathbf{I}_2\right]$$

The vector of derivatives $\mathbf{a}_2$ for displacement $u_2$ is calculated in the same way:

$$\mathbf{a}_2=\mathbf{C}_2\mathbf{u}_2 \text{ where } \mathbf{u}_2=\left[u_2^0 \quad u_2^1 \quad \dots \quad u_2^m\right]^{\mathrm{T}}, \text{and } \mathbf{C}_2=\mathbf{C}_1 \text{ if } \mathbf{W}_2=\mathbf{W}_1$$

By combining derivatives for $u_1$ and $u_2$, we can obtain

$$\mathbf{a}=\mathbf{C}\mathbf{u}_E,\quad \mathbf{u}_E=\left[u_1^0 \quad u_2^0 \quad u_1^1 \quad u_2^1 \quad \dots \quad u_1^m \quad u_2^m\right]^{\mathrm{T}} \tag{2.7}$$

In this way, we have successfully used the GFD method to relate the displacement-gradients at Point $P_0$ to nodal displacements.

It should be noted that for problems involving rupture and fragmentation, the interaction between two adjacent Points needs to be cut off when the internal boundary is cracked, by removing them from each other's support domain (shown in Section 3.3). Thus it is possible that for Point $P_0$, with the development of dense cracks, only one other Point $P_1$ is left inside its support domain, (shown in Figure 4). In that case, the GFD method is incapable to compute displacement-derivatives since the matrix $\mathbf{A}^{\mathrm{T}}\mathbf{W}_1\mathbf{A}$ in Eq. (2.5) will be rank-deficient. Instead, we directly calculate the vector of derivatives $\mathbf{a}_1$ by a simple difference method (Eq. (2.8)) in this case, where $L$ is the distance between the two Points. This simple formula is often used in Moving Particle Semi-implicit Method and other meshless particle methods.

$$\mathbf{a}_1 = \begin{bmatrix} \dfrac{\partial u_1}{\partial x_1} \\ \dfrac{\partial u_1}{\partial x_2} \end{bmatrix}_{P_0} = \begin{bmatrix} \dfrac{u_1^1 - u_1^0}{L} \cdot \dfrac{x_1^1 - x_1^0}{L} \\ \dfrac{u_1^1 - u_1^0}{L} \cdot \dfrac{x_2^1 - x_2^0}{L} \end{bmatrix} = \begin{bmatrix} -\dfrac{x_1^1 - x_1^0}{L^2} & \dfrac{x_1^1 - x_1^0}{L^2} \\ -\dfrac{x_2^1 - x_2^0}{L^2} & \dfrac{x_2^1 - x_2^0}{L^2} \end{bmatrix} \begin{bmatrix} u_1^0 \\ u_1^1 \end{bmatrix} \tag{2.8}$$

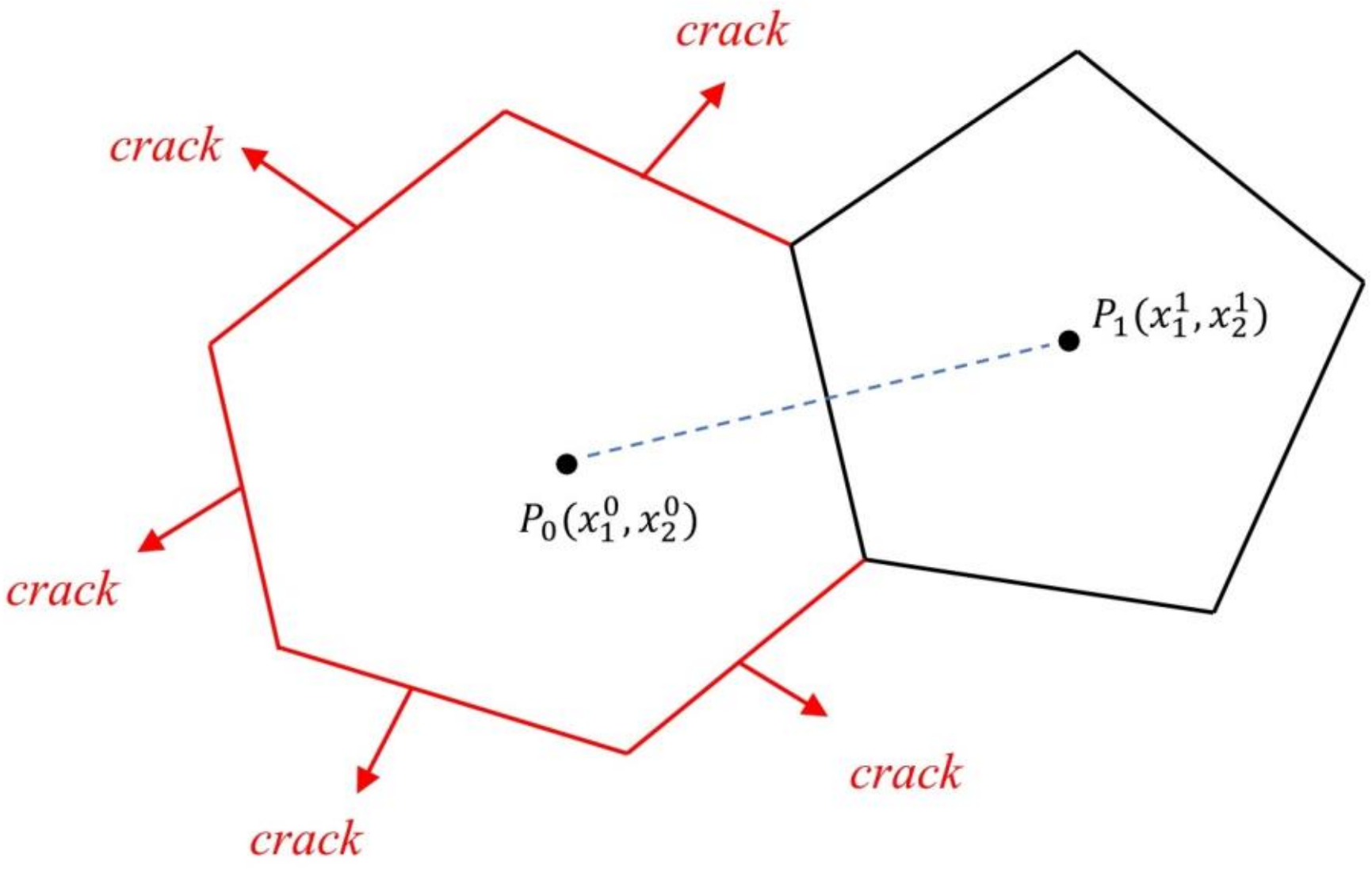

Figure 4. The support of $P_0$ when all the subdomains boundaries but one are cracked

2.2 The Compactly-Supported Radial Basis Function method

The local interpolation of $u_1$ at Point $P_0$ using CSRBFs can be written as

$$u_1(\mathbf{x}) = \sum_{i=0}^{m} R_i(\mathbf{x})\alpha_i + \sum_{j=0}^{q} S_j(\mathbf{x})\beta_j \tag{2.9}$$

where $R_i$ is the compactly-supported radial basis function, and $S_j$ is the polynomial basis function. $m$+1 is the number of all Points located in the support domain of Point $P_0$ and also the number of CSRBFs, and $q$+1 is the number of polynomial basis functions. $\alpha_i$ and $\beta_j$ are constant coefficients yet to be solved. Here, the definition of the support for each Point in the CSRBF method is the same as that of the GFD method.

In a matrix-vector form, Eq. (2.9) can be rewritten as

$$u_1 = \mathbf{R}^{\mathrm{T}}\boldsymbol{\alpha} + \mathbf{S}^{\mathrm{T}}\boldsymbol{\beta}$$

where

$$\mathbf{R}^{\mathrm{T}} = \begin{bmatrix} R_0(\mathbf{x}) & R_1(\mathbf{x}) & \dots & R_m(\mathbf{x}) \end{bmatrix}$$

$$\mathbf{S}^{\mathrm{T}} = \begin{bmatrix} S_0(\mathbf{x}) & S_1(\mathbf{x}) & \dots & S_q(\mathbf{x}) \end{bmatrix}$$

$$\boldsymbol{\alpha}^{\mathrm{T}} = \begin{bmatrix} \alpha_0 & \alpha_1 & \dots & \alpha_m \end{bmatrix}$$

$$\boldsymbol{\beta}^{\mathrm{T}} = \begin{bmatrix} \beta_0 & \beta_1 & \dots & \beta_q \end{bmatrix}$$

Considering first-order polynomial basis functions, $\mathbf{S}^{\mathrm{T}} = \begin{bmatrix} 1 & x_1 & x_2 \end{bmatrix}$ and for second-order polynomial basis functions, $\mathbf{S}^{\mathrm{T}} = \begin{bmatrix} 1 & x_1 & x_2 & x_1^2 & x_2^2 & x_1x_2 \end{bmatrix}$. In this paper, the following Compactly-Supported Radial Basis Function is used:

$$R_i(\mathbf{x}_j) = \begin{cases} \left(1 - \dfrac{d_{ij}}{r}\right)^3 \left(1 + 3\dfrac{d_{ij}}{r}\right) & d_{ij} \le r \\ 0 \quad d_{ij} > r \end{cases}$$

where $d_{ij}$ is the distance between Point $P_i$ and Point $P_j$, and $r$ is the radius of the support domain. Of course, other CSRBFs may also be employed.

In order to determine the coefficient vectors $\boldsymbol{\alpha}$ and $\boldsymbol{\beta}$, Eq. (2.9) is prescribed to be satisfied at all the Points inside the support domain of Point $P_0$, then we have

$$\mathbf{R}_m\boldsymbol{\alpha}+\mathbf{S}_q\boldsymbol{\beta}=\mathbf{u}_1 \tag{2.10}$$

where

$$\mathbf{R}_m=\begin{bmatrix} R_0(\mathbf{x}_0) & R_1(\mathbf{x}_0) & \dots & R_m(\mathbf{x}_0) \\ R_0(\mathbf{x}_1) & R_1(\mathbf{x}_1) & \dots & R_m(\mathbf{x}_1) \\ \dots & \dots & \dots & \dots \\ R_0(\mathbf{x}_m) & R_1(\mathbf{x}_m) & \dots & R_m(\mathbf{x}_m) \end{bmatrix}$$

$$\mathbf{S}_q=\begin{bmatrix} S_0(\mathbf{x}_0) & S_1(\mathbf{x}_0) & \dots & S_q(\mathbf{x}_0) \\ S_0(\mathbf{x}_1) & S_1(\mathbf{x}_1) & \dots & S_q(\mathbf{x}_1) \\ \dots & \dots & \dots & \dots \\ S_0(\mathbf{x}_m) & S_1(\mathbf{x}_m) & \dots & S_q(\mathbf{x}_m) \end{bmatrix}$$

$$(\mathbf{u}_1)^{\mathrm{T}}=\begin{bmatrix} u_1^0 & u_1^1 & \dots & u_1^m \end{bmatrix}$$

Besides, another $q+1$ constraint conditions are employed to solve unknown coefficients [10]:

$$\mathbf{S}_q^{\mathrm{T}}\boldsymbol{\alpha}=0 \tag{2.11}$$

Combining Eq. (2.10) and Eq. (2.11), the following equation in a matrix form is obtained

$$\mathbf{G}\boldsymbol{\alpha}_e=\hat{\mathbf{u}}_1 \tag{2.12}$$

where

$$\mathbf{G}=\begin{bmatrix} \mathbf{R}_m & \mathbf{S}_q \\ \mathbf{S}_q^{\mathrm{T}} & 0 \end{bmatrix} \qquad \boldsymbol{\alpha}_e=\begin{bmatrix} \boldsymbol{\alpha} \\ \boldsymbol{\beta} \end{bmatrix}$$

$$(\hat{\mathbf{u}}_1)^{\mathrm{T}}=\begin{bmatrix} u_1^0 & u_1^1 & \dots & u_1^m & 0 & 0 & \dots & 0 \end{bmatrix}_{1\times(m+q+2)}$$

Therefore, the unknown coefficient vectors $\boldsymbol{\alpha}$ and $\boldsymbol{\beta}$ are solved as

$$\boldsymbol{\alpha}_e=\mathbf{G}^{-1}\hat{\mathbf{u}}_1 \tag{2.13}$$

Substituting Eq. (2.13) into Eq. (2.9), we obtain the interpolation of $u_1(\mathbf{x})$ at Point $P_0$ with CSRBFs and polynomial basis functions.

$$u_1 = \begin{bmatrix} \mathbf{R}^{\mathrm{T}} & \mathbf{S}^{\mathrm{T}} \end{bmatrix} \boldsymbol{\alpha}_e = \begin{bmatrix} \mathbf{R}^{\mathrm{T}} & \mathbf{S}^{\mathrm{T}} \end{bmatrix} \mathbf{G}^{-1} \hat{\mathbf{u}}_1 = \begin{bmatrix} \phi_0 & \cdots & \phi_m \end{bmatrix} \mathbf{u}_1 = \boldsymbol{\Phi}^{\mathrm{T}} \mathbf{u}_1 \tag{2.14}$$

Then, the gradients of $u_1(\mathbf{x})$ can be derived

$$\begin{aligned}
\frac{\partial u_1}{\partial x_1} &= \begin{bmatrix} \frac{\partial \mathbf{R}^{\mathrm{T}}}{\partial x_1} & \frac{\partial \mathbf{S}^{\mathrm{T}}}{\partial x_1} \end{bmatrix} \mathbf{G}^{-1} \hat{\mathbf{u}}_1 = \begin{bmatrix} \phi_{0,1} & \cdots & \phi_{m,1} \end{bmatrix} \mathbf{u}_1 \\
\frac{\partial u_1}{\partial x_2} &= \begin{bmatrix} \frac{\partial \mathbf{R}^{\mathrm{T}}}{\partial x_2} & \frac{\partial \mathbf{S}^{\mathrm{T}}}{\partial x_2} \end{bmatrix} \mathbf{G}^{-1} \hat{\mathbf{u}}_1 = \begin{bmatrix} \phi_{0,2} & \cdots & \phi_{m,2} \end{bmatrix} \mathbf{u}_1 \\
\frac{\partial^2 u_1}{\partial x_1^{\ 2}} &= \begin{bmatrix} \frac{\partial^2 \mathbf{R}^{\mathrm{T}}}{\partial x_1^{\ 2}} & \frac{\partial^2 \mathbf{S}^{\mathrm{T}}}{\partial x_1^{\ 2}} \end{bmatrix} \mathbf{G}^{-1} \hat{\mathbf{u}}_1 = \begin{bmatrix} \phi_{0,11} & \cdots & \phi_{m,11} \end{bmatrix} \mathbf{u}_1 \\
\frac{\partial^2 u_1}{\partial x_2^{\ 2}} &= \begin{bmatrix} \frac{\partial^2 \mathbf{R}^{\mathrm{T}}}{\partial x_2^{\ 2}} & \frac{\partial^2 \mathbf{S}^{\mathrm{T}}}{\partial x_2^{\ 2}} \end{bmatrix} \mathbf{G}^{-1} \hat{\mathbf{u}}_1 = \begin{bmatrix} \phi_{0,22} & \cdots & \phi_{m,22} \end{bmatrix} \mathbf{u}_1 \\
\frac{\partial^2 u_1}{\partial x_1 \partial x_2} &= \begin{bmatrix} \frac{\partial^2 \mathbf{R}^{\mathrm{T}}}{\partial x_1 \partial x_2} & \frac{\partial^2 \mathbf{S}^{\mathrm{T}}}{\partial x_1 \partial x_2} \end{bmatrix} \mathbf{G}^{-1} \hat{\mathbf{u}}_1 = \begin{bmatrix} \phi_{0,12} & \cdots & \phi_{m,12} \end{bmatrix} \mathbf{u}_1
\end{aligned} \tag{2.15}$$

Since the gradients of displacement $u_2(\mathbf{x})$ are calculated in the same way, we can summarize the results as below,

$$\mathbf{a} = \mathbf{C}\mathbf{u}_E \tag{2.16}$$

where for linear trial functions,

$$\mathbf{C} = \begin{bmatrix} \phi_{0,1} & 0 & \cdots & \phi_{m,1} & 0 \\ \phi_{0,2} & 0 & \cdots & \phi_{m,2} & 0 \\ 0 & \phi_{0,1} & \cdots & 0 & \phi_{m,1} \\ 0 & \phi_{0,2} & \cdots & 0 & \phi_{m,2} \end{bmatrix}$$

and for quadratic trial functions,

$$\mathbf{C} = \begin{bmatrix} \phi_{0,1} & 0 & \cdots & \phi_{m,1} & 0 \\ \phi_{0,2} & 0 & \cdots & \phi_{m,2} & 0 \\ \phi_{0,11} & 0 & \cdots & \phi_{m,11} & 0 \\ \phi_{0,22} & 0 & \cdots & \phi_{m,22} & 0 \\ \phi_{0,12} & 0 & \cdots & \phi_{m,12} & 0 \\ 0 & \phi_{0,1} & \cdots & 0 & \phi_{m,1} \\ 0 & \phi_{0,2} & \cdots & 0 & \phi_{m,2} \\ 0 & \phi_{0,11} & \cdots & 0 & \phi_{m,11} \\ 0 & \phi_{0,22} & \cdots & 0 & \phi_{m,22} \\ 0 & \phi_{0,12} & \cdots & 0 & \phi_{m,12} \end{bmatrix}$$

$$\mathbf{u}_E = \begin{bmatrix} u_1^0 & u_2^0 & u_1^1 & u_2^1 & \dots & u_1^m & u_2^m \end{bmatrix}^{\mathrm{T}}$$

As we mentioned in Section 2.1, for a Point, if only one other Point is within its support domain because of the development of dense cracks, we employ the simple difference method to calculate derivatives.

It should be noted that in addition to the GFD and CSRBF methods, plenty of other approaches are also available to calculate the gradients of displacement such as Moving Least Square, Smoothed Particle Hydrodynamics, Moving Particle Semi-implicit methods, …etc. However, the difference about the development of trial functions between the present FPM and other meshless methods is that: in the current FPM, displacement-gradients are assumed to be constant over each subdomain, and are thus only calculated at each Point.

## 2.3 The derivation of trial and test functions

Finally, by substituting Eq. (2.7) or Eq. (2.16) into Eq. (2.3), the relation between $\mathbf{u}^h$ and $\mathbf{u}_E$ is obtained in Eq. (2.17), where the matrix $\mathbf{N}$ is called the shape function of $\mathbf{u}^h$ in $E_0$.

$$\mathbf{u}^h = \mathbf{N}\mathbf{u}_E \tag{2.17}$$

$$\begin{aligned}\mathbf{N} &= \begin{bmatrix} \mathbf{h}(x_1, x_2) & 0 \\ 0 & \mathbf{h}(x_1, x_2) \end{bmatrix} \mathbf{C} + \mathbf{I}_3 \\ &= \begin{bmatrix} N_0 & 0 & N_1 & 0 & \dots & N_m & 0 \\ 0 & N_0 & 0 & N_1 & \dots & 0 & N_m \end{bmatrix}_{2\times(2m+2)}\end{aligned}$$

$$\text{where } \mathbf{I}_3 = \begin{bmatrix} 1 & 0 & 0 & \dots & 0 \\ 0 & 1 & 0 & \dots & 0 \end{bmatrix}_{2\times(2m+2)}$$

According to the Eq. (2.1), the corresponding strain $\boldsymbol{\varepsilon}^h$ and stress $\boldsymbol{\sigma}^h$ in terms of $\mathbf{u}_E$ are given in Eq. (2.18),

$$
\boldsymbol{\varepsilon}^h = \begin{bmatrix} \varepsilon_{11}^h \\ \varepsilon_{22}^h \\ 2\varepsilon_{12}^h \end{bmatrix} = \begin{bmatrix} \frac{\partial}{\partial x_1} & 0 \\ 0 & \frac{\partial}{\partial x_2} \\ \frac{\partial}{\partial x_2} & \frac{\partial}{\partial x_1} \end{bmatrix} \mathbf{u}^h = \mathbf{B}\mathbf{u}_E
$$

$$
\boldsymbol{\sigma}^h = \begin{bmatrix} \sigma_{11}^h \\ \sigma_{22}^h \\ \sigma_{12}^h \end{bmatrix} = \mathbf{D}\boldsymbol{\varepsilon}^h = \mathbf{D}\mathbf{B}\mathbf{u}_E \tag{2.18}
$$

where **D** is the stress-strain matrix. In this paper, we consider the material to be isotropic for simplicity.

$$
\mathbf{D} = \frac{\bar{E}}{1-\bar{v}^2} \begin{bmatrix} 1 & \bar{v} & 0 \\ \bar{v} & 1 & 0 \\ 0 & 0 & (1-\bar{v})/2 \end{bmatrix}
$$

$$
\text{where } \bar{E} = \begin{cases} E & (\text{for plane stress}) \\ \dfrac{E}{1-v^2} & (\text{for plane strain}) \end{cases}
$$

$$
\bar{v} = \begin{cases} v & (\text{for plane stress}) \\ \dfrac{v}{1-v} & (\text{for plane strain}) \end{cases}
$$

Following the same procedure, we can derive $\mathbf{u}^h$ in each subdomain $E_i \in \Omega$. Eventually, the displacement vector $\mathbf{u}^h$ in the entire domain can be obtained. The corresponding test function $\mathbf{v}$ is prescribed to possess the same shape as the trial function in each subdomain, in the present FPM based on the Galerkin weak form.

Reviewing the process of constructing trial and test functions, we can see that no continuity requirements exist on the internal boundary between two contiguous neighboring subdomains. In other words, these two contiguous subdomains have their own trial and test function values on their common internal boundary. Therefore, only simple, local, polynomial, Point-based and piecewise-continuous trial and test functions are employed in the problem domain $\Omega$.

To illustrate the discontinuity of trial functions and shape functions, a 2D example is shown here. We assume that 25 Points are scattered irregularly in a $1\times1$ square. The nodal displacement-gradients are related to nodal displacements using the generalized finite difference method. For simplicity, trial functions are linear in this example. The graphical representation of all the shape functions about Point 13 (the subscripts in Eq. (2.16) equal 13) is given in Figure 5. The corresponding trial function of $u_1$ simulating the exponential function $e^{-10(x_1-0.5)^2-10(x_2-0.5)^2}$ is shown in Figure 6.

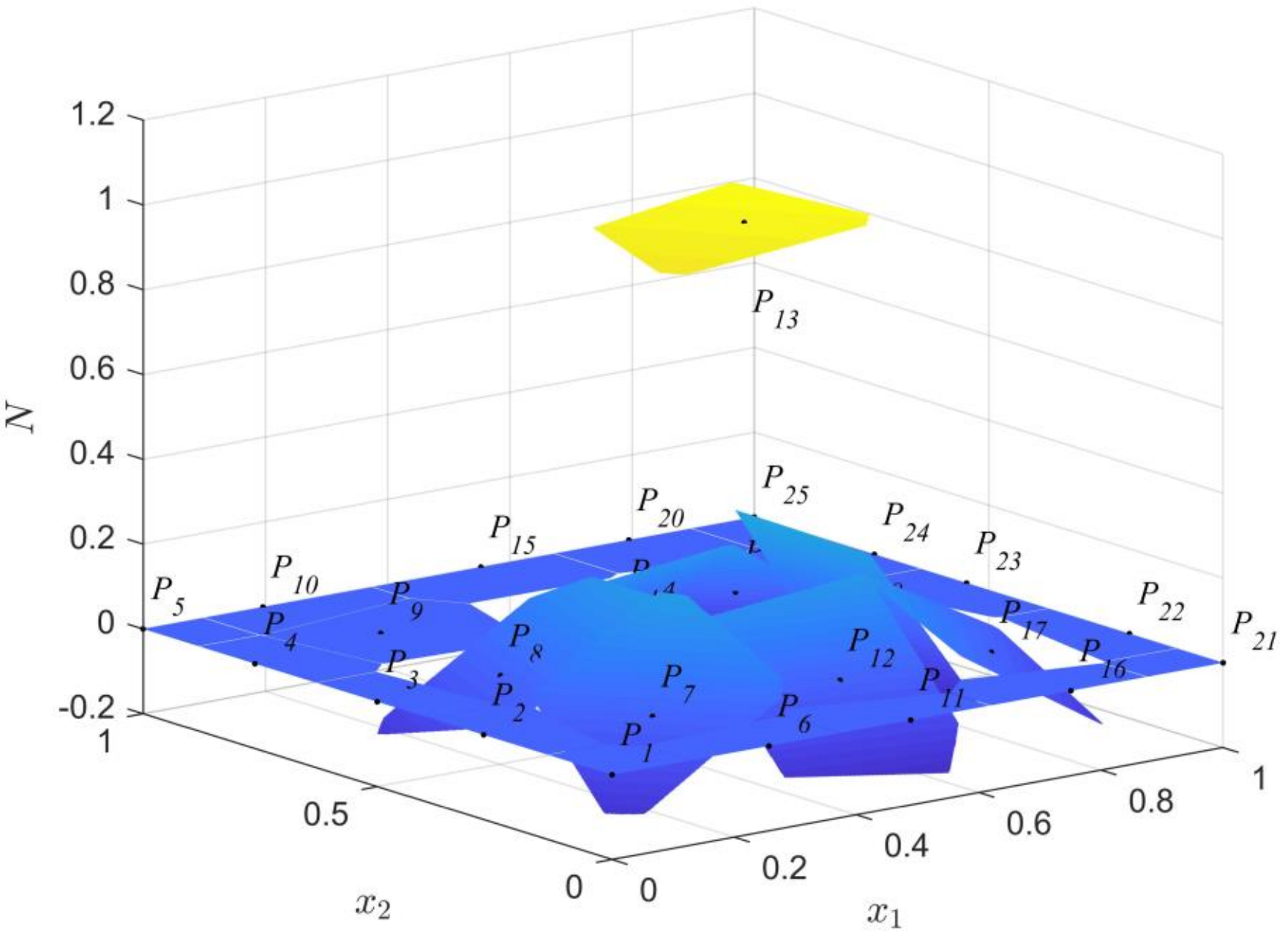

Figure 5. The shape functions about Point 13

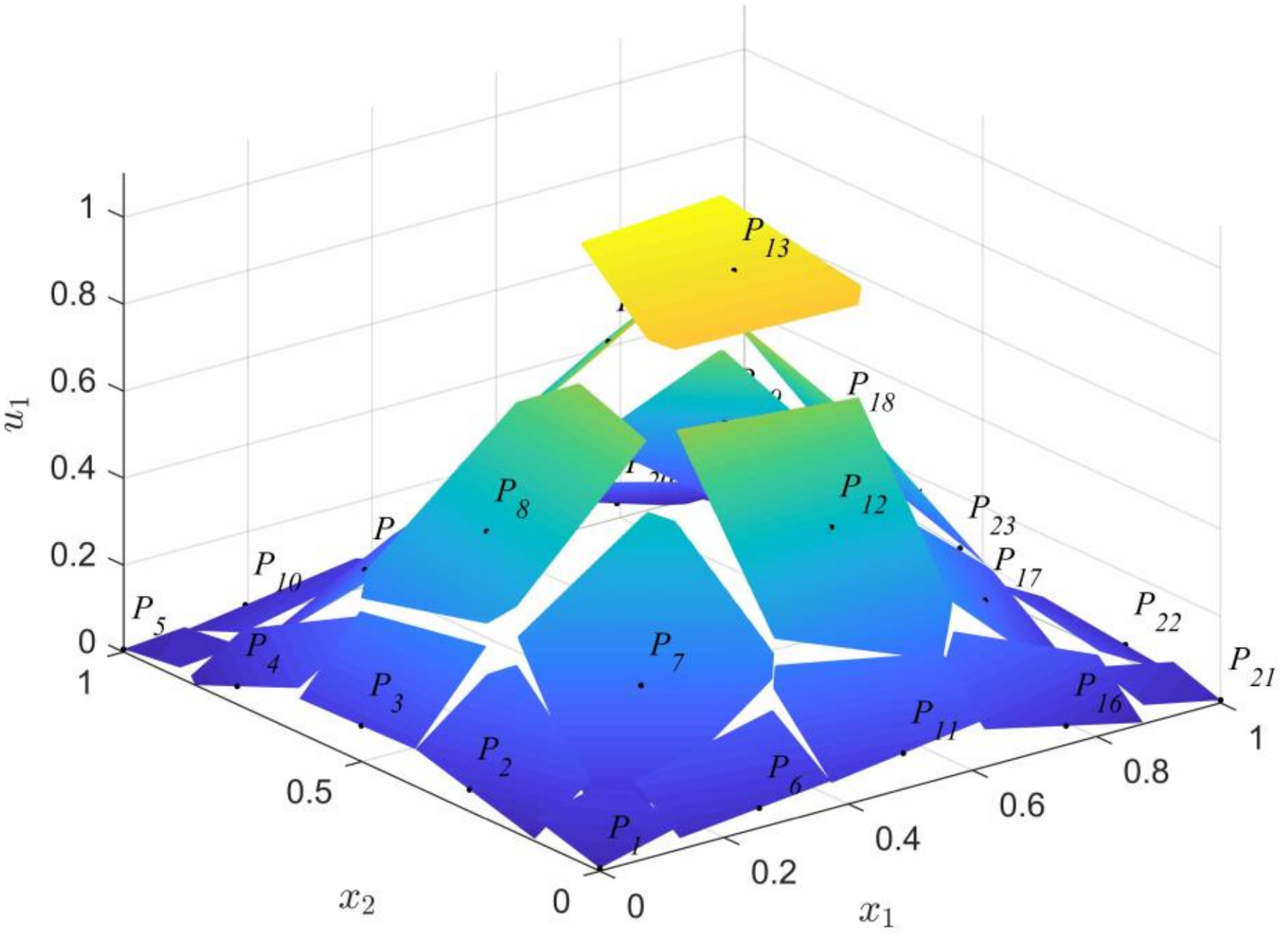

Figure 6. The trial function simulating an exponential function

Unfortunately, because of this discontinuity of trial and test functions, if we directly use the trial and test functions in the traditional Galerkin weak form which is widely used in the FEM, EFG and other numerical methods, the solution will be inconsistent, inaccurate and cannot pass the patch tests [8]. In order to solve this inconsistency problem, Numerical Flux Corrections are introduced to the FPM.

# 3. Numerical Flux Corrections and Algorithmic Implementation

## 3.1 Interior Penalty (IP) Numerical Flux Corrections

Numerical Fluxes, frequently used in Discontinuous Galerkin FEM Methods, are employed in FPM to resolve the inconsistency caused by the discontinuity of trial and test functions. A variety of Numerical Fluxes have been developed in literature. In this study, the Interior Penalty Numerical Fluxes which can lead to consistent and stable results with symmetric global stiffness matrices, are used.

The governing equations of linear elasticity in 2D have been shown in Eq. (2.1). We multiply the second equation by the test function **v** and then integrate it on a generic subdomain $E$ by parts,

$$\int_E \sigma_{ij}(\mathbf{u}) v_{i,j} d\Omega - \int_{\partial E} \sigma_{ij}(\mathbf{u}) n_j v_i d\Gamma = \int_E f_i v_i d\Omega \tag{3.1}$$

where $\partial E$ is the boundary of the subdomain $E$, $\mathbf{n}$ is the unit vector outward to $\partial E$.

For every subdomain $E \in \Omega$, Eq. (3.1) should be satisfied. Therefore, we sum Eq. (3.1) over all subdomains.

$$\sum_{E\in\Omega} \int_E \sigma_{ij}(\mathbf{u}) v_{i,j} d\Omega - \sum_{E\in\Omega} \int_{\partial E} \sigma_{ij}(\mathbf{u}) n_j v_i d\Gamma = \sum_{E\in\Omega} \int_E f_i v_i d\Omega \tag{3.2}$$

Considering the symmetry of stress tensor ($\sigma_{ij} = \sigma_{ji}$), we can transform the first term of Eq. (3.2) into the following form:

$$\sum_{E\in\Omega} \int_E \sigma_{ij}(\mathbf{u}) v_{i,j} d\Omega = \sum_{E\in\Omega} \int_E \sigma_{ij}(\mathbf{u}) \varepsilon_{ij}(\mathbf{v}) d\Omega \tag{3.3}$$

Let $\Gamma$ denote the set of all external and internal boundaries and $\Gamma_h = \Gamma - \Gamma_t - \Gamma_u$ stands for the set of all internal boundaries. For convenience, we rewrite the second term in Eq. (3.2) with the jump operator [] and average operator {}.

$$\begin{aligned} \sum_{E\in\Omega} \int_{\partial E} \sigma_{ij}(\mathbf{u}) n_j v_i d\Gamma &= \sum_{e\in\Gamma_h} \int_e \left( \left[\sigma_{ij}(\mathbf{u}) n_j^e\right]\{v_i\} + \{\sigma_{ij}(\mathbf{u}) n_j^e\}[v_i] \right) d\Gamma \\ &+ \sum_{e\in\Gamma_u} \int_e \{\sigma_{ij}(\mathbf{u}) n_j^e\}[v_i] d\Gamma + \sum_{e\in\Gamma_t} \int_e \{\sigma_{ij}(\mathbf{u}) n_j^e\}[v_i] d\Gamma \end{aligned} \tag{3.4}$$

When $e \in \Gamma_h$ (assuming $e$ is shared by subdomains $E_1$ and $E_2$), $n_j^e$ is a unit vector normal to $e$ and points from $E_1$ to $E_2$. The average {} and jump [] operator for any quantity $w$, at an internal boundary are defined as Eq. (3.5).

$$[w] = w\big|_e^{E_1} - w\big|_e^{E_2}, \qquad \{w\} = \frac{1}{2}\left( w\big|_e^{E_1} + w\big|_e^{E_2} \right) \tag{3.5}$$

When $e \in \partial\Omega$, $n_j^e$ is pointing outward to $\partial\Omega$ and the average {} and jump [] operator are defined as $[w] = w\big|_e, \{w\} = w\big|_e$. It should be noted that for two neighboring subdomains, no matter which one is chosen as $E_1$, Eq (3.4) remains unchanged.

If $\sigma_{ij}(\mathbf{u})$ is the exact solution in an intact domain, there should be no jump on internal boundaries. In other words, $\left[\sigma_{ij}(\mathbf{u}) n_j^e\right] = 0$ should be satisfied in a strong or a weak-sense, i.e., tractions are reciprocated at internal intact boundaries. Besides, with the traction boundary condition $\sigma_{ij}(\mathbf{u})n_j = \bar{t}_i$, Eq. (3.4) can be rewritten as below.

$$\sum_{E\in\Omega}\int_{\partial E}\sigma_{ij}(\mathbf{u})n_j v_i d\Gamma = \sum_{e\in\Gamma_h\cup\Gamma_u}\int_e \{\sigma_{ij}(\mathbf{u})n_j^e\}[v_i]d\Gamma + \sum_{e\in\Gamma_t}\int_e \bar{t}_i v_i d\Gamma \tag{3.6}$$

Eventually, we substitute Eq. (3.3), Eq. (3.6) into Eq. (3.2), and add two boundary integrals $-\sum_{e\in\Gamma_h\cup\Gamma_u}\int_e \{\sigma_{ij}(\mathbf{v})n_j^e\}[u_i]d\Gamma$ and $\sum_{e\in\Gamma_h\cup\Gamma_u}\frac{\eta}{h_e}\int_e [u_i][v_i]d\Gamma$. Note these two terms should vanish for the exact solution, as the displacement jump should be zero at internal intact boundaries. Then we obtain the weak form with Interior Penalty Numerical Flux Corrections for linear elasticity,

$$\begin{aligned}
&\sum_{E\in\Omega}\int_E \sigma_{ij}(\mathbf{u})\varepsilon_{ij}(\mathbf{v})d\Omega - \sum_{e\in\Gamma_h\cup\Gamma_u}\int_e \{\sigma_{ij}(\mathbf{u})n_j^e\}[v_i]d\Gamma \\
&-\sum_{e\in\Gamma_h\cup\Gamma_u}\int_e \{\sigma_{ij}(\mathbf{v})n_j^e\}[u_i]d\Gamma + \sum_{e\in\Gamma_h\cup\Gamma_u}\frac{\eta}{h_e}\int_e [u_i][v_i]d\Gamma \\
&= \sum_{E\in\Omega}\int_E f_i v_i d\Omega + \sum_{e\in\Gamma_t}\int_e \bar{t}_i v_i d\Gamma - \sum_{e\in\Gamma_u}\int_e \sigma_{ij}(\mathbf{v})n_j^e \bar{u}_i d\Gamma + \sum_{e\in\Gamma_u}\frac{\eta}{h_e}\int_e \bar{u}_i v_i d\Gamma
\end{aligned} \tag{3.7}$$

where $h_e$ is an edge-dependent parameter and is taken to be equal to the length of the boundary segment in this paper; $\eta$ is a positive scalar independent of the edge size. It should be noted that with IP Numerical Flux Corrections, the method is only stable when the penalty parameter $\eta$ is large enough [11]. A discussion about the effect of the penalty parameter is given in Section 4 and more information can be found in [12].

We can find that in Eq. (3.7), displacement boundary conditions are imposed weakly. If we impose $u_i = \bar{u}_i$ strongly at the boundary points, Eq. (3.7) can be simplified as follows.

$$\begin{aligned}
&\sum_{E\in\Omega}\int_E \sigma_{ij}(\mathbf{u})\varepsilon_{ij}(\mathbf{v})d\Omega - \sum_{e\in\Gamma_h}\int_e \{\sigma_{ij}(\mathbf{u})n_j^e\}[v_i]d\Gamma \\
&-\sum_{e\in\Gamma_h}\int_e \{\sigma_{ij}(\mathbf{v})n_j^e\}[u_i]d\Gamma + \sum_{e\in\Gamma_h}\frac{\eta}{h_e}\int_e [u_i][v_i]d\Gamma \\
&= \sum_{E\in\Omega}\int_E f_i v_i d\Omega + \sum_{e\in\Gamma_t}\int_e \bar{t}_i v_i d\Gamma
\end{aligned} \tag{3.8}$$

For brevity, we can rewrite Eq. (3.8) in a matrix-vector form,

$$
\begin{aligned}
&\sum_{E\in\Omega}\int_E \boldsymbol{\varepsilon}_v^{\mathrm{T}}\boldsymbol{\sigma}_u d\Omega-\sum_{e\in\Gamma_h}\int_e [\mathbf{v}]^{\mathrm{T}}\{\mathbf{n}_e\boldsymbol{\sigma}_u\}d\Gamma \\
&-\sum_{e\in\Gamma_h}\int_e \{\mathbf{n}_e\boldsymbol{\sigma}_v\}^{\mathrm{T}}[\mathbf{u}]d\Gamma+\sum_{e\in\Gamma_h}\frac{\eta}{h_e}\int_e [\mathbf{v}]^{\mathrm{T}}[\mathbf{u}]d\Gamma \\
&=\sum_{E\in\Omega}\int_E \mathbf{v}^{\mathrm{T}}\mathbf{f}d\Omega+\sum_{e\in\Gamma_t}\int_e \mathbf{v}^{\mathrm{T}}\bar{\mathbf{t}}d\Gamma
\end{aligned}
\tag{3.9}
$$

where

$$
\boldsymbol{\sigma}=\begin{bmatrix}\sigma_{11}\\ \sigma_{22}\\ \sigma_{12}\end{bmatrix}\quad
\boldsymbol{\varepsilon}=\begin{bmatrix}\varepsilon_{11}\\ \varepsilon_{22}\\ 2\varepsilon_{12}\end{bmatrix}\quad
\mathbf{u}=\begin{bmatrix}u_1\\ u_2\end{bmatrix}\quad
\mathbf{v}=\begin{bmatrix}v_1\\ v_2\end{bmatrix}
$$

$$
\mathbf{n}_e=\begin{bmatrix}n_1^e & 0 & n_2^e\\ 0 & n_2^e & n_1^e\end{bmatrix}\quad
\mathbf{f}=\begin{bmatrix}f_1\\ f_2\end{bmatrix}\quad
\bar{\mathbf{t}}=\begin{bmatrix}\bar{t}_1\\ \bar{t}_2\end{bmatrix}
$$

Compared with the traditional Galerkin weak form [1], the Eq. (3.9) involves 3 extra boundary integrals on the left side, while the others stay identical. These additional boundary integrals are the contributions of the Interior Penalty Numerical Flux Corrections.

## 3.2 Point and Boundary Stiffness Matrices

This section will concentrate on the algorithmic implementation of the FPM. In Section 2, we have obtained the shape function **N** for $\mathbf{u}^h$ and **v**, **B** for $\boldsymbol{\varepsilon}$, **DB** for $\boldsymbol{\sigma}$. By substituting them into the first term of Eq. (3.9), we derive the Point Stiffness Matrix $\mathbf{K}_E$, which is defined as the contribution of each Point to the global stiffness matrix.

$$
\mathbf{K}_E=\int_E \mathbf{B}^{\mathrm{T}}\mathbf{D}\mathbf{B}d\Omega,\quad \forall E\in\Omega \tag{3.10}
$$

For the boundary integrals, the corresponding boundary stiffness matrix $\mathbf{K}_h$ is defined as below. The subscripts 1 and 2 denote which subdomain these shape functions belong to.

$$
\begin{aligned}
\mathbf{K}_h = & -\frac{1}{2}\int_e \mathbf{N}_1^{\mathrm{T}}\mathbf{n}_e\mathbf{D}\mathbf{B}_1 d\Gamma - \frac{1}{2}\int_e \mathbf{B}_1^{\mathrm{T}}\mathbf{D}\mathbf{n}_e^{\mathrm{T}}\mathbf{N}_1 d\Gamma + \frac{\eta}{h_e}\int_e \mathbf{N}_1^{\mathrm{T}}\mathbf{N}_1 d\Gamma \\
& -\frac{1}{2}\int_e \mathbf{N}_1^{\mathrm{T}}\mathbf{n}_e\mathbf{D}\mathbf{B}_2 d\Gamma + \frac{1}{2}\int_e \mathbf{B}_1^{\mathrm{T}}\mathbf{D}\mathbf{n}_e^{\mathrm{T}}\mathbf{N}_2 d\Gamma - \frac{\eta}{h_e}\int_e \mathbf{N}_1^{\mathrm{T}}\mathbf{N}_2 d\Gamma \\
& +\frac{1}{2}\int_e \mathbf{N}_2^{\mathrm{T}}\mathbf{n}_e\mathbf{D}\mathbf{B}_1 d\Gamma - \frac{1}{2}\int_e \mathbf{B}_2^{\mathrm{T}}\mathbf{D}\mathbf{n}_e^{\mathrm{T}}\mathbf{N}_1 d\Gamma - \frac{\eta}{h_e}\int_e \mathbf{N}_2^{\mathrm{T}}\mathbf{N}_1 d\Gamma \\
& +\frac{1}{2}\int_e \mathbf{N}_2^{\mathrm{T}}\mathbf{n}_e\mathbf{D}\mathbf{B}_2 d\Gamma + \frac{1}{2}\int_e \mathbf{B}_2^{\mathrm{T}}\mathbf{D}\mathbf{n}_e^{\mathrm{T}}\mathbf{N}_2 d\Gamma + \frac{\eta}{h_e}\int_e \mathbf{N}_2^{\mathrm{T}}\mathbf{N}_2 d\Gamma \\
& \forall\ e \in \Gamma_h
\end{aligned}
\tag{3.11}
$$

When linear trial functions are employed for $\mathbf{u}^h$, the shape function **B** is constant and **N** is linear in each subdomain. Therefore, the integral for submatrix $\mathbf{K}_E$ can be calculated just multiplying the integrand by the area of the corresponding subdomain. For integrals on boundaries, the numerical integration and direct analytic computation are both effective. In this paper, 1 Point Gauss Integration method is used for boundary integrals as it is sufficient for linear trial functions which are used.

When higher-order trial functions are employed, in order to calculate the submatrix $\mathbf{K}_E$, Hammer integration will be effective if we divide the polynomial subdomain into several triangles. Specifically, for FPM with second-order trial functions, three Points Hammer integration for each triangle will be accurate enough.

In the FPM, the global stiffness matrix **K** is obtained by assembling all the submatrices $\mathbf{K}_E$ and $\mathbf{K}_h$. This assembling process is the same as what is done in the FEM. Eventually, the FPM will lead to a linear system of equations with a sparse, symmetric and positive definitive global stiffness matrix:

$$
\mathbf{Kq} = \mathbf{Q} \tag{3.12}
$$

where **K** is the global stiffness matrix, **q** is the vector with nodal DoFs, **Q** is the load vector and computed through Eq. (3.13).

$$
\mathbf{Q} = \sum_{E\in\Omega}\int_E \mathbf{N}^{\mathrm{T}}\mathbf{f}\, d\Omega + \sum_{e\in\Gamma_t}\int_e \mathbf{N}^{\mathrm{T}}\bar{\mathbf{t}}\, d\Gamma \tag{3.13}
$$

The entire algorithm is thus summarized as follows.

1. Loop over all Points in the problem domain $\Omega$

a. Compute shape functions

(**C** from Eq. (2.7) or (2.16), **N** from Eq. (2.17), **B** from Eq. (2.18))

b. Compute the Point Stiffness Matrix $\mathbf{K}_E$ from Eq. (3.10)

2. Loop over all internal boundaries

   a. Compute shape functions

   b. Compute the Boundary Stiffness Matrix $\mathbf{K}_h$ from Eq. (3.11)

3. Assemble all $\mathbf{K}_E$ and $\mathbf{K}_h$ to obtain the global stiffness matrix **K**

4. Compute the load vector **Q** from Eq. (3.13)

5. Solve the assembled Eq (3. 12) to obtain the nodal DoFs vector **q**

6. Postprocess to obtain displacements, strains, and stresses at each Point

### 3.3 Simulations of Crack Initiation & Propagation in FPM

With the simple discontinuous Point-based trial and test functions of FPM, we simulate development of cracks by cut off interactions between two neighboring Points when the internal boundary of them is cracked. As the support of each Point contains only its neighbors, we only need to set the internal boundary into two traction-free boundaries and remove the two Points from each other's support (shown in Figure 7(b)) when their internal boundary is cracked. Specifically, for the Point on one side of the crack, we prescribe that it will not be included in the support of the Point located on the other side of the crack. This is consistent with the visibility criterion used in meshless methods [2]. If the support domain is defined as a circle, according to the visibility criterion, as shown in Figure 7(a), the Point located in the shadow is removed from the center Point's support domain.

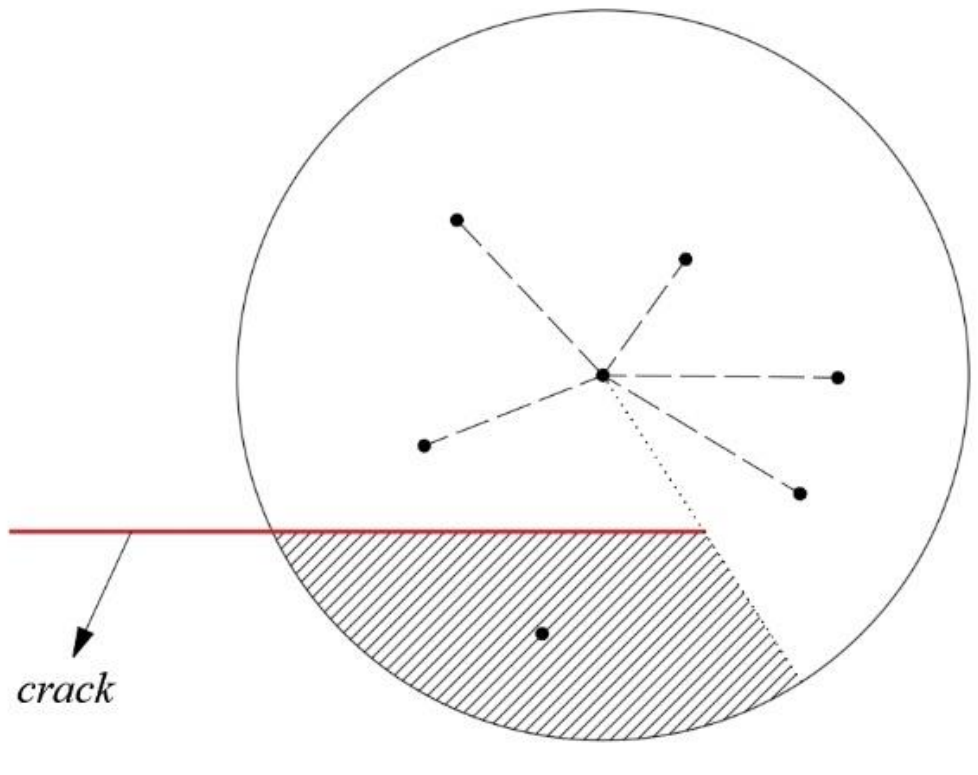

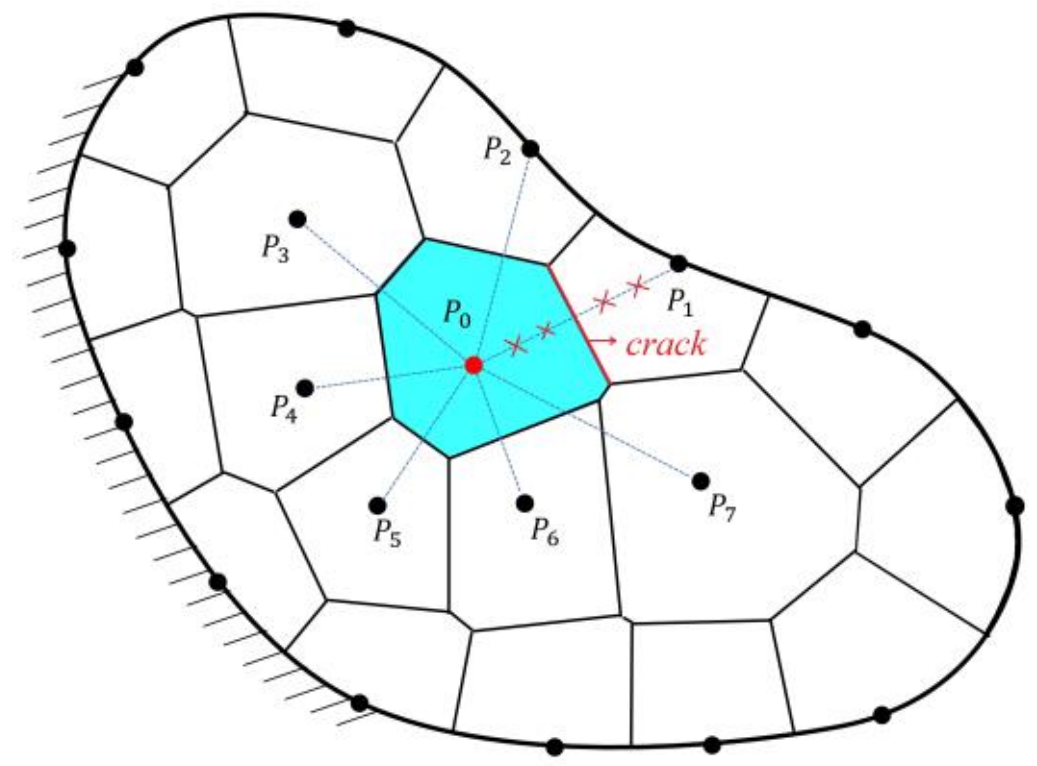

(a) (b)

Figure 7 (a). The visibility criterion when the support domain is circular

(b) The visibility criterion when the support domain contains neighboring Points

Getting back to the Eq. (3.9), in addition to the adjustment of the supports, we only need to delete the terms related to the IP numerical fluxes of an internal boundary, this will lead to changes of $\mathbf{K}_h$ and $\mathbf{K}_E$ for a few points and internal boundaries near the newly developed crack-segment. In this way, the stiffness matrix $\mathbf{K}$ will only need to be adjusted slightly and there is no change for the terms on the right-hand side. Moreover, the number of the DoFs, the dimensions of the global stiffness matrix and the load vector will stay the same. This is much simpler than remeshing, or deleting finite elements, to model the development of cracks. Examples involving crack initiation and propagation will be given in Section 4.

It should be noted that the FPM proposed in this study is natural for simulating phenomena related to cracks, such as damage, fracture and fragmentation. However, the traditional-continuum-physics-based criteria for crack initiation and propagation is a different matter, completely apart from the numerical method of FPM itself. For different problems and for different materials, various criteria for crack initiation and propagation have been developed. Appropriate criteria should be used based on which realistic engineering problem is to be solved. In this study, we only adopt certain criteria for simulation to demonstrate the power of the FPM numerical method itself. Detailed discussions and judgments on various criteria of crack developments are beyond the scope of the current paper.

# 4. Numerical Examples

In this section, a variety of problems are solved with the FPM. In order to estimate the errors of numerical results conveniently, we define two relative errors $r_u$ and $r_E$ with the displacement $L^2$ norm and the energy norm, respectively.

$$r_u = \frac{\left\| \mathbf{u}^h - \mathbf{u}^{\text{exact}} \right\|_{L^2}}{\left\| \mathbf{u}^{\text{exact}} \right\|_{L^2}}$$

$$r_E = \frac{\left\| E^h - E^{\text{exact}} \right\|_E}{\left\| E^{\text{exact}} \right\|_E} \tag{4.1}$$

where $\left\| \mathbf{u} \right\|_{L^2} = \left( \int_\Omega \mathbf{u}^{\mathrm{T}} \mathbf{u} d\Omega \right)^{1/2}$, $\left\| E \right\|_E = \left( \frac{1}{2} \int_\Omega \boldsymbol{\varepsilon}^{\mathrm{T}} \mathbf{D} \boldsymbol{\varepsilon} d\Omega \right)^{1/2}$

If without specific statements, examples in this paper are solved by FPM with linear trial functions.

## 4.1 Patch Test

In this subsection, we design the following patch test in a unit square domain (shown in (c)

Figure 8) to examine the consistency of the FPM. A Plane Stress condition is considered, with the exact displacements and stresses prescribed as below

$$\mathbf{u} = \begin{bmatrix} u_1 \\ u_2 \end{bmatrix} = \begin{bmatrix} x_1 + x_2 \\ x_1 + x_2 \end{bmatrix}$$

$$\boldsymbol{\sigma} = \begin{bmatrix} \sigma_{11} \\ \sigma_{22} \\ \sigma_{12} \end{bmatrix} = \begin{bmatrix} \dfrac{E}{1-v} \\ \dfrac{E}{1-v} \\ \dfrac{E}{1+v} \end{bmatrix} \tag{4.2}$$

$$x_1, x_2 \in (0,1)$$

Traction boundary conditions are prescribed on the edges, according to Eq. (4.2), and three displacement constraints are enforced to prevent rigid-body movements.

Since the solutions are linear for the displacements, when linear trial functions are employed for displacements, the numerical solutions of $\mathbf{u}^h$ and $\boldsymbol{\sigma}^h$ should reproduce those in Eq. (4.2).

The distributions of points in 3 different patterns are given in (c)

Figure 8. In these three cases, no matter whether the points are scattered uniformly or randomly, the present FPM is accurate enough to pass the patch tests (see Table II). This demonstrates the consistency of the present FPM.

Table II. Relative errors for patch tests

| | 9 regular points | 9 random points | 25 random points |
|---|---|---|---|
| $r_u$ (GFD) | $6.23\times10^{-7}$ | $5.66\times10^{-7}$ | $6.80\times10^{-7}$ |
| $r_E$ (GFD) | $1.48\times10^{-7}$ | $2.59\times10^{-7}$ | $2.51\times10^{-7}$ |
| $r_u$ (CSRBF) | $5.91\times10^{-7}$ | $4.36\times10^{-7}$ | $5.44\times10^{-7}$ |
| $r_E$ (CSRBF) | $2.38\times10^{-7}$ | $2.50\times10^{-7}$ | $3.71\times10^{-7}$ |

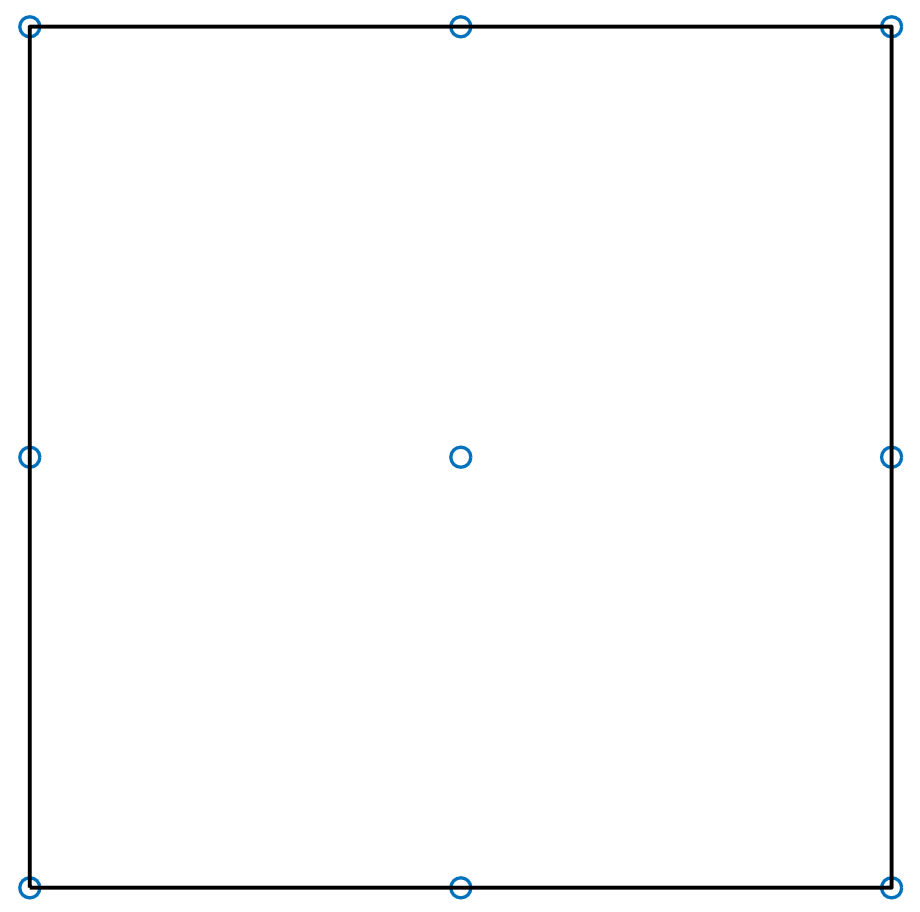

(a)

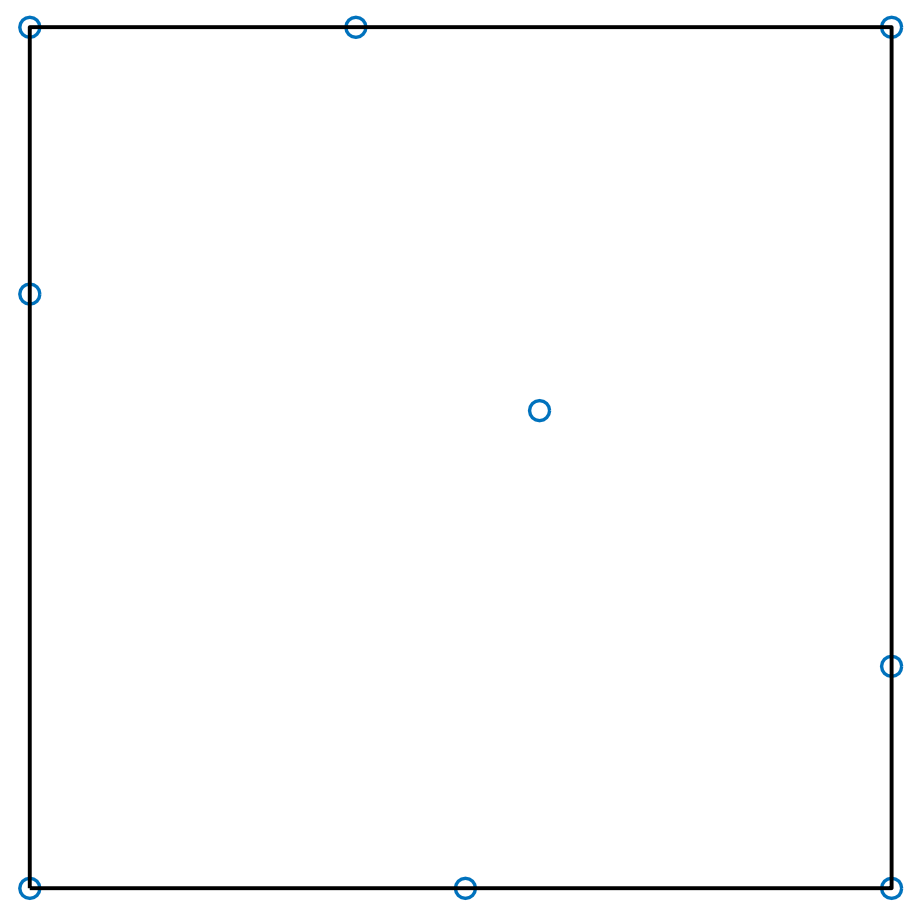

(b)

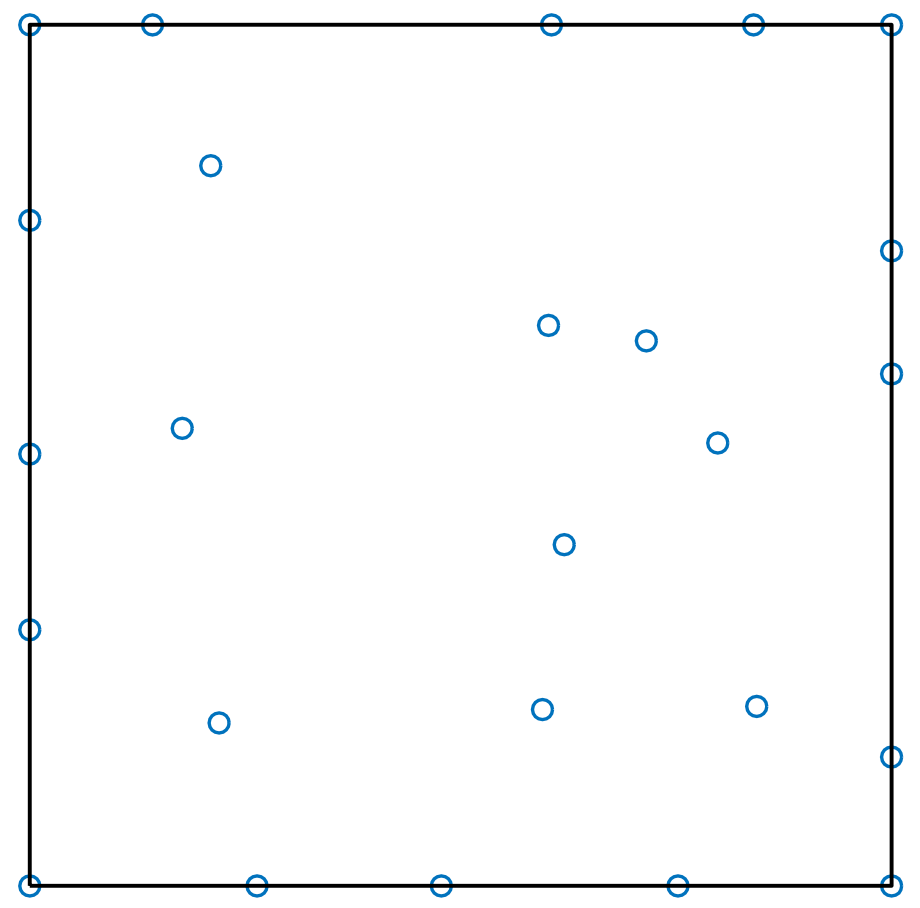

(c)

Figure 8. Three distributions of Points for patch tests

(a) 9 regular points (b)9 irregular points (c) 25 irregular points

## 4.2 A Cantilever Beam

In this section, we employ the FPM to solve a cantilever beam problem with a parabolic-shear traction at one end (shown in Figure 9). The corresponding analytical solutions of displacements and stresses for the Plane Stress case are given in [13].

$$u_1 = -\frac{P}{6EI}\left(x_2 - \frac{H}{2}\right)\left[3x_1\left(2L - x_1\right) + \left(2 + v\right)x_2\left(x_2 - H\right)\right] \tag{4.3}$$

$$u_2 = \frac{P}{6EI}\left[ x_1^{\,2}\left(3L - x_1\right) + 3v\left(L - x_1\right)\left(x_2 - \frac{H}{2}\right)^2 + \frac{4+5v}{4} H^2 x_1 \right]$$

$$\begin{bmatrix} \sigma_{11} \\ \sigma_{22} \\ \sigma_{12} \end{bmatrix} = \begin{bmatrix} -\frac{P}{I}\left(L - x_1\right)\left(x_2 - \frac{H}{2}\right) \\ 0 \\ -\frac{Px_2}{2I}\left(x_2 - H\right) \end{bmatrix}$$

where

$$I = \frac{H^3}{12}$$

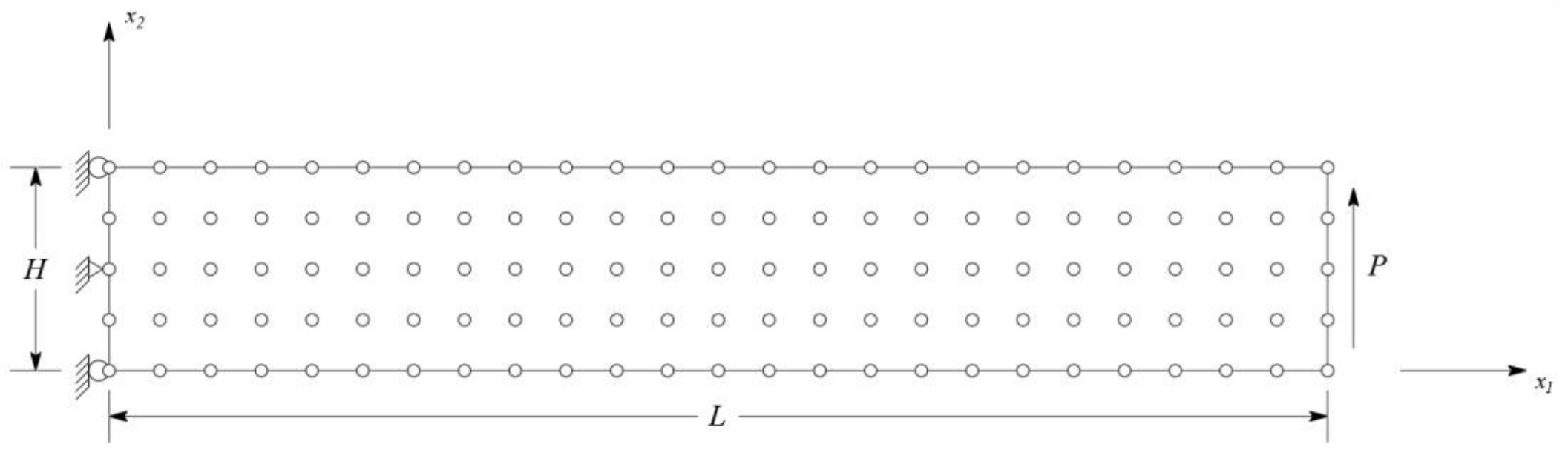

Figure 9. The cantilever beam with a parabolic-shear load

Specifically, we prescribe that $P = 1$, $E = 1\times10^5$, $H = 1$ and $L = 8$. Poisson's ratio $v$ is prescribed as 0.3 and the penalty parameter $\eta = E$ . Displacement boundary conditions are imposed on the left and right edges of the beam, and traction boundary conditions are prescribed on the upper and bottom edges, according to Eq. (4.3). With 891 points distributed in the beam either uniformly or randomly, the comparisons between numerical solutions $u_1$ , $\sigma_{11}$ derived from linear trial functions (GFD or CSRBF) and analytical solutions along the line $x_1 = L/2$ are given in Figure 10(a) and (b), respectively.

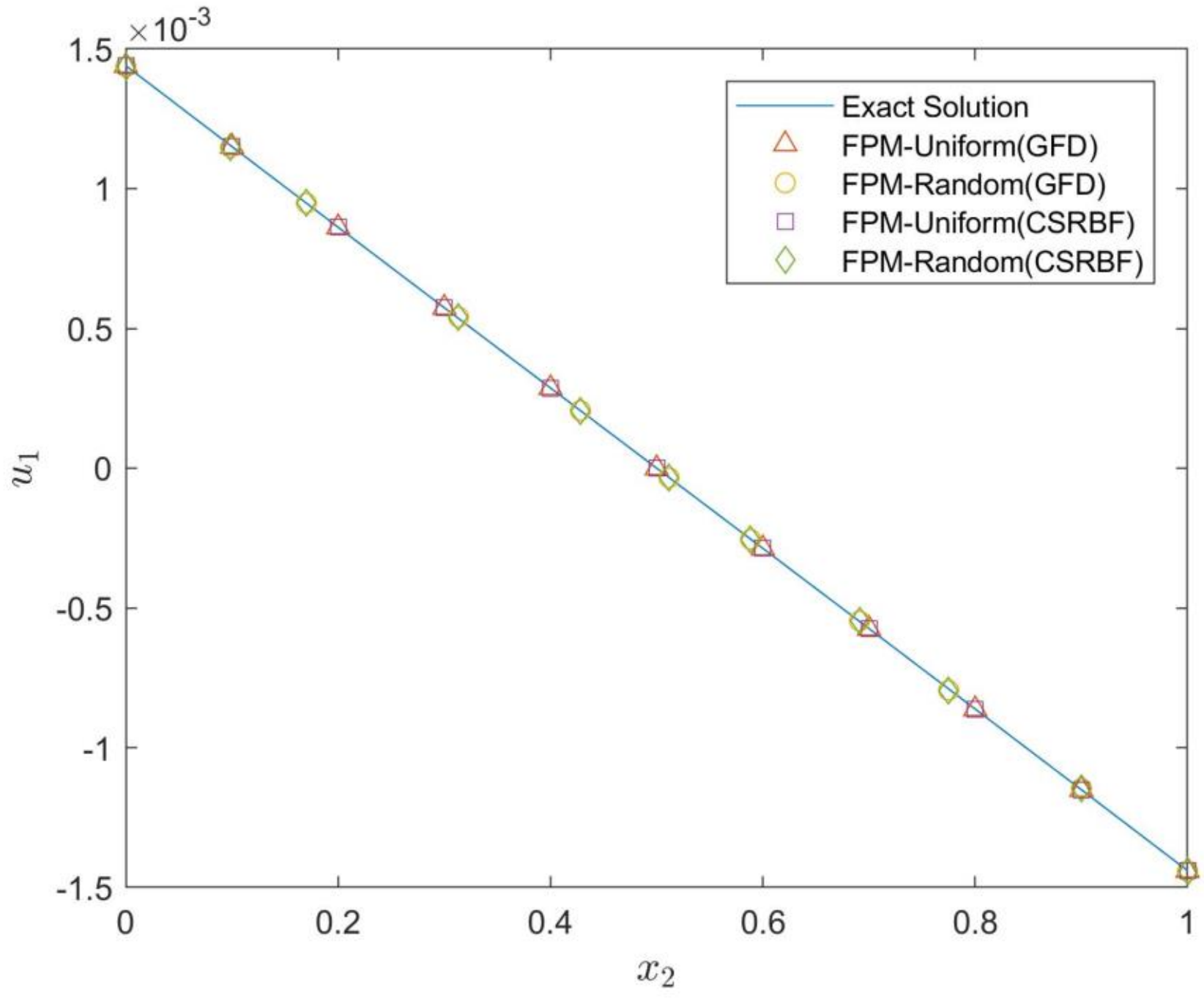

(a)

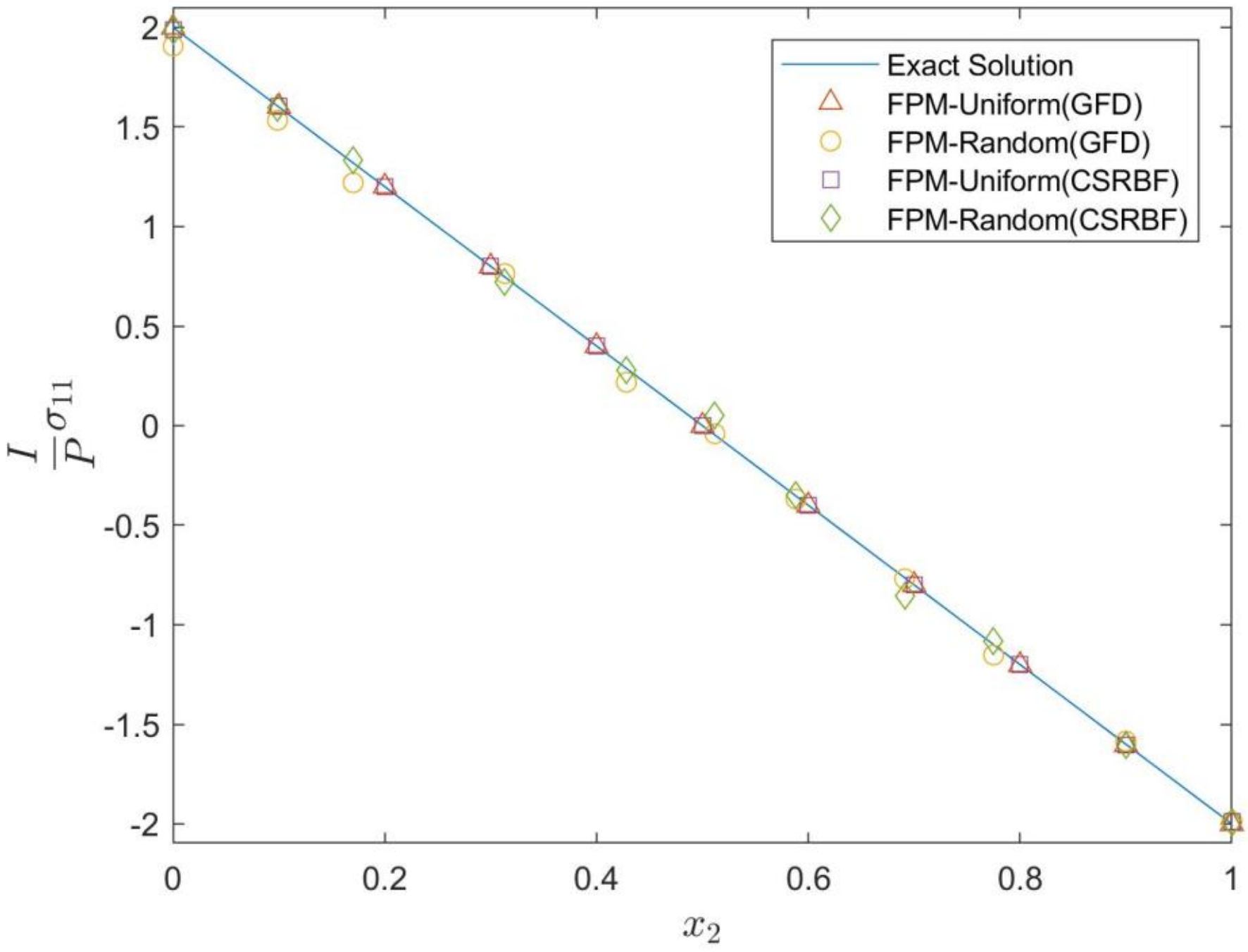

(b)

Figure 10(a). Numerical solutions of $u_1$ along $x_1 = L/2$

(b). Numerical solutions of $\sigma_{11}$ along $x_1 = L/2$

To demonstrate the convergence of the FPM, uniformly distributed sets of $41\times6$, $81\times11$ and $161\times21$ points are used to solve the problem, respectively. Besides, the Poisson's ratio is set alternatively as 0.3 or 0.4999 to test whether the FPM can be used to model nearly incompressible materials.

The relations between $h$ (the distance of two neighboring Points in $x_1$ direction) and the relative errors $r_u$, $r_E$ are shown in Figure 11 for FPM with linear trial functions which are derived from the GFD or CSRBF. For FPM with quadratic trial functions, the support domain of any Point is set to contain not only its neighboring Points but also adjacent Points around these neighboring Points. Relative errors for the computed solutions by FPM are given in Figure 12. Also, the corresponding convergence rate $R$ is given in Figure 11 and Figure 12. Compared with three-point triangular finite elements, whose convergence rates are 2 and 1 for displacements and the strain energy, respectively [1], the present FPM with linear trial functions shows a better performance in the convergence rate for the strain energy.

For the traditional FEM, when materials are nearly incompressible, volume locking leads to much smaller solutions of the displacement fields. However, from Figure 11 and Figure 12, it is obvious that the FPM performs well when $v = 0.4999$, which means that the FPM is a locking-free method for nearly incompressible materials at least for this problem

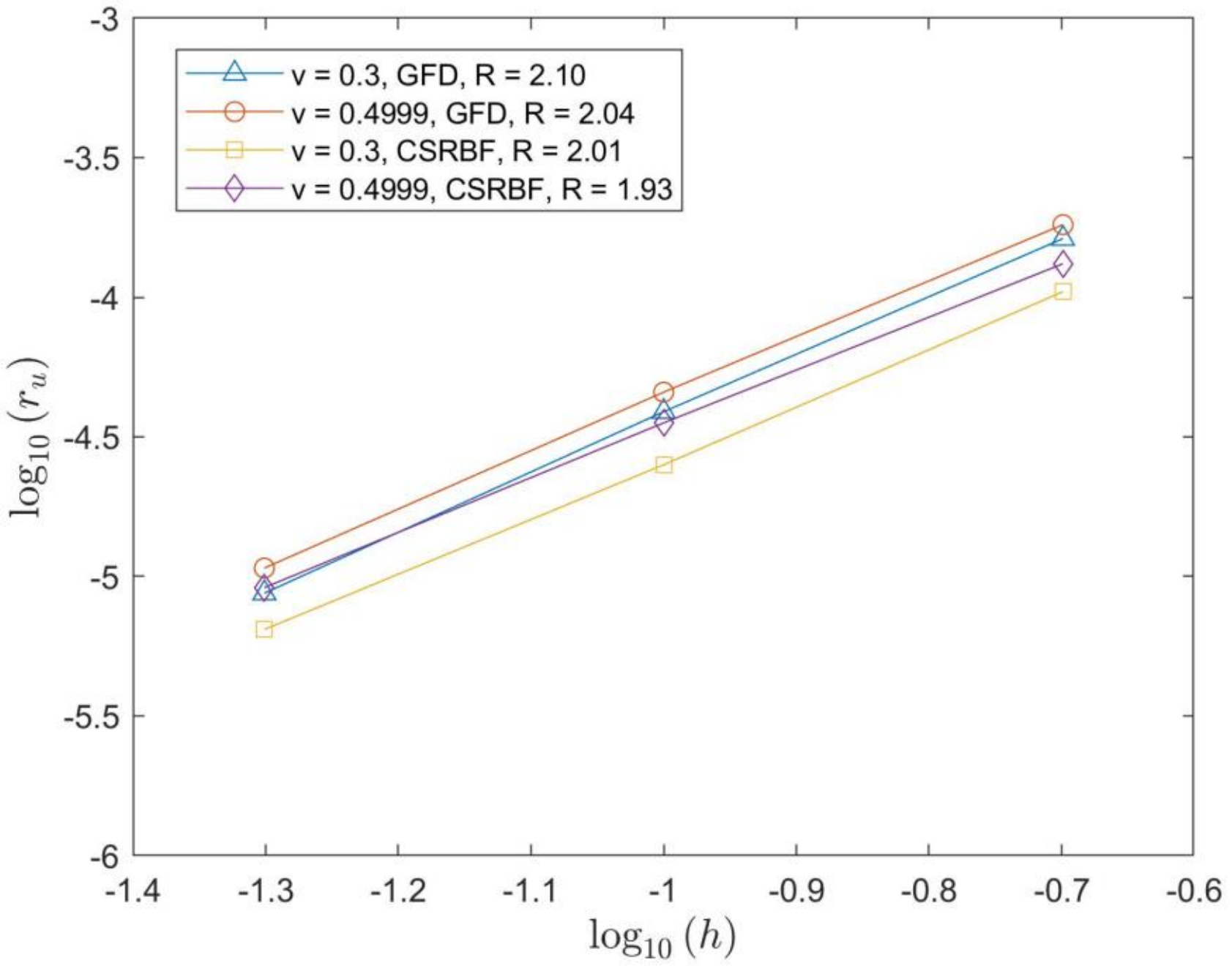

(a)

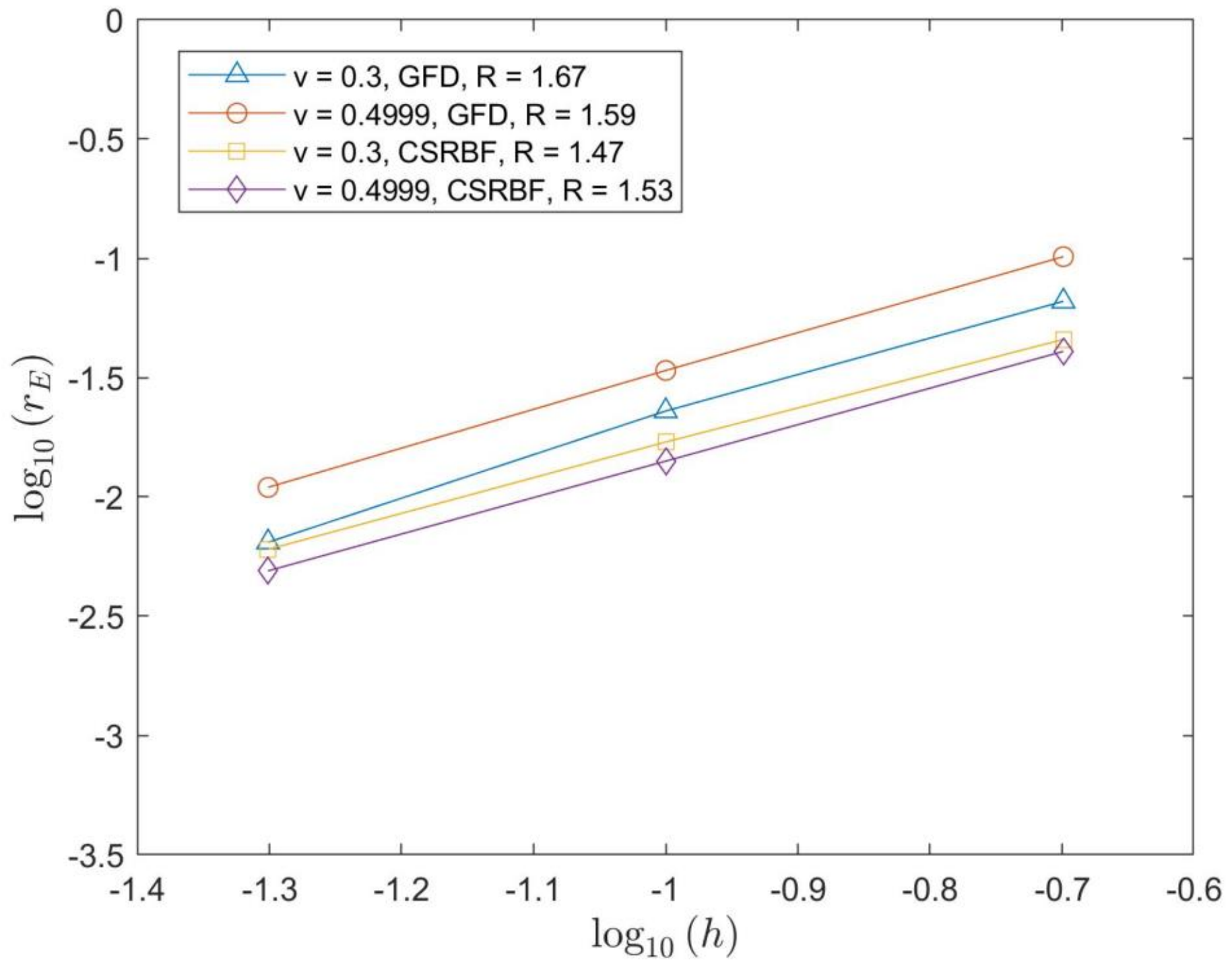

(b)

Figure 11(a). Relative errors and convergence rates for $r_u$ (linear trial functions)

(b). Relative errors and convergence rates for $r_E$ (linear trial functions)

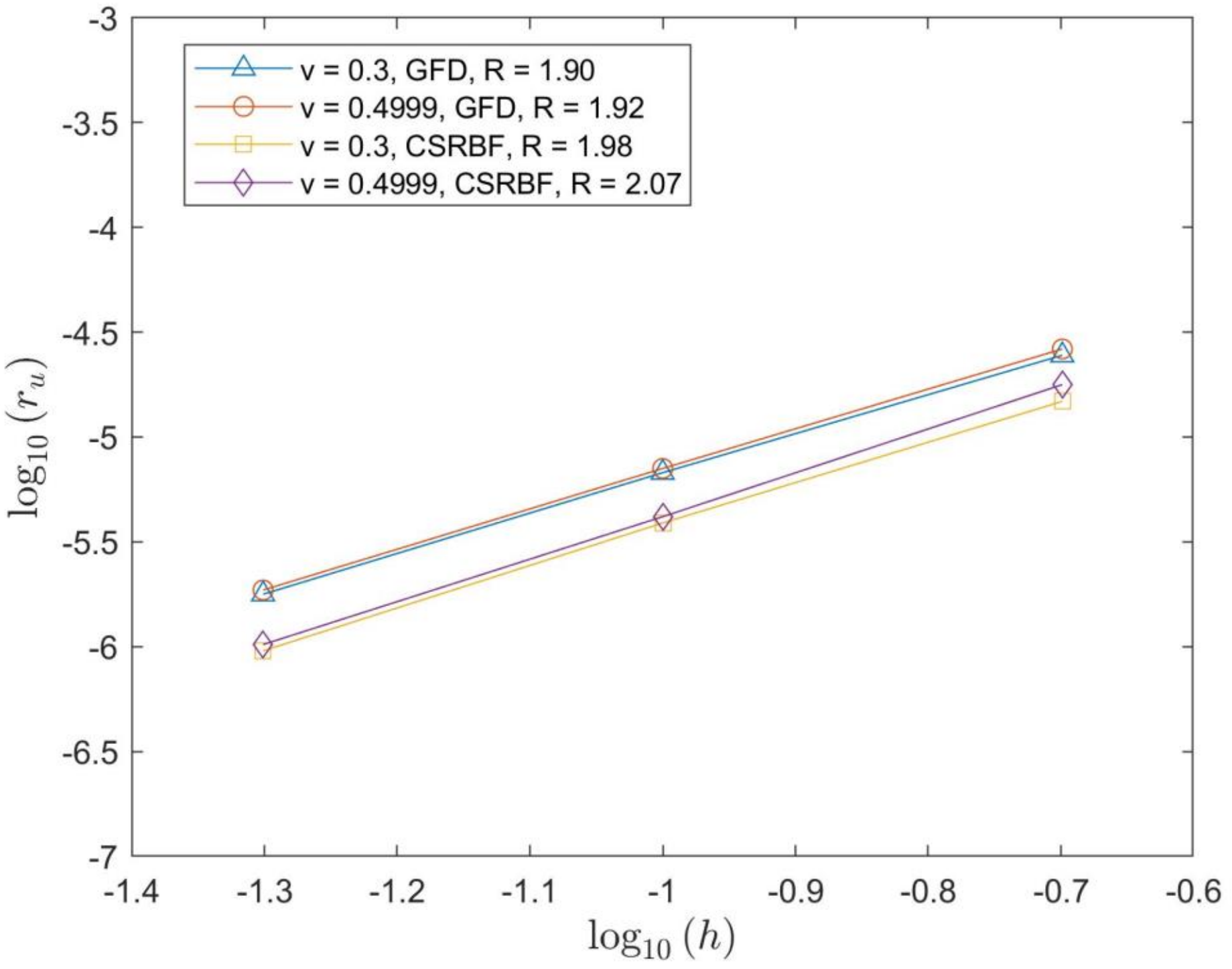

(a)

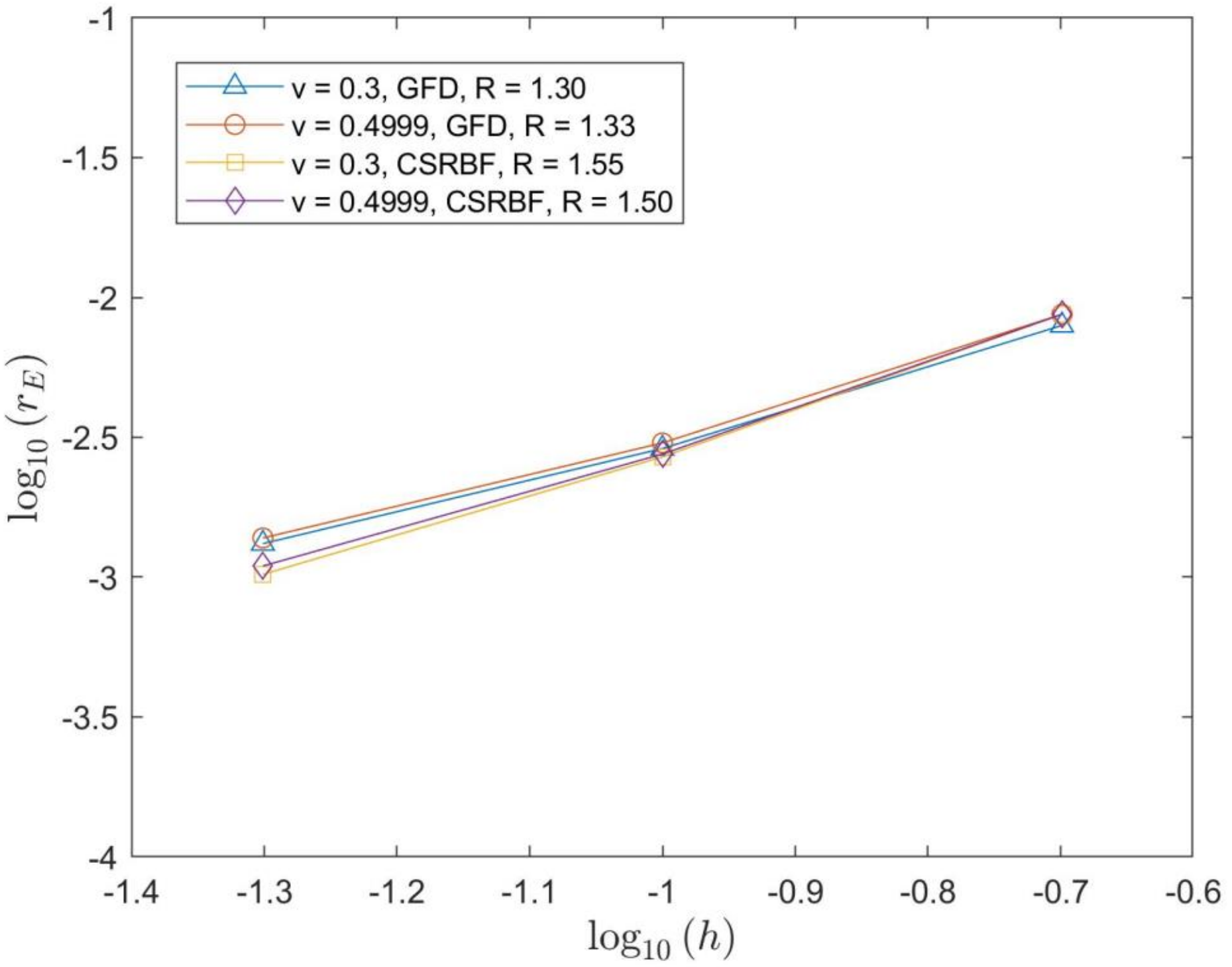

(b)

Figure 12(a). Relative errors and convergence rates for $r_u$ (quadratic trial functions)

(b). Relative errors and convergence rates for $r_E$ (quadratic trial functions)

We can see that in Section 4.1, FPMs with both the GFD and the CSRBF can pass patch tests, and in Section 4.2, relative errors and convergence rates of these two methods are similar. However, reviewing Eqs. (2.5) and (2.12), the GFD needs to invert a $2\times2$ matrix ($\mathbf{A}^{\mathrm{T}}\mathbf{W}\mathbf{A}$) for linear trial functions, while the CSRBF needs to invert a $(m+q+2)\times(m+q+2)$ matrix (**G**). Moreover, the CSRBF needs additional procedures (Eq. (2.15)) to compute displacement-derivatives at Points. For these considerations, the following numerical examples are only given by FPM with the GFD method to relate nodal displacement-derivatives to nodal DoFs, and linear trial functions are used due to its simplicity.

For this cantilever beam problem, Figure 13 further shows the relation between relative errors and the penalty parameter $\eta$ in Eq. (3.8). For this case, 161×21 points are distributed uniformly in the beam.

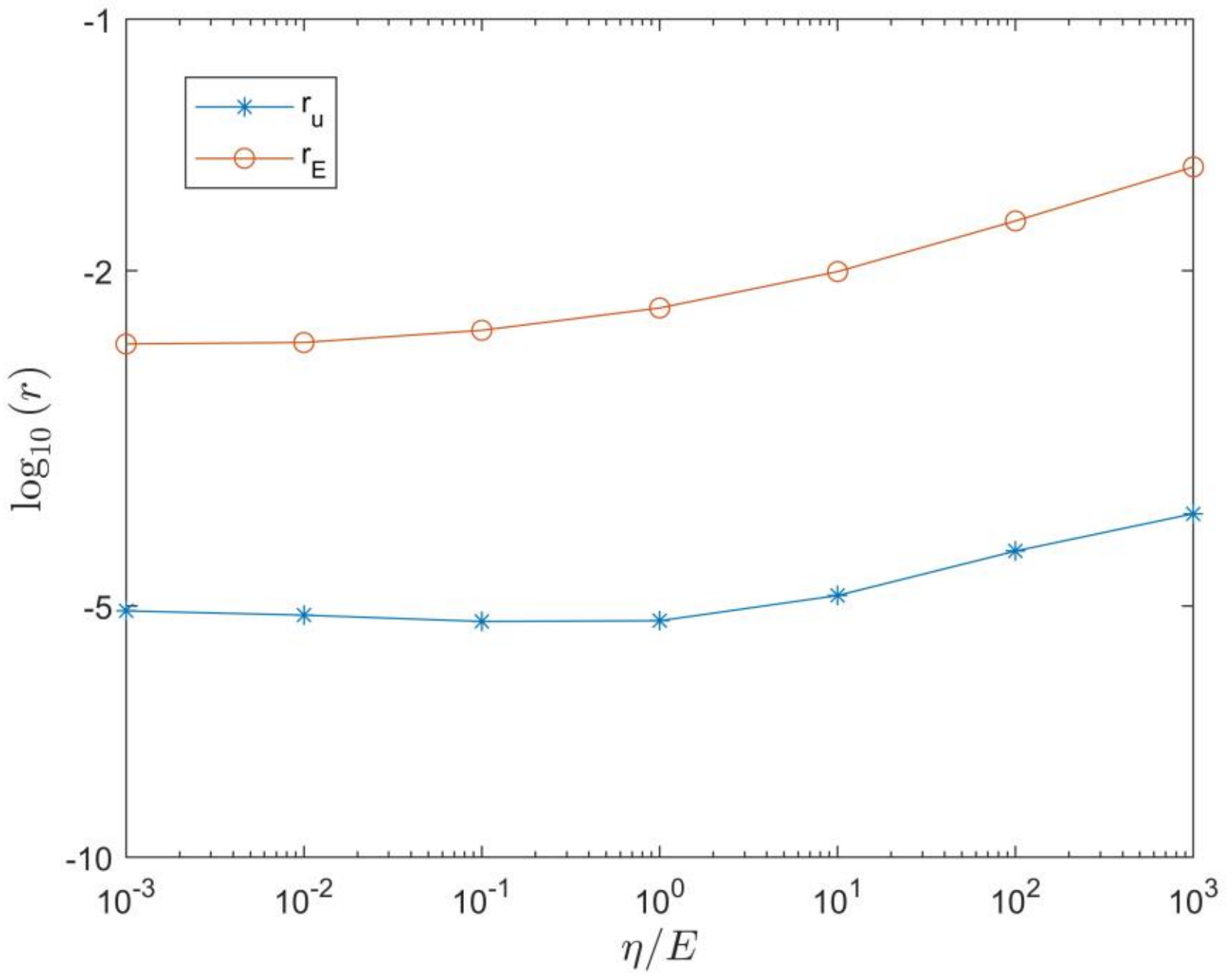

Figure 13. The relations between relative errors and $\eta$

A larger penalty coefficient can result in smaller jumps of displacements on internal boundaries, but also will increase the condition number of the global matrix, thus affecting the precision of solutions [12]. From Figure 13, we can find that with the penalty parameter changing from a small value to a large one, the relative errors stay steady at first and are then gradually increased. Based on the fact that the penalty parameter needs to be large enough to maintain the stability of the method, it is suggested to set it within the range from $10^{-2} \times E$ to $10^{2} \times E$.

It should be noticed that, the penalty parameters used for FPM can be smaller than those for discontinuous Galerkin FEM. This is because the trial functions in FPM between two neighboring subdomains are neither conforming, nor independent. This is a very important difference between FPM and discontinuous Galerkin FEM. In DG method, the neighboring elements are having entirely discontinuous and entirely independent trial functions, i.e. the DoFs in one element do not affect the field solutions in neighboring elements. In contrast, for FPM, the DoFs in this subdomain will influence the field solutions in neighboring subdomains. For this reason, one may need a large penalty factor in DG method, but one needs only a smaller penalty factor for FPM.

4.3 A Cook's Skew Beam

In this section, the FPM is employed to simulate the Cook's skew beam problem [14]. As shown in Figure 14, the beam is fixed at the left end, with its right end subjected to a uniform shear traction $F = 1/16$. The material is isotropic linear elastic with Young's modulus $E = 1.0$ and Poisson's ration $v = 0.3$.

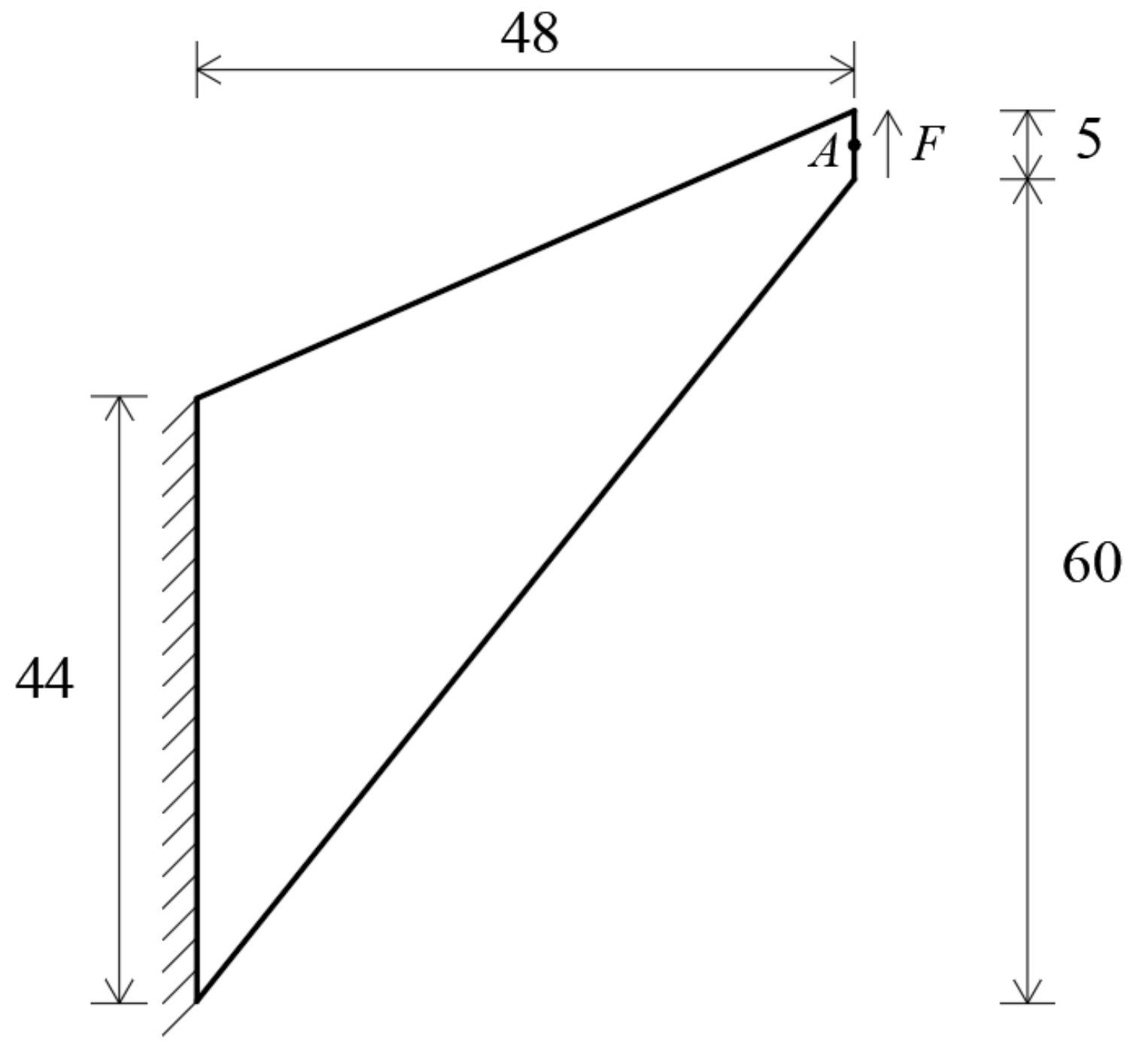

Figure 14. Cook's skew beam

This structure is simulated by FPM using linear and quadratic trial functions, respectively, and for comparison it is also modeled in ABAQUS by FEM with the same mesh (shown in Figure 15), using CPS4 and CPS8 elements.

Numerical results of $u_2$ at Point A ($u_2^A$) obtained by linear FPM, quadratic FPM, FEM with CPS4 and CPS8 elements are given in Table III. For comparison, the numerical solution of $u_2^A$ obtained by FEM with 6370 CPS8 elements is considered as the reference value.

Numerical solutions of $\sigma_{11}$ obtained by FPM and FEM are illustrated in Figure 15. For comparison, Numerical solutions obtained by FEM with 6370 CPS8 elements are also shown here.

Table III. Numerical solutions of $u_2^A$

| | FPM | | | FEM | | Reference Value |
|---|---|---|---|---|---|---|
| | Linear | Quadratic | Linear (with very distorted mesh) | CPS4 | CPS8 | |
| $u_2^A$ | 16.277 | 20.078 | 16.519 | 14.284 | 19.777 | 19.869 |

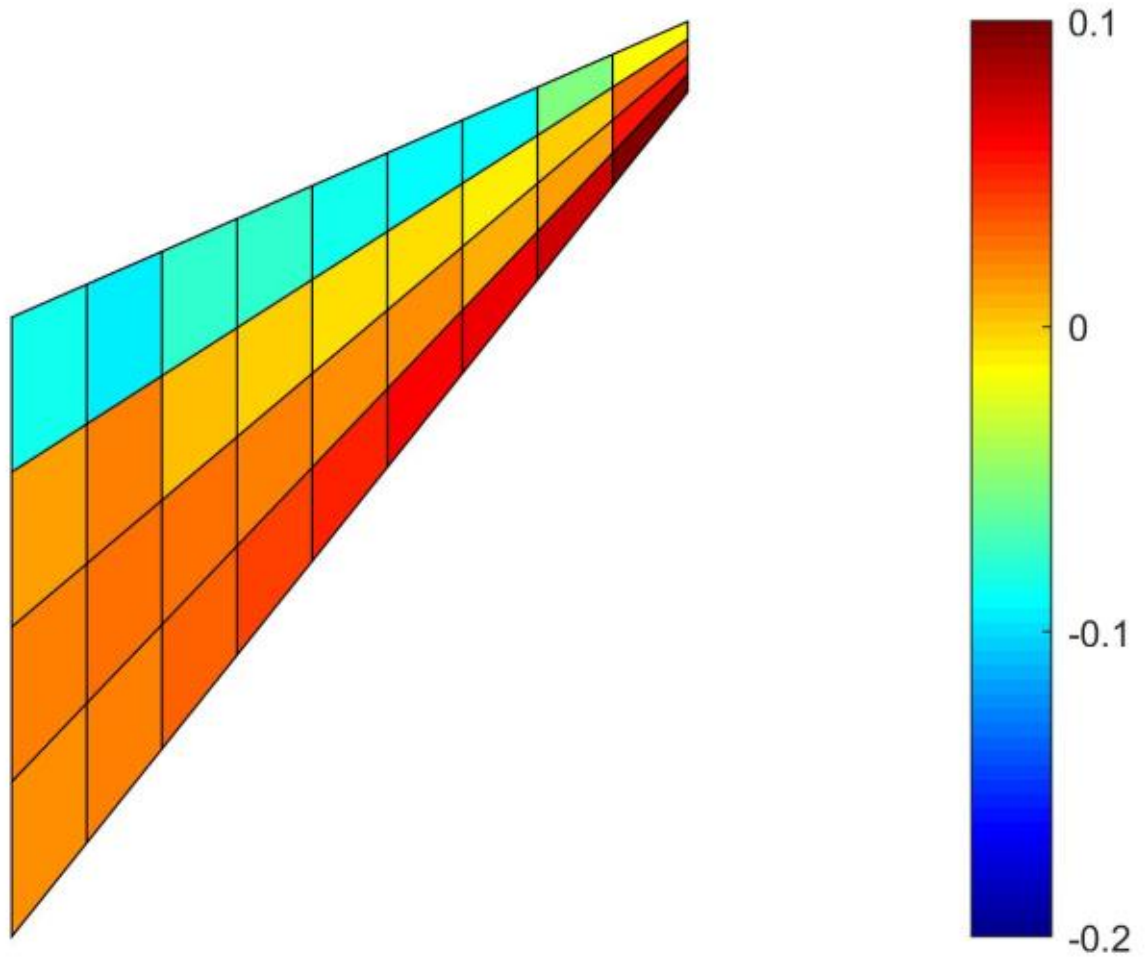

(a)

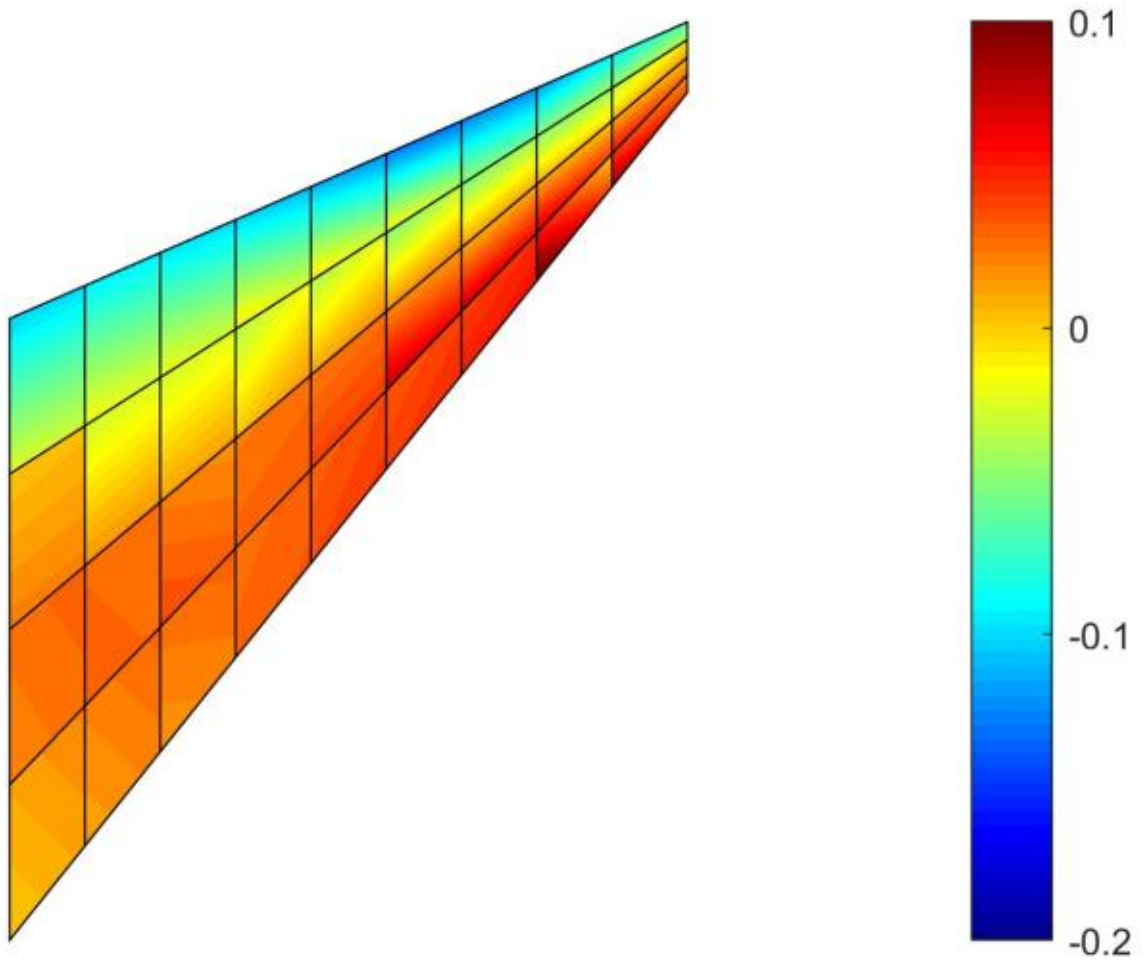

(b)

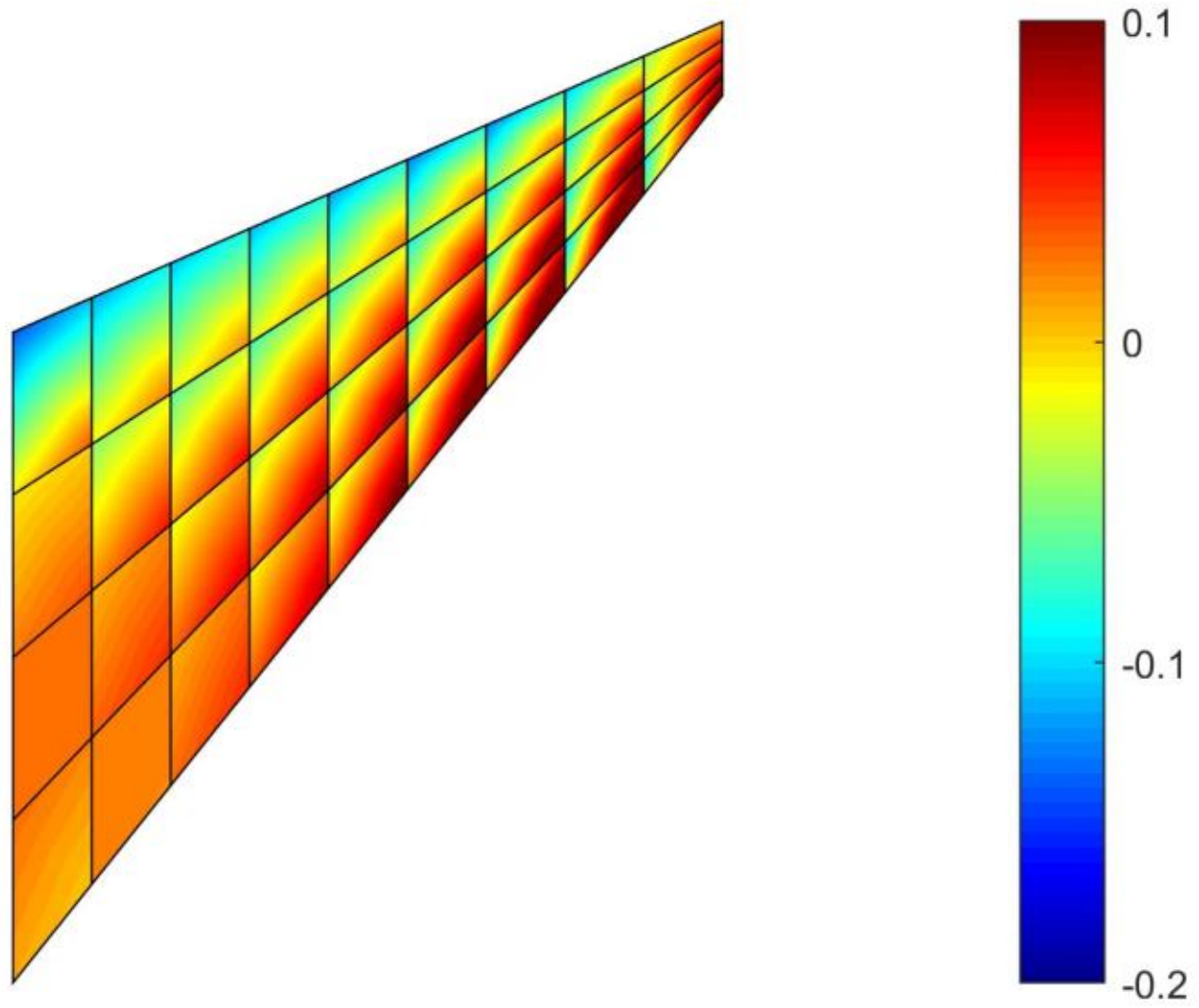

(c)

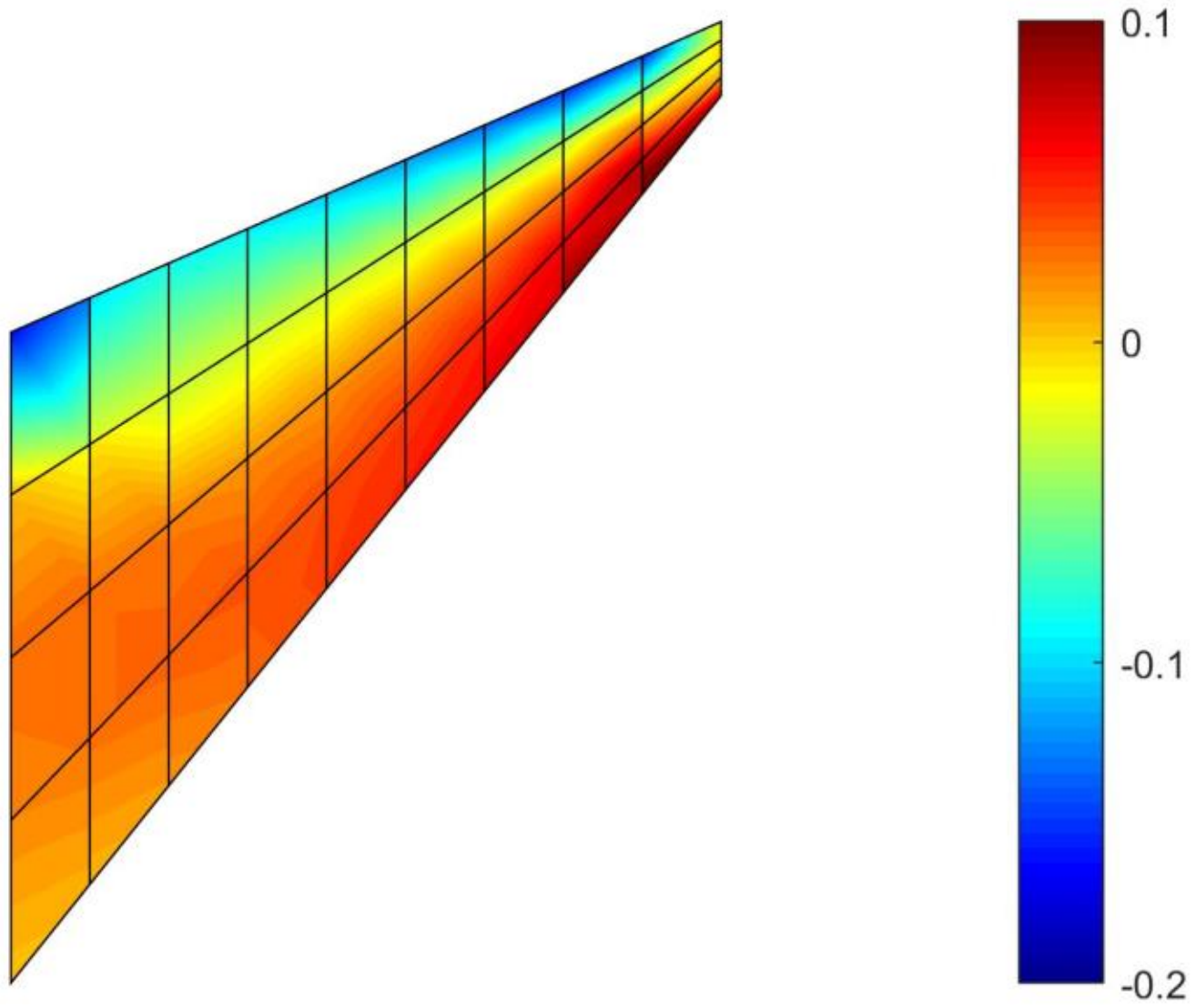

(d)

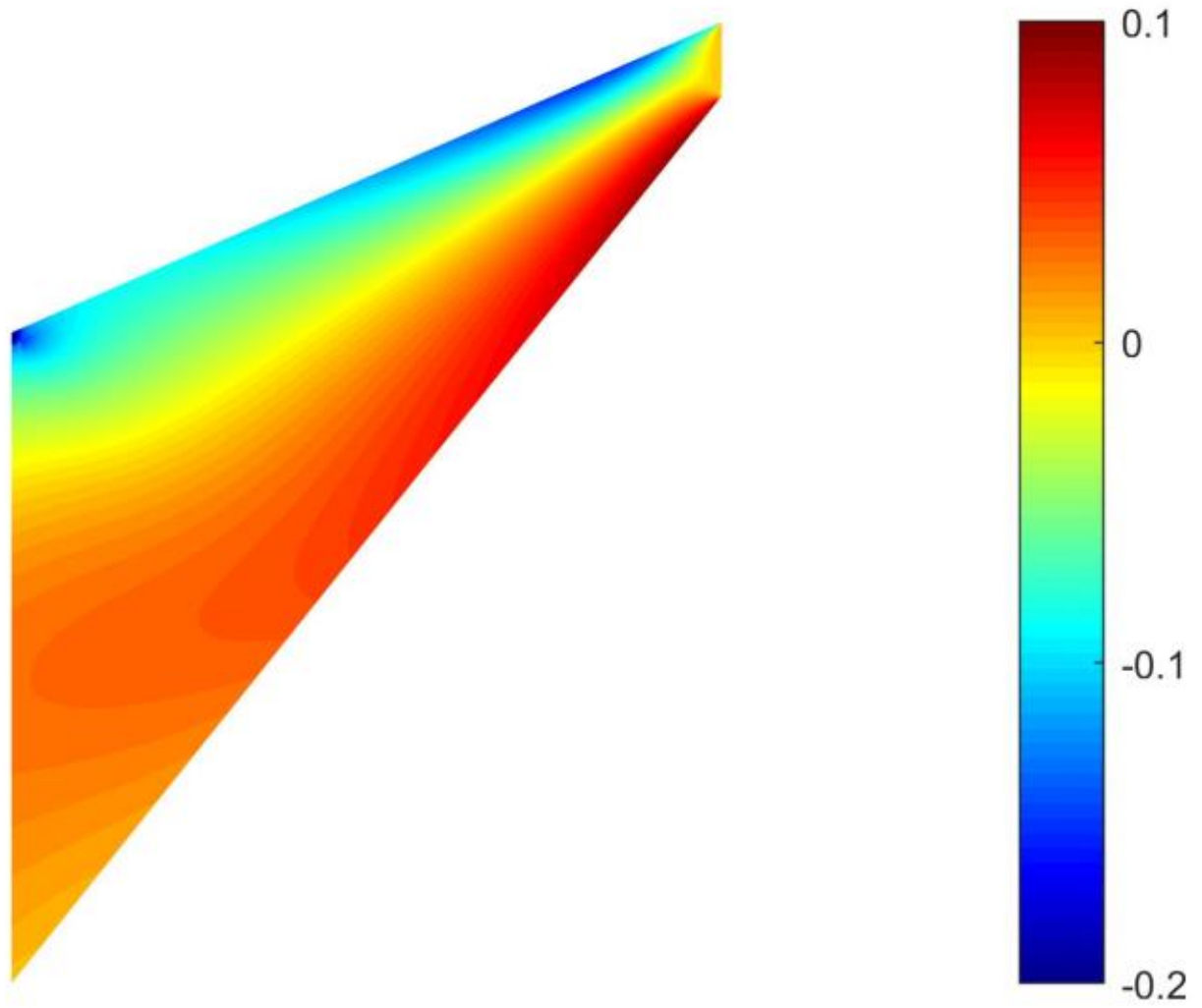

(e)

Figure 15. Numerical solutions of $\sigma_{11}$

(a)Linear FPM (b) Quadratic FPM

(c)FEM with CPS4 elements (d) FEM with CPS8 elements (e) Reference solution

We can see that for this cook's skew beam problem, quadratic FPM and FEM with CPS8 elements both lead to satisfactory solutions with few elements, but linear FPM performs much better than FEM with CPS4 elements (the standard four-node quadrilateral element). This is because shape functions in FPM are established by Point-based methods rather than element-based ones. Therefore, FPM is expected to be much less sensitive to mesh-distortion. As a demonstration, the Cook's skew beam is also modeled by FPM with one concave quadrilateral element inside, as shown in (b)

Figure 16. Numerical solutions of $u_2$ and $\sigma_{11}$ obtained by FPM with this very distorted mesh are illustrated (b)

Figure 16. Besides, the value of $u_2$ at Point A is given in Table III. We can see that even with a concave element inside, the FPM can still give satisfactory results. But for FEM, a concave element will lead to wrong solutions which is not allowed in FEM

software.

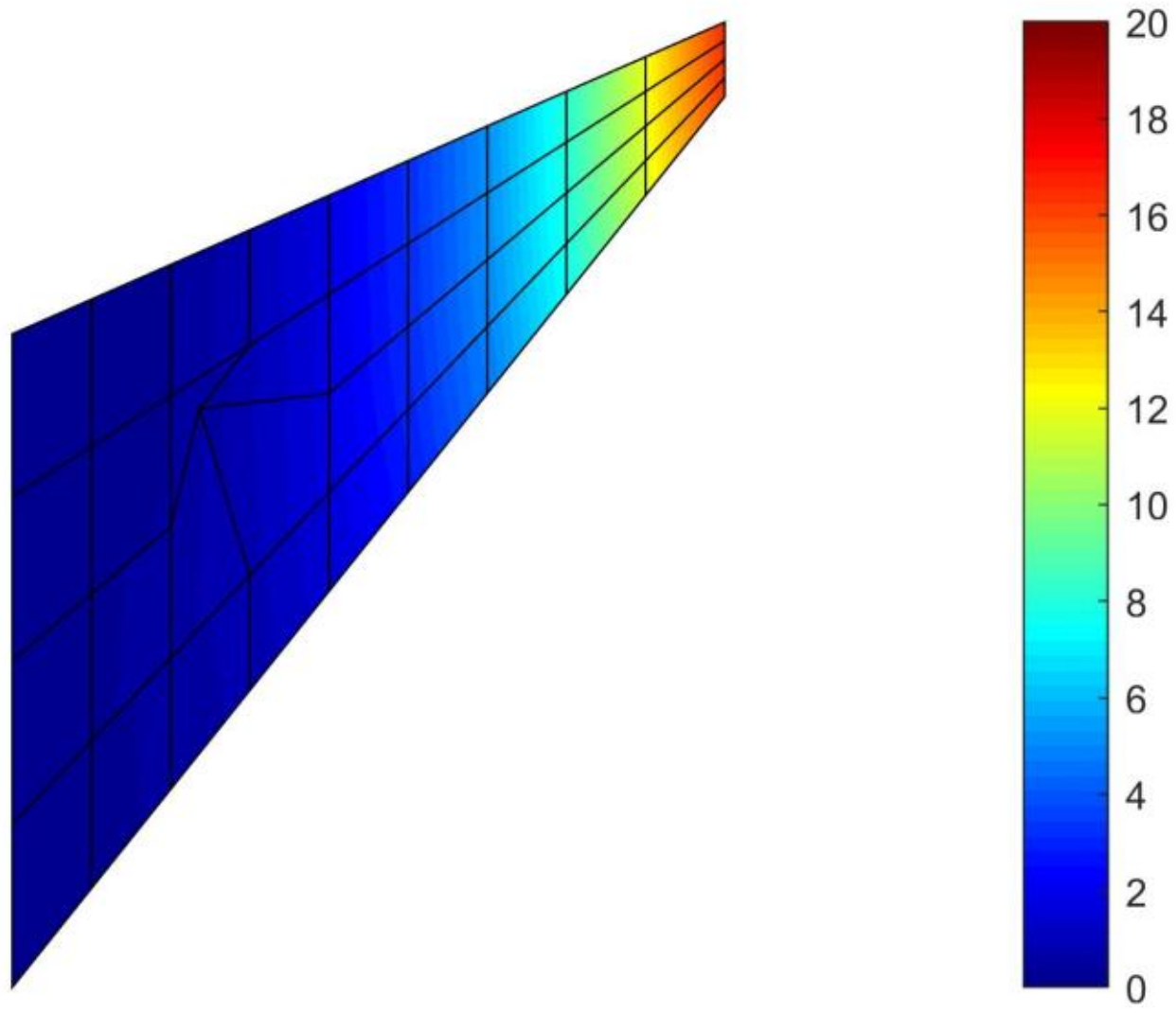

(a)

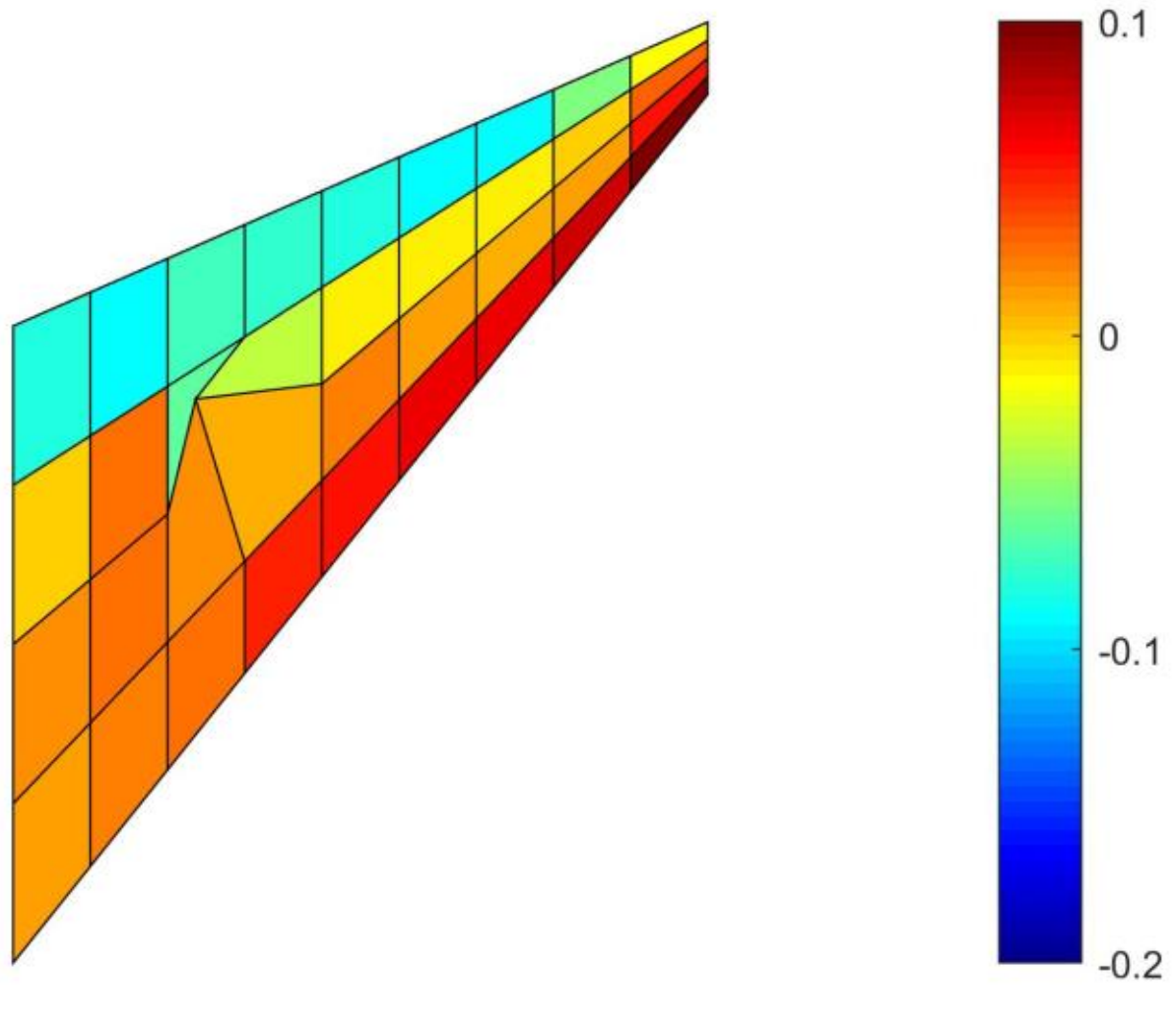

(b)

Figure 16. Numerical solutions by FPM with linear trial functions

(a) $u_2$ (b) $\sigma_{11}$

4.4 A Ring with Radial Tensile Traction

In this subsection, a ring with radial tensile traction is solved by the FPM (shown in Figure 17(a)). The ring is defined as $\{(x_1, x_2) | a^2 \le x_1^2 + x_2^2 \le b^2\}$ and it is subjected to a uniform radial tension. Since the ring is symmetric in geometry, we only model the upper right quarter (shown in Figure 17(b)). Symmetry boundary conditions are imposed on the left and bottom edges, which means $u_1 = 0$, $t_2 = 0$ for the left edge and $u_2 = 0$, $t_1 = 0$ for the bottom edge. Traction boundary conditions are imposed at *r* = *b*, according to the tensile traction *p*. The edge, *r* = *a*, is set to be a traction-free boundary.

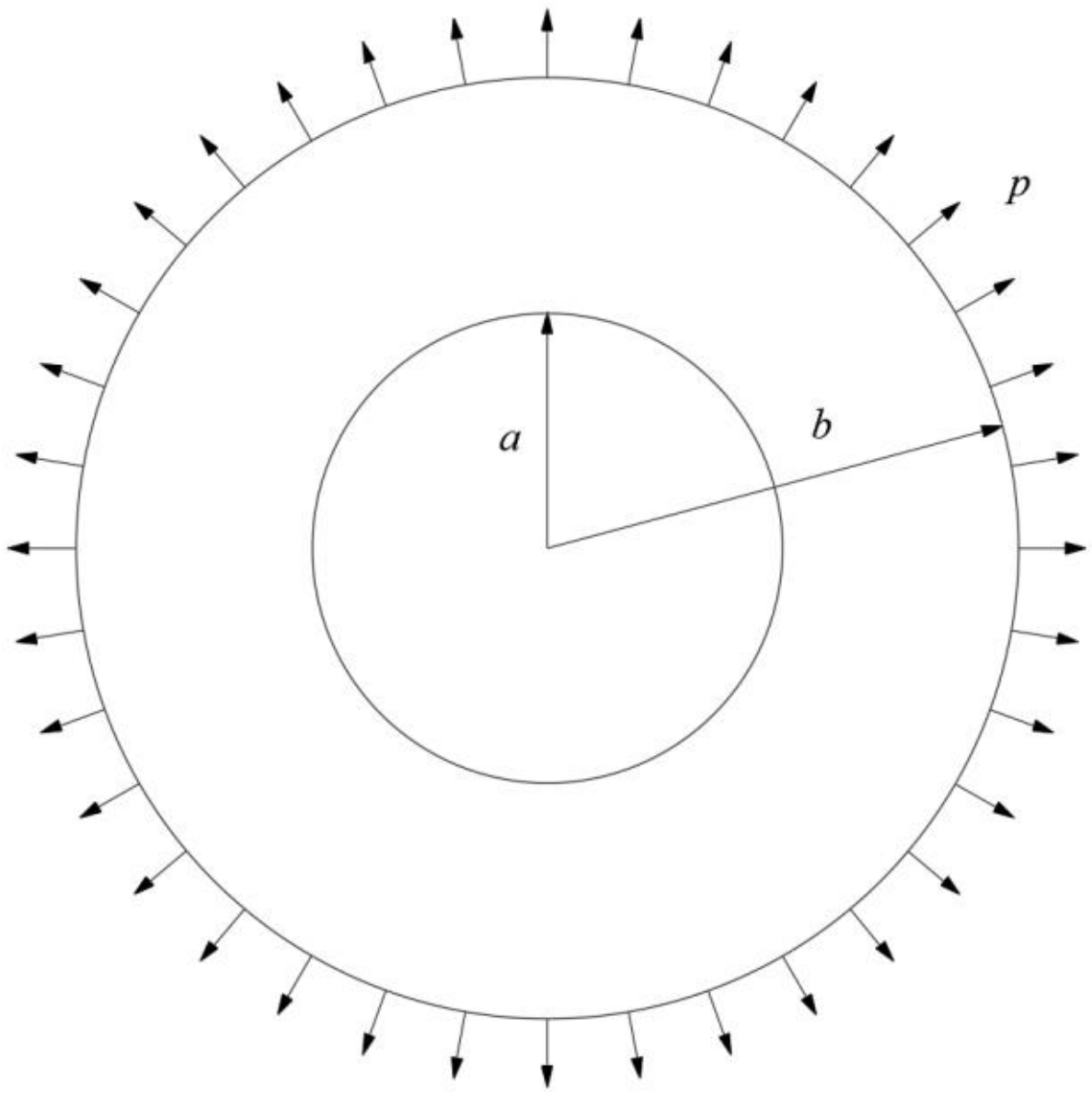

(a)

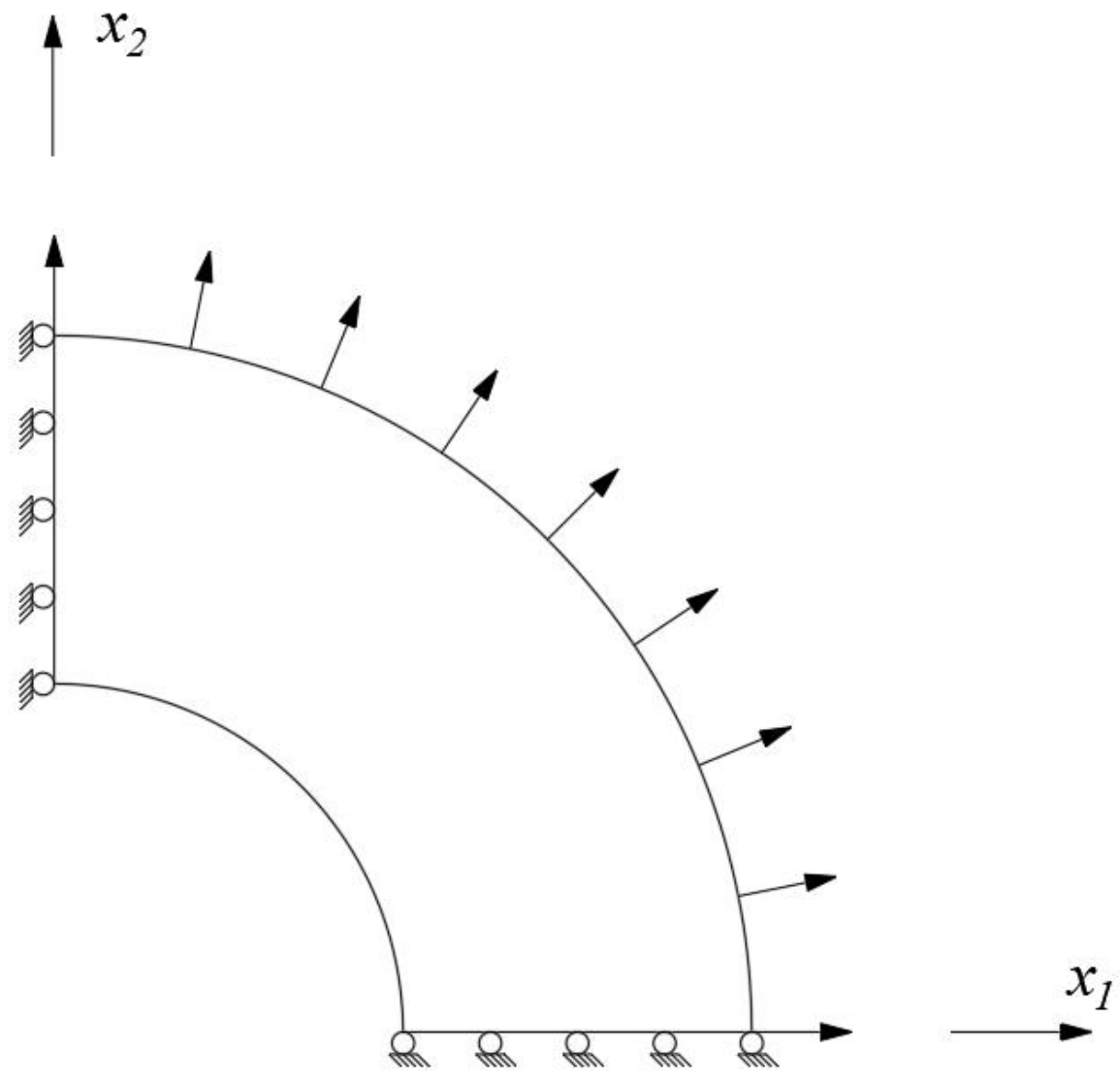

(b)

Figure 17(a). The ring with radial tension

(b). The ring on fourfold symmetry

Specifically, we prescribe that $a = 1$, $b = 2$ and $p = 1$ for the ring. The exact solutions for stresses and displacements are given in Eqs. (4.4) and (4.5), respectively, where $(r,\theta)$ are the polar coordinates and $\theta$ is anticlockwise measured from the positive $x_1$ -axis. The problem is solved considering a Plane Stress condition with $E = 1\times10^5$ and $v = 0.3$. The penalty parameter is set as $\eta = E$.

$$\sigma_{rr} = \frac{b^2}{b^2 - a^2}\left(1 - \frac{a^2}{r^2}\right)p$$

$$\sigma_{\theta\theta} = \frac{b^2}{b^2 - a^2}\left(1 + \frac{a^2}{r^2}\right)p \tag{4.4}$$

$$\sigma_{r\theta} = 0$$

$$u_r = \frac{1}{E}\left(\frac{(1-v)b^2 p}{b^2 - a^2} r + \frac{(1+v)a^2 b^2}{b^2 - a^2}\frac{1}{r}\right) \tag{4.5}$$

$$u_\theta = 0$$

To study the convergence of the FPM for displacements and the strain energy, regularly distributed sets of 15×11, 15×16 and 15×21 points are considered (shown in Figure 18(a)). The relation between $h$, defined as the longest distance between two neighbouring Points, and relative errors are given in Figure 19, where $R$ stands for the convergence rate.

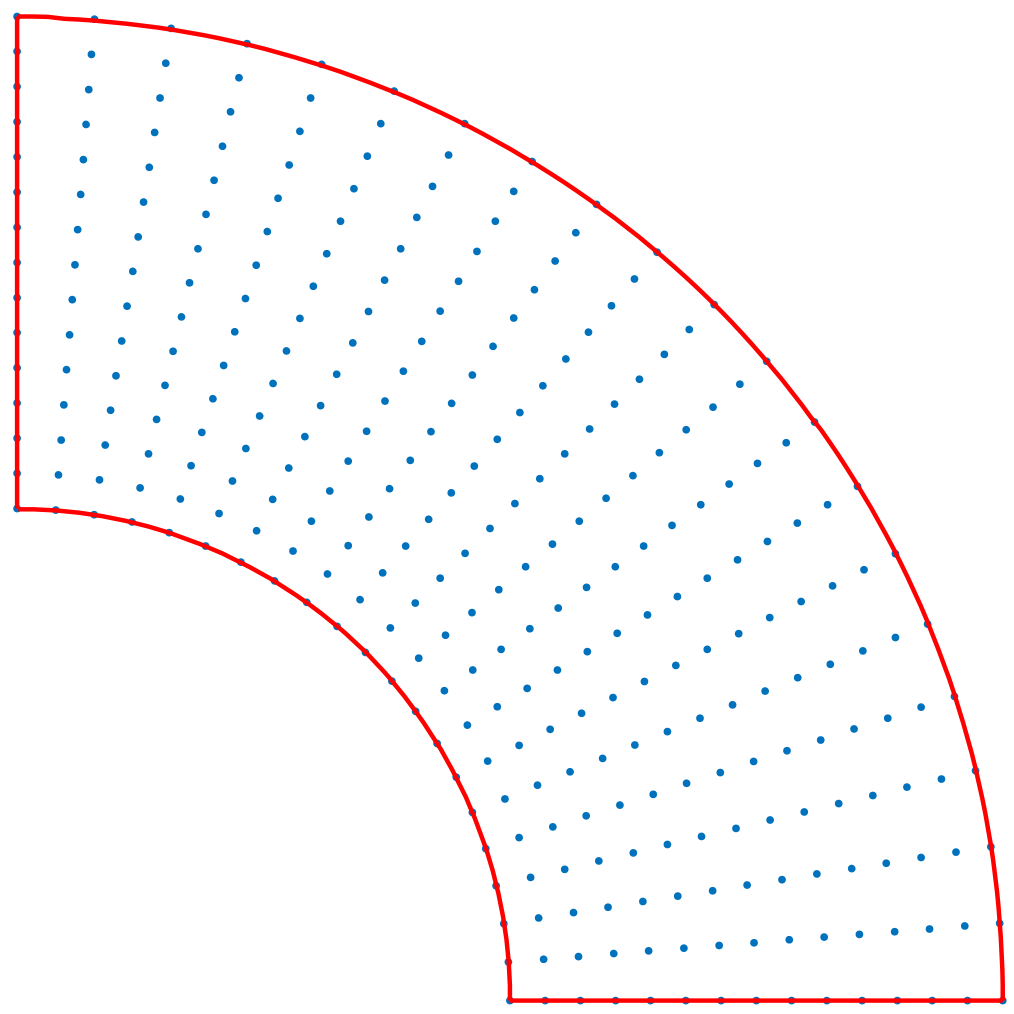

(a)

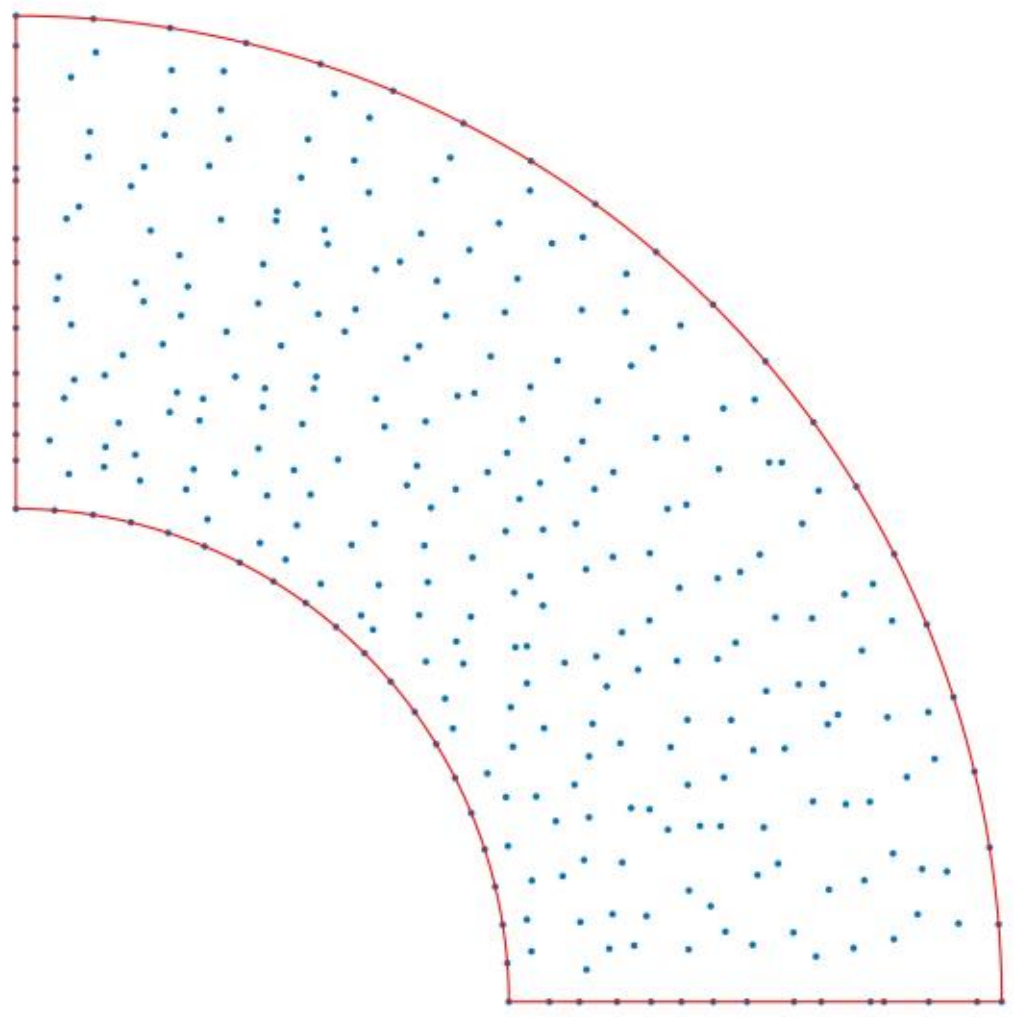

(b)

Figure 18(a). The regular distribution of 15×21 points

(b). The irregular distribution of 15×21 points

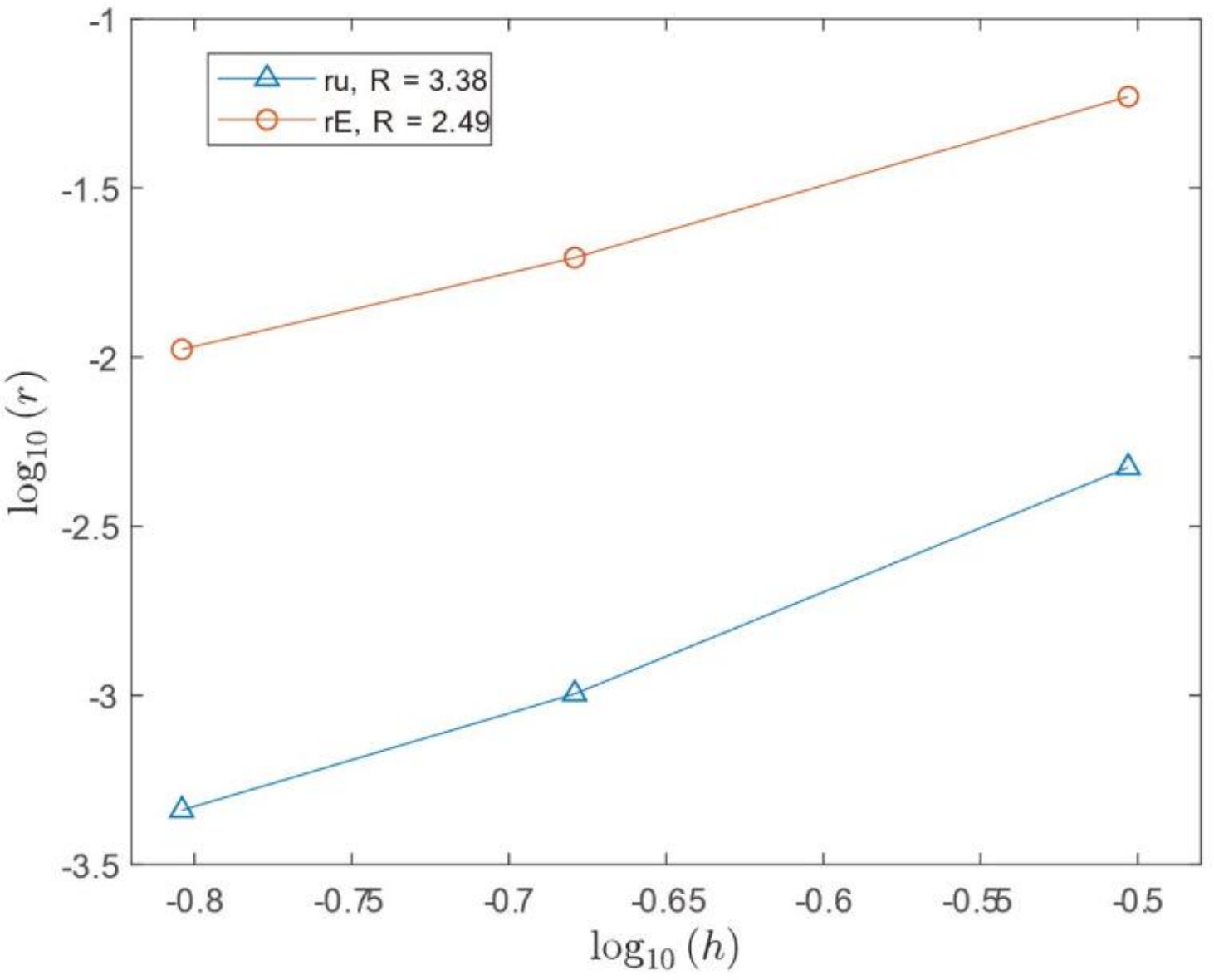

Figure 19. Relative errors and convergence rates for the problem

Additionally, numerical solutions of $u_r$ and $\sigma_{rr}$ along $\theta = 0$ with 15×21 Points distributed regularly and randomly are shown in Figure 20(a) and (b), respectively. The irregular distribution of Points is illustrated in Figure 18(b). It can be seen that the FPM can solve the displacements and stresses with a satisfactory accuracy for both the situations.

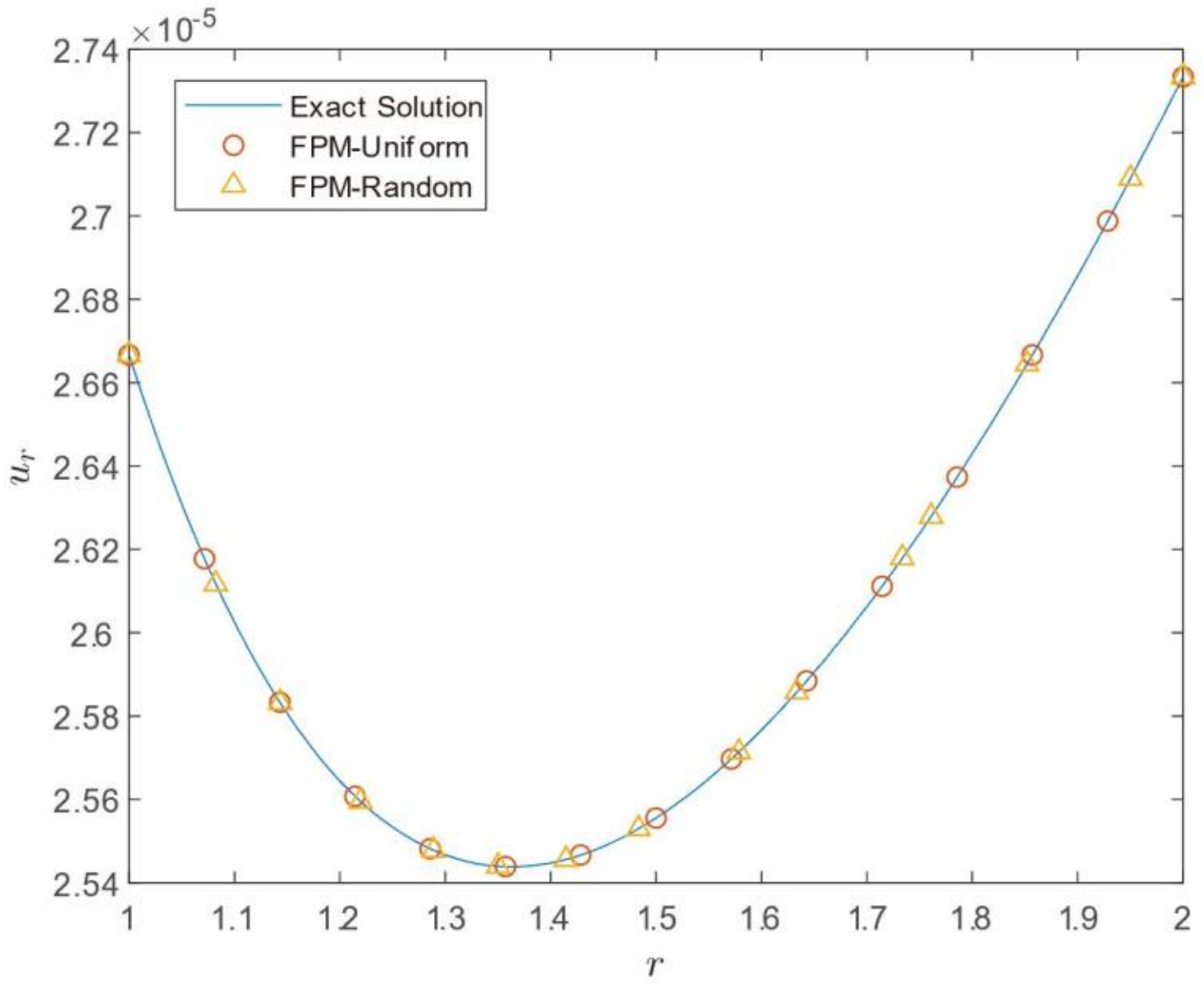

(a)

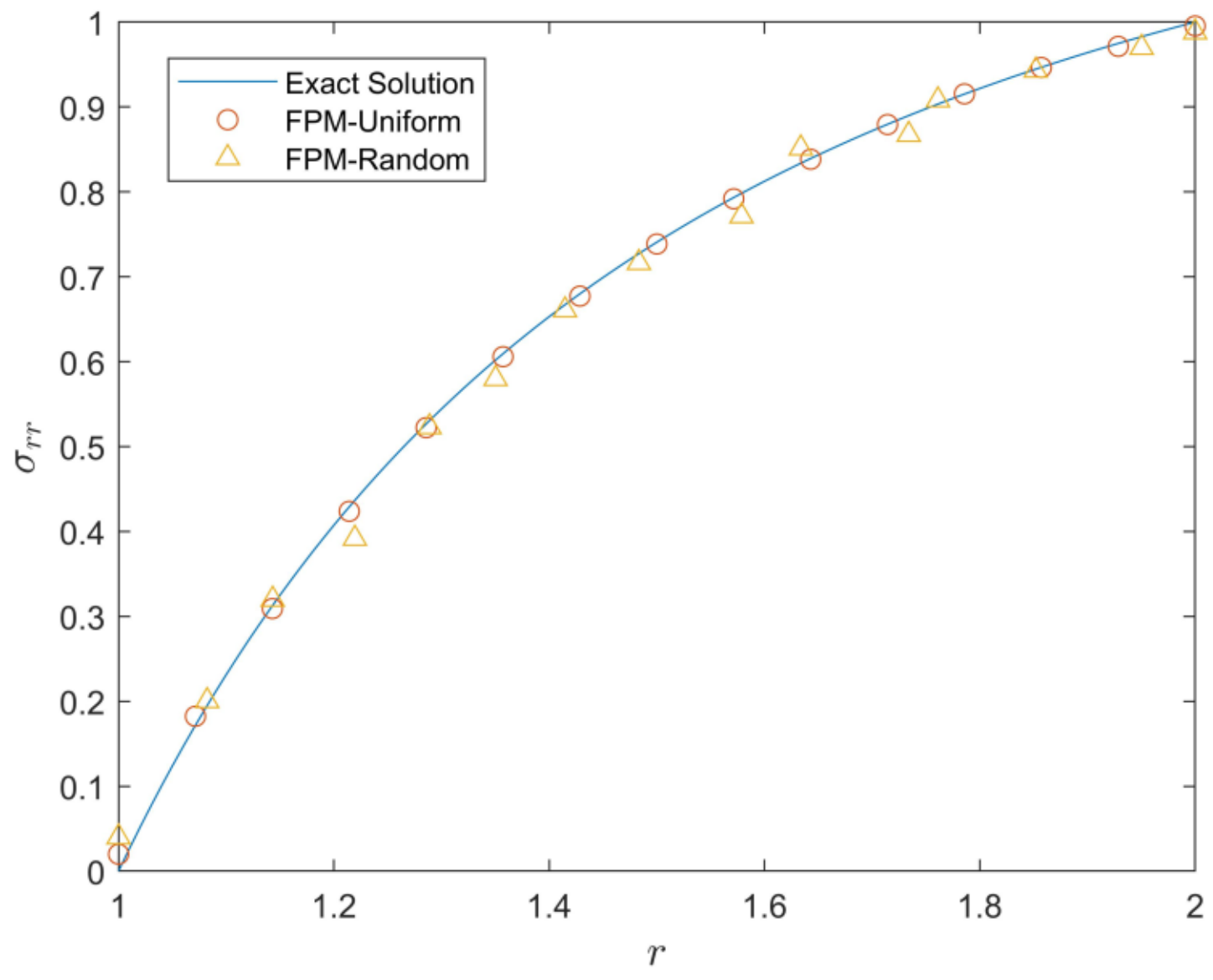

(b)

Figure 20(a). Numerical solutions of $u_r$ along $\theta = 0$

(b). Numerical solutions of $\sigma_{rr}$ along $\theta = 0$

### 4.5 An Infinite Plate with A Circular Hole

In this subsection, we employ the FPM to model an infinite plate with a circular hole. As shown in Figure 21(a), the circular hole (radius equals $a$) is located at the plate's center and a uniform tensile stress $p$ is imposed in the $x_1$ direction at infinity. The exact solutions for stresses and displacements are given in Eqs. (4.6) and (4.7) respectively.

$$\sigma_{11} = p\left[1 - \frac{a^2}{r^2}\left(\frac{3}{2}\cos 2\theta + \cos 4\theta\right) + \frac{3a^4}{2r^4}\cos 4\theta\right]$$

$$\sigma_{12} = p\left[-\frac{a^2}{r^2}\left(\frac{1}{2}\sin 2\theta + \sin 4\theta\right) + \frac{3a^4}{2r^4}\sin 4\theta\right] \tag{4.6}$$

$$\sigma_{22} = p\left[-\frac{a^2}{r^2}\left(\frac{1}{2}\cos 2\theta - \cos 4\theta\right) - \frac{3a^4}{2r^4}\cos 4\theta\right]$$

$$u_1 = \frac{1+v}{E}\, p\,\left(\frac{1}{1+v} r\cos\theta + \frac{2}{1+v}\frac{a^2}{r}\cos\theta \;+ \frac{1}{2}\frac{a^2}{r}\cos 3\theta - \frac{1}{2}\frac{a^4}{r^3}\cos 3\theta\right) \tag{4.7}$$

$$u_2 = \frac{1+v}{E} p \left( \frac{-v}{1+v} r \sin\theta - \frac{1-v}{1+v} \frac{a^2}{r} \sin\theta + \frac{1}{2} \frac{a^2}{r} \sin 3\theta - \frac{1}{2} \frac{a^4}{r^3} \sin 3\theta \right)$$

Based on the symmetry of the problem, we simplify the model by considering a quarter of the plate, as shown in Figure 21(b). Symmetry boundary conditions $u_1 = 0$, $t_2 = 0$ at $x_1 = 0$ and $u_2 = 0$, $t_1 = 0$ at $x_2 = 0$, are imposed. Displacement boundary conditions are imposed on the upper side ($x_2 = 2$) and right side ($x_1 = 2$) according to Eq. (4.7). The edge at $r = 1$ is set to be a traction-free boundary.

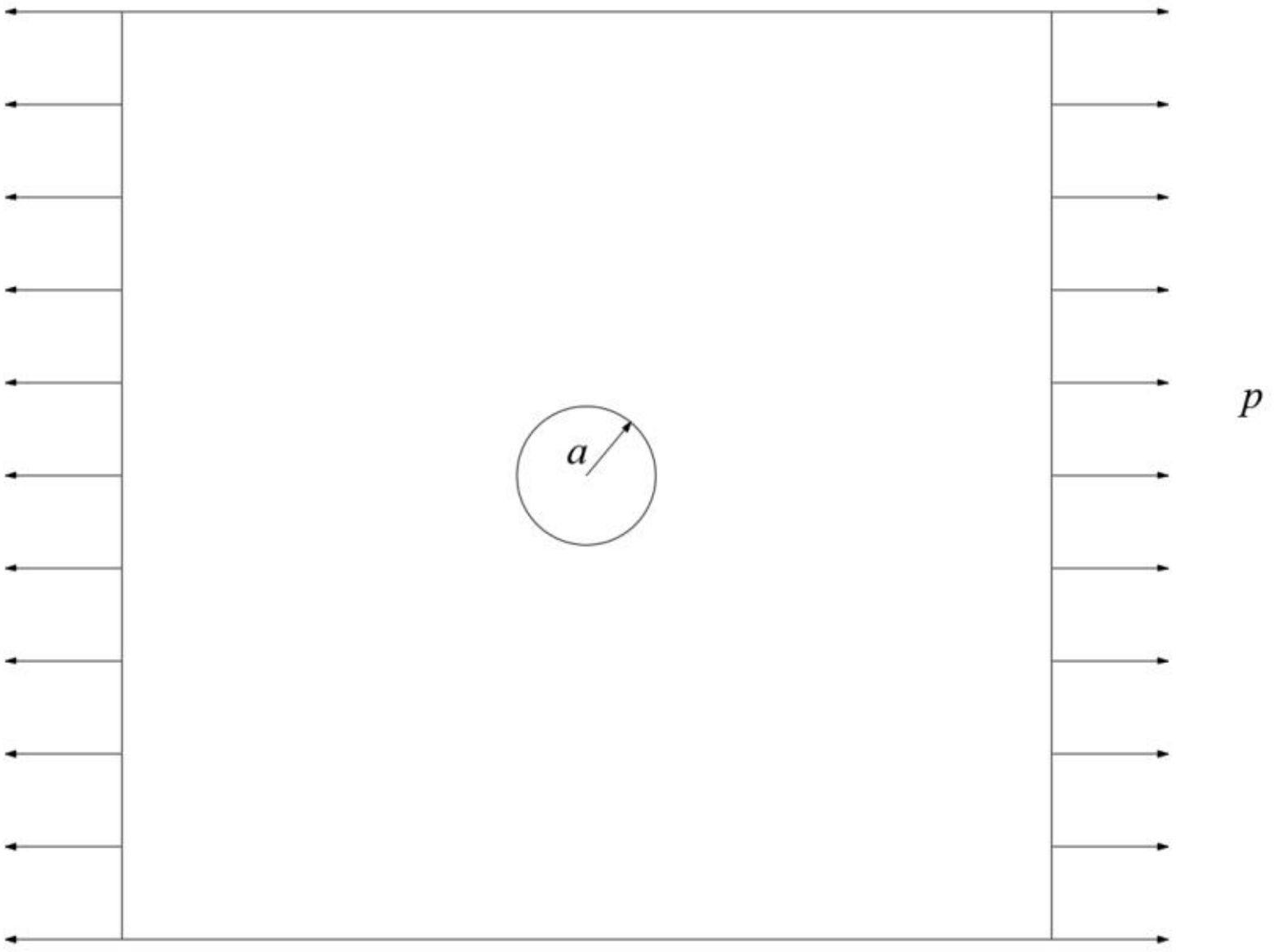

(a)

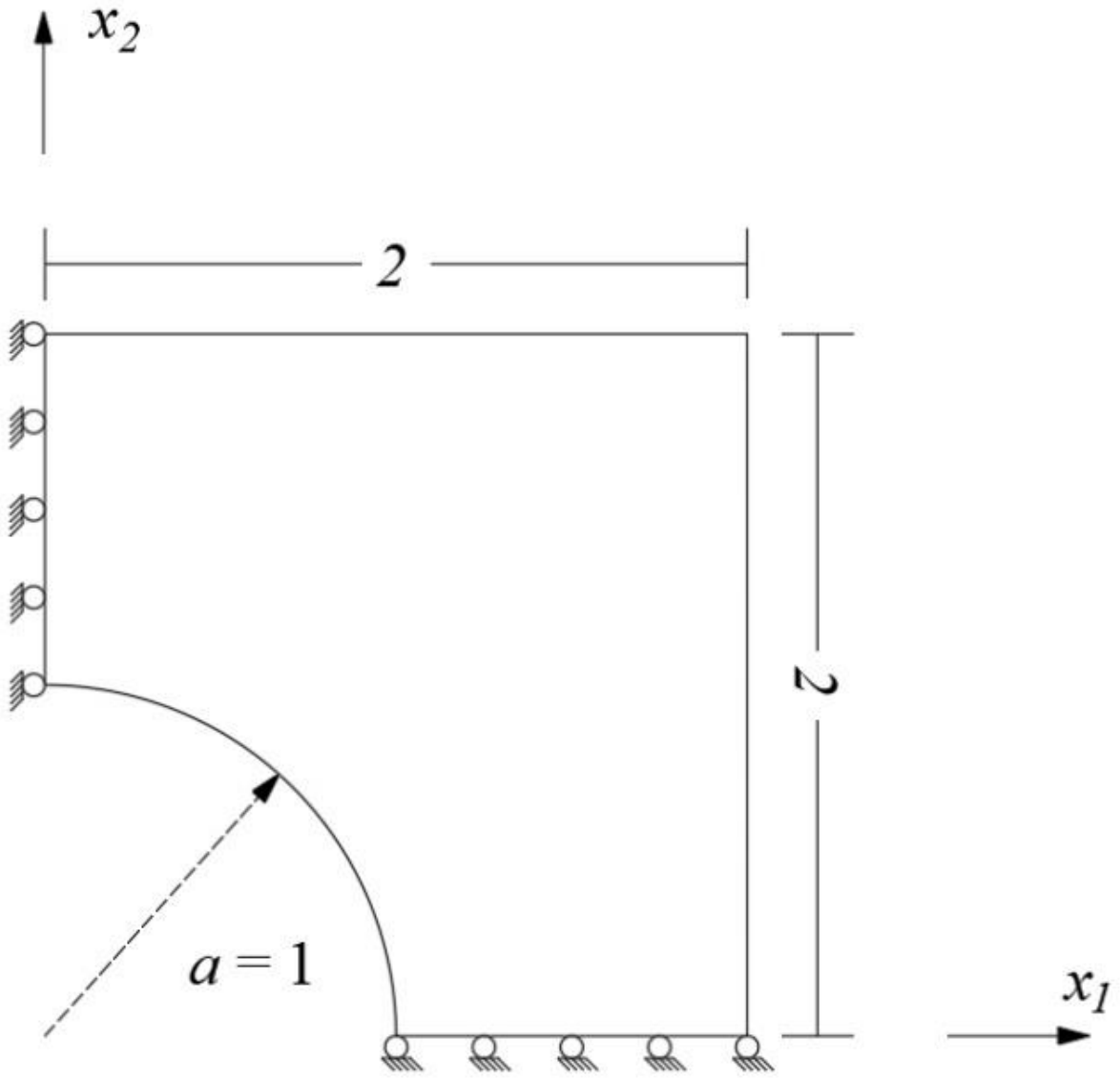

(b)

Figure 21(a). An infinite plate with a circular hole under remote tension

(b). The simplified model after considering symmetry

The problem is solved considering a Plane Stress condition with $E = 1$ and $v = 0.3$. The penalty coefficient $\eta$ is set to be $E$. There are 805 Points distributed randomly in the domain (shown in Figure 22).

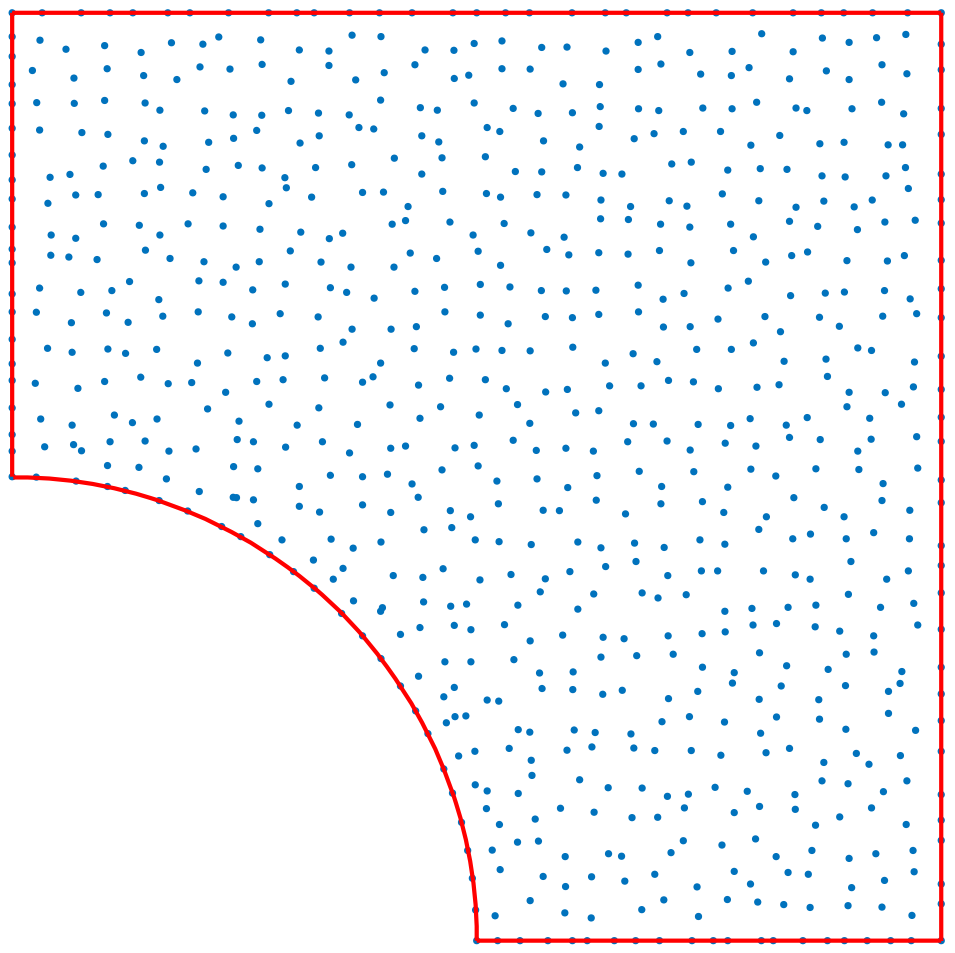

Figure 22. The Random distribution of Points

The numerical solution of $u_1$ and the corresponding error compared to the exact solution are given in Figure 23(a) and (b), respectively. Relative errors are $r_u = 3.83 \times 10^{-4}$ and $r_E = 1.84 \times 10^{-2}$.

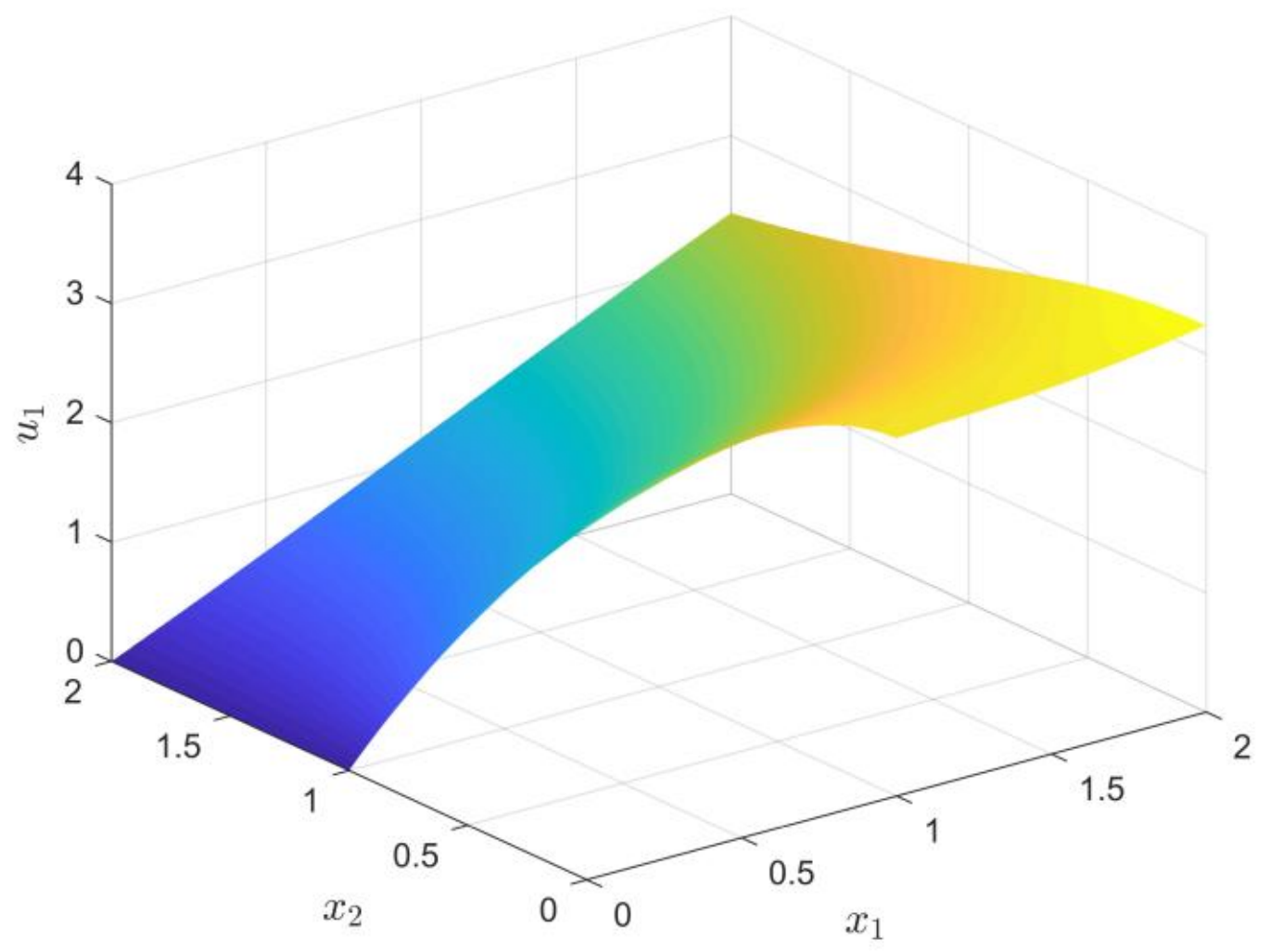

(a)

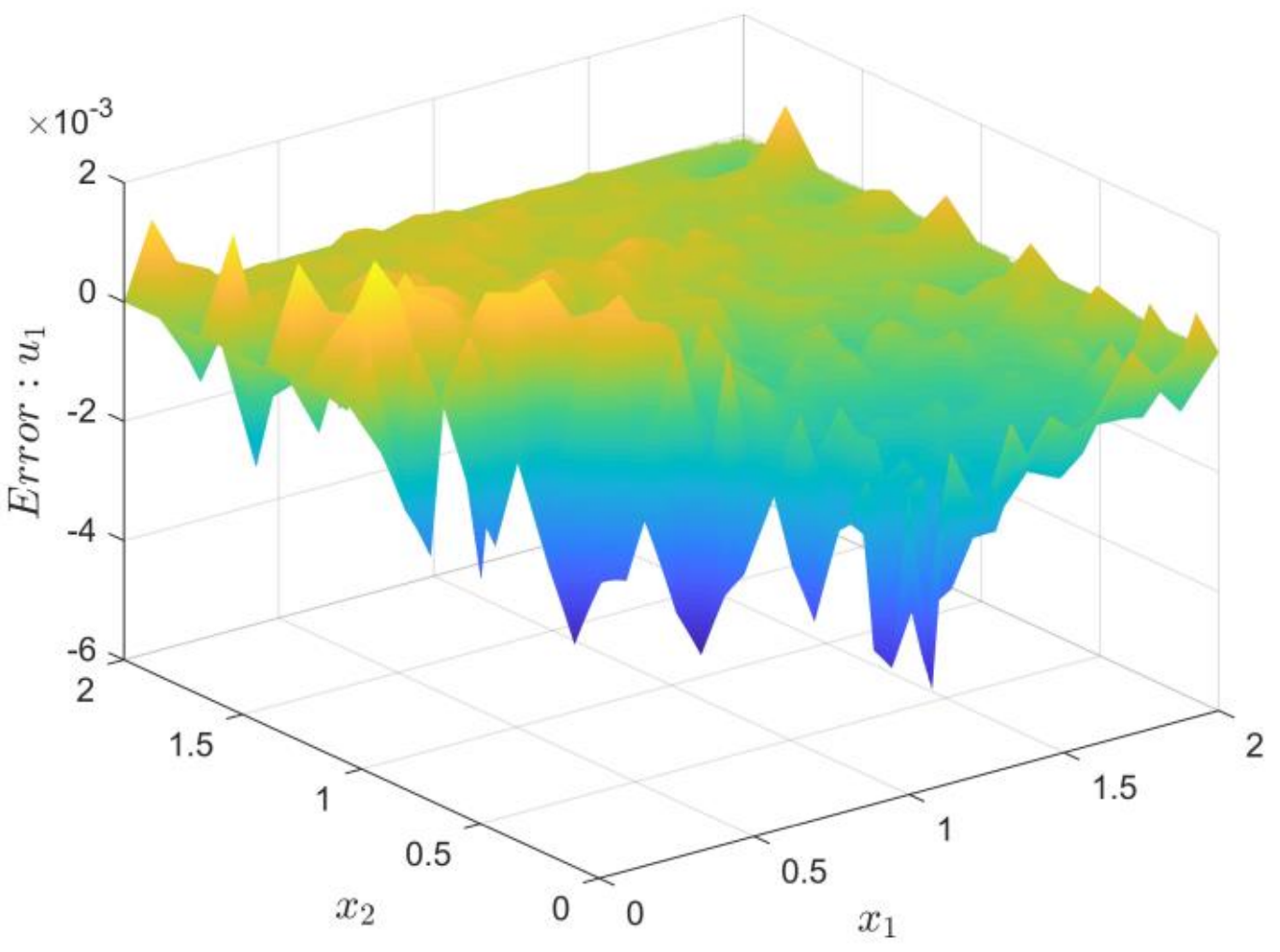

(b)

Figure 23(a). The numerical solution of $u_1$

(b). The error of $u_1$

Additionally, the numerical solutions of $u_2$ and $\sigma_{11}$ along $x_1 = 0$, as compared to the exact solutions are present in Figure 24(a) and (b), respectively. We can see that the FPM gives satisfactory solutions to this stress concentration problem.

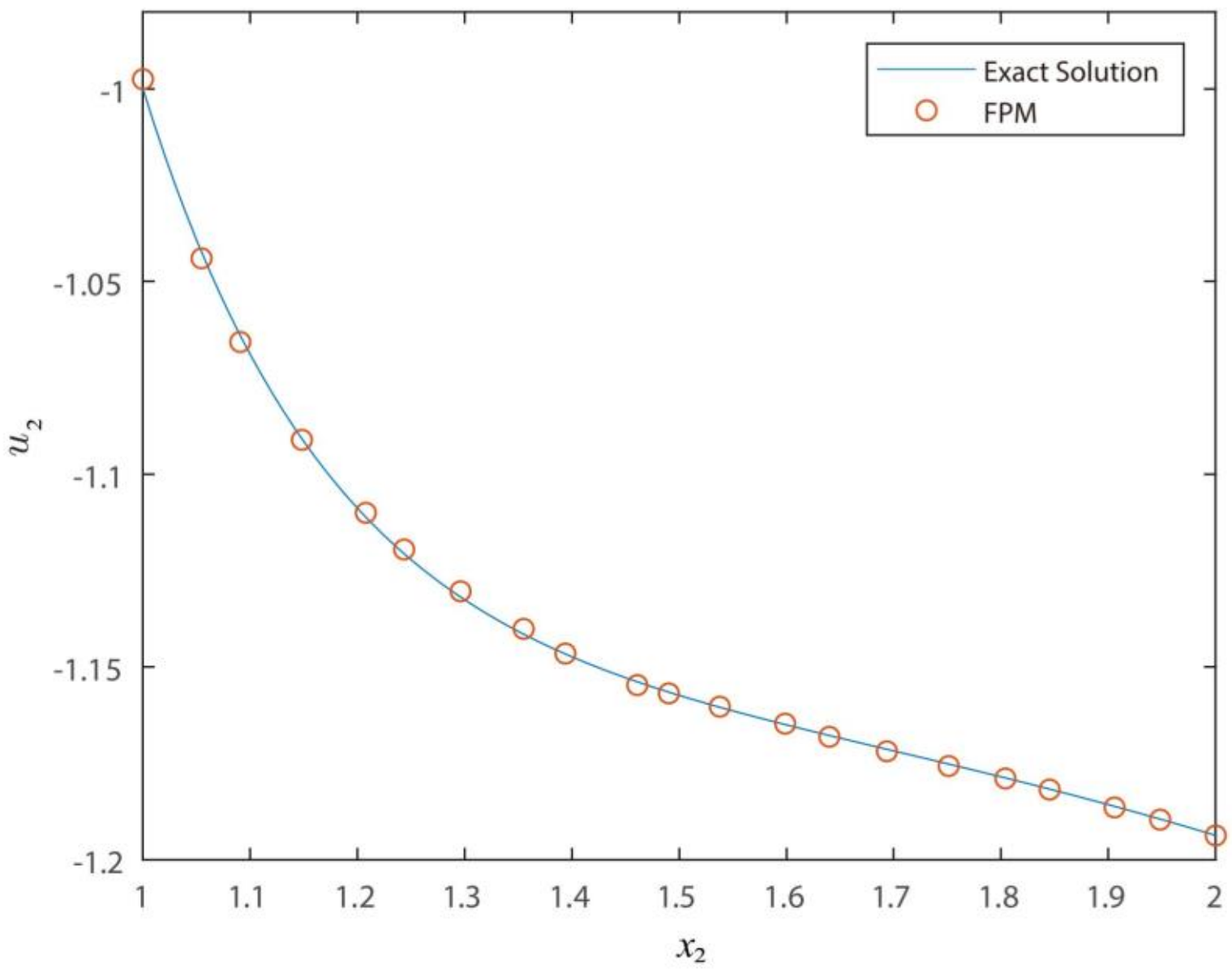

(a)

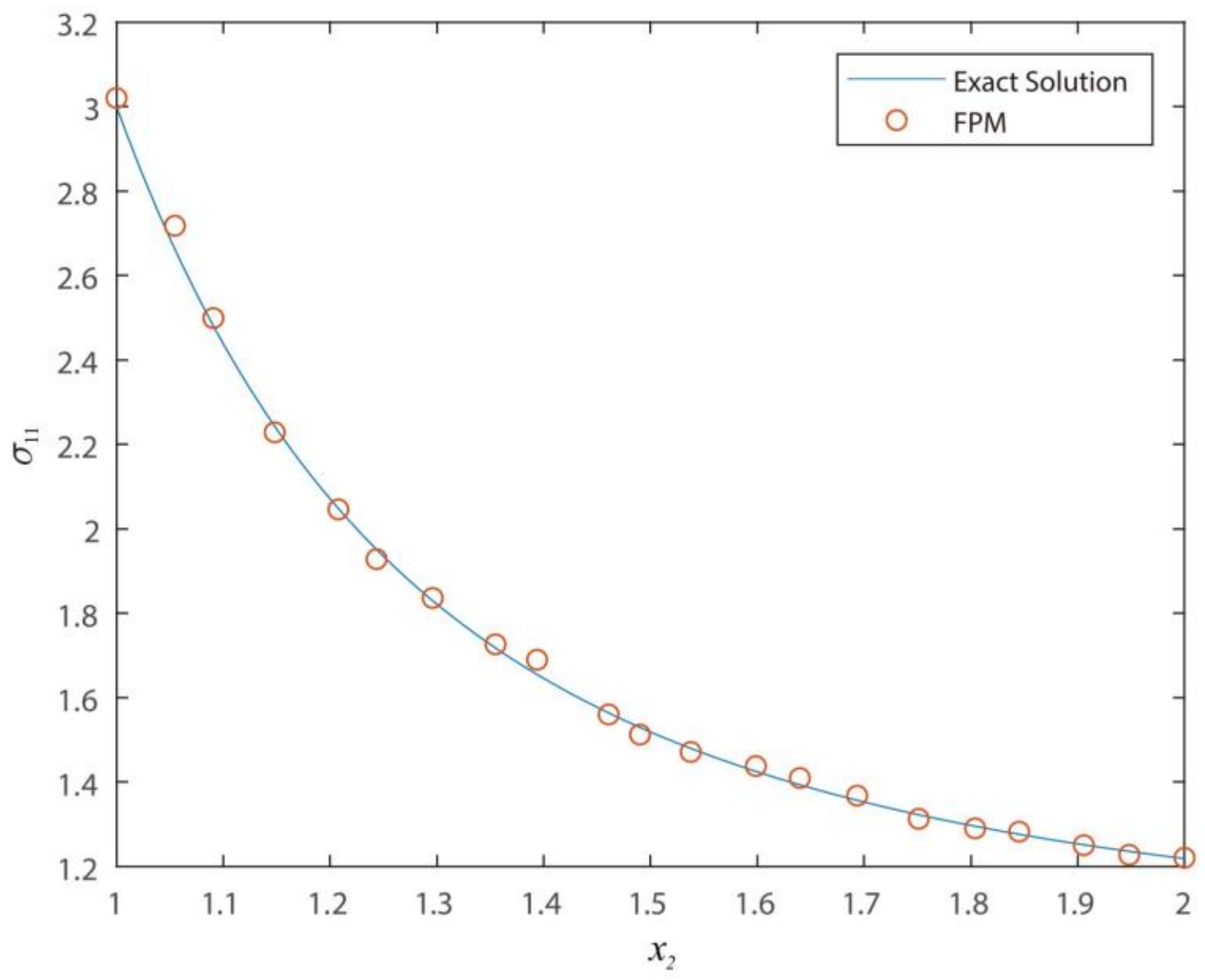

(b)

Figure 24(a). The numerical solution of $u_2$ along $x_1 = 0$

(b). The numerical solution of $\sigma_{11}$ along $x_1 = 0$

In order to demonstrate the influence of mesh on FPM, stress concentration factors are recalculated with 43, 217, 793 Points uniformly scattered in the domain, when the penalty parameter is set to be 1. Results are given in Table IV. It can be seen more Points lead to more accurate solutions. And in order to demonstrate the influence of the penalty parameter, stress concentration factors with different penalty parameters when the number of Points is equal to 793 are given in Table V. We can see that the results remain stable when the penalty parameter varies from 0.01 to 100.

Table IV. Stress Concentration Factors with different mesh

| The number of Points | 43 | 217 | 793 |
|---|---|---|---|
| Stress Concentration Factor | 2.833 | 2.967 | 2.974 |

Table V. Stress Concentration Factors with different penalty parameters

| Penalty Parameter | 0.01 | 0.1 | 1 | 10 | 100 |
|---|---|---|---|---|---|
| Stress Concentration Factor | 2.981 | 2.982 | 2.974 | 2.995 | 2.980 |

4.6 An Infinite Plate with A Pre-Existing Mode-I Crack

In this subsection, we apply the FPM to model a pre-existing crack. Specifically, a mode-I crack problem is considered. The analytical displacement and stress fields near the crack tip for a mode-I crack are given in Eq. (4.8) and Eq. (4.9), respectively [15], where $(r,\theta)$ are polar coordinates measured from the crack tip and $K_I$ is the mode-I stress intensity factor.

$$\begin{Bmatrix} u_1 \\ u_2 \end{Bmatrix} = \frac{K_I}{2\mu}\sqrt{\frac{r}{2\pi}} \begin{Bmatrix} \cos\frac{\theta}{2}\left[\kappa - 1 + 2\sin^2\frac{\theta}{2}\right] \\ \sin\frac{\theta}{2}\left[\kappa + 1 - 2\cos^2\frac{\theta}{2}\right] \end{Bmatrix} \tag{4.8}$$

$$\text{where } \kappa = \frac{3-\bar{v}}{1+\bar{v}} \quad \mu = \frac{\bar{E}}{2(1+\bar{v})}$$

$$\begin{Bmatrix} \sigma_{11} \\ \sigma_{22} \\ \sigma_{12} \end{Bmatrix} = \frac{K_I}{\sqrt{2\pi r}}\cos\frac{\theta}{2} \begin{Bmatrix} 1 - \sin\frac{\theta}{2}\sin\frac{3\theta}{2} \\ 1 + \sin\frac{\theta}{2}\sin\frac{3\theta}{2} \\ \sin\frac{\theta}{2}\cos\frac{3\theta}{2} \end{Bmatrix} \tag{4.9}$$

As shown in Figure 25, a single edge-cracked square plate with width = $b$ and crack length = $a$ is studied. Displacement boundary conditions are imposed on its four sides according to Eq. (4.8) with $K_I$ prescribed as 1. This problem is analyzed considering a Plane Stress condition. We prescribe that $b = 2a = 10$, $E = 1$, $v = 0.3$. The penalty coefficient $\eta$ is set to be equal to $E$.

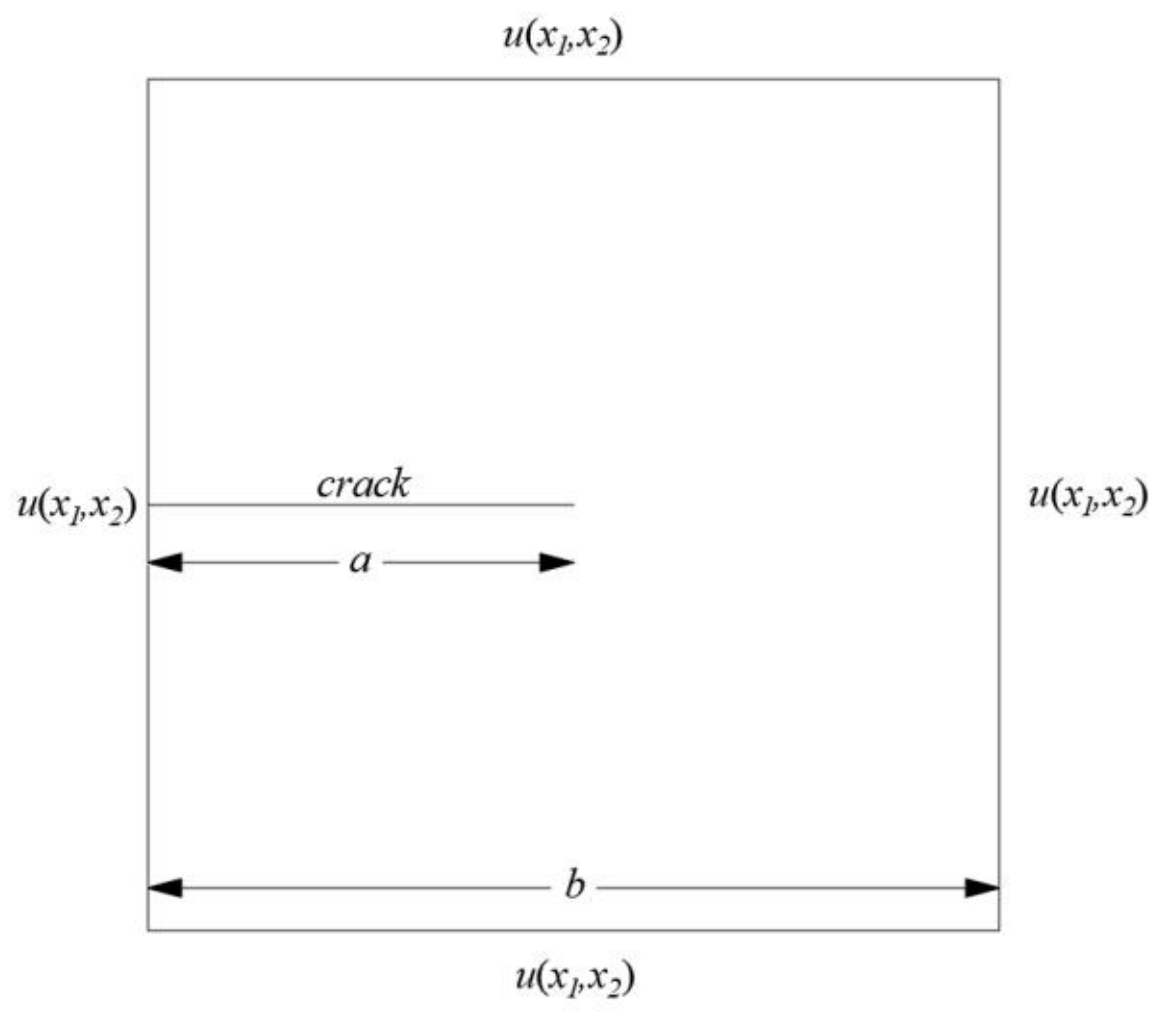

Figure 25. An edge-cracked square plate

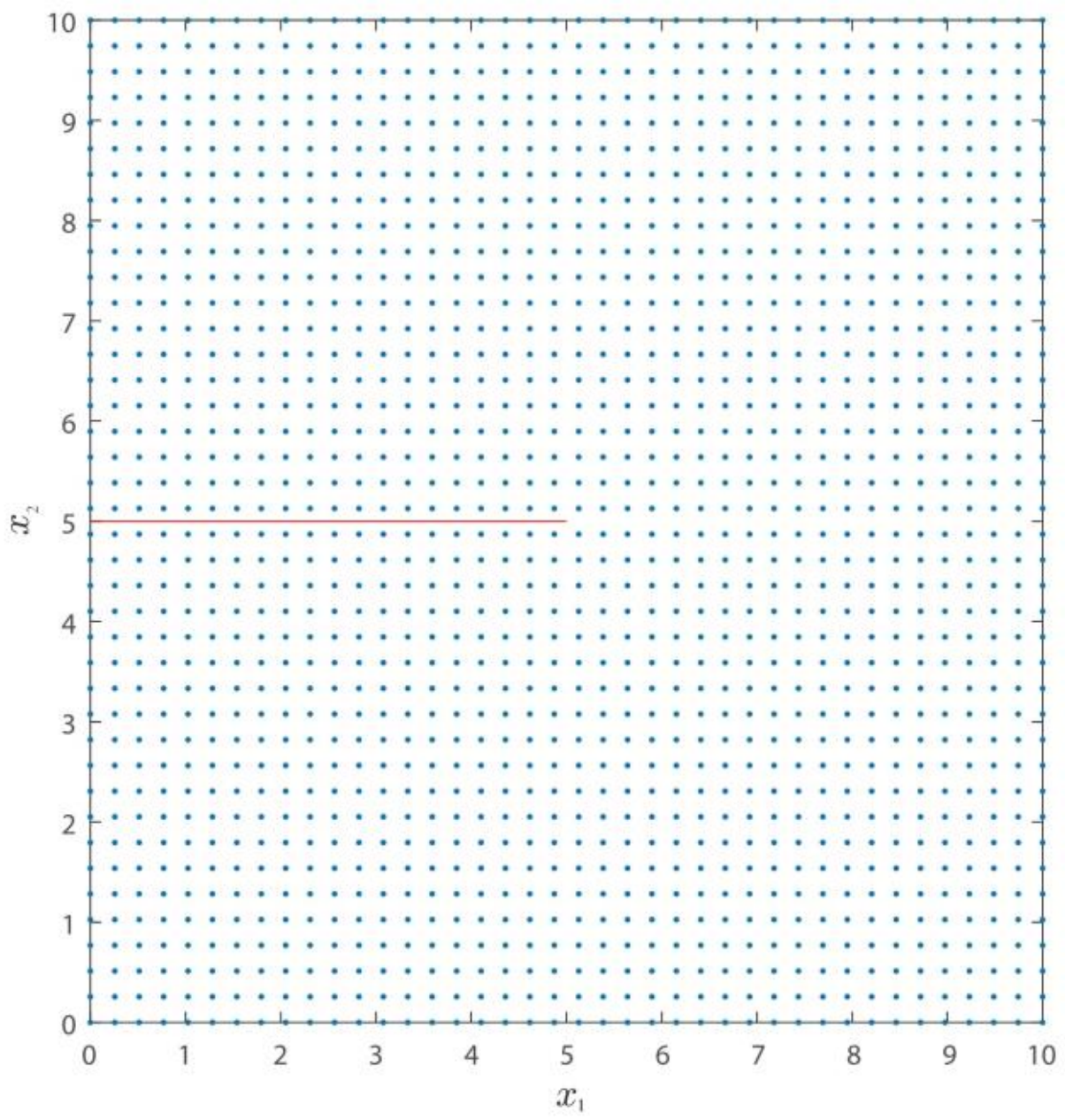

Figure 26. The regular distribution of Points

40×40 Points are scattered regularly in the domain (shown in Figure 26). The numerical solutions of $u_2$ and $\sigma_{22}$ along $x_1 = 5 + 0.5h$ ($h$ is the distance between 2 neighboring Points) as compared to exact solutions are demonstrated in Figure 27(a) and (b), respectively. It can be seen that the numerical solutions obtained by the FPM can achieve excellent accuracy.

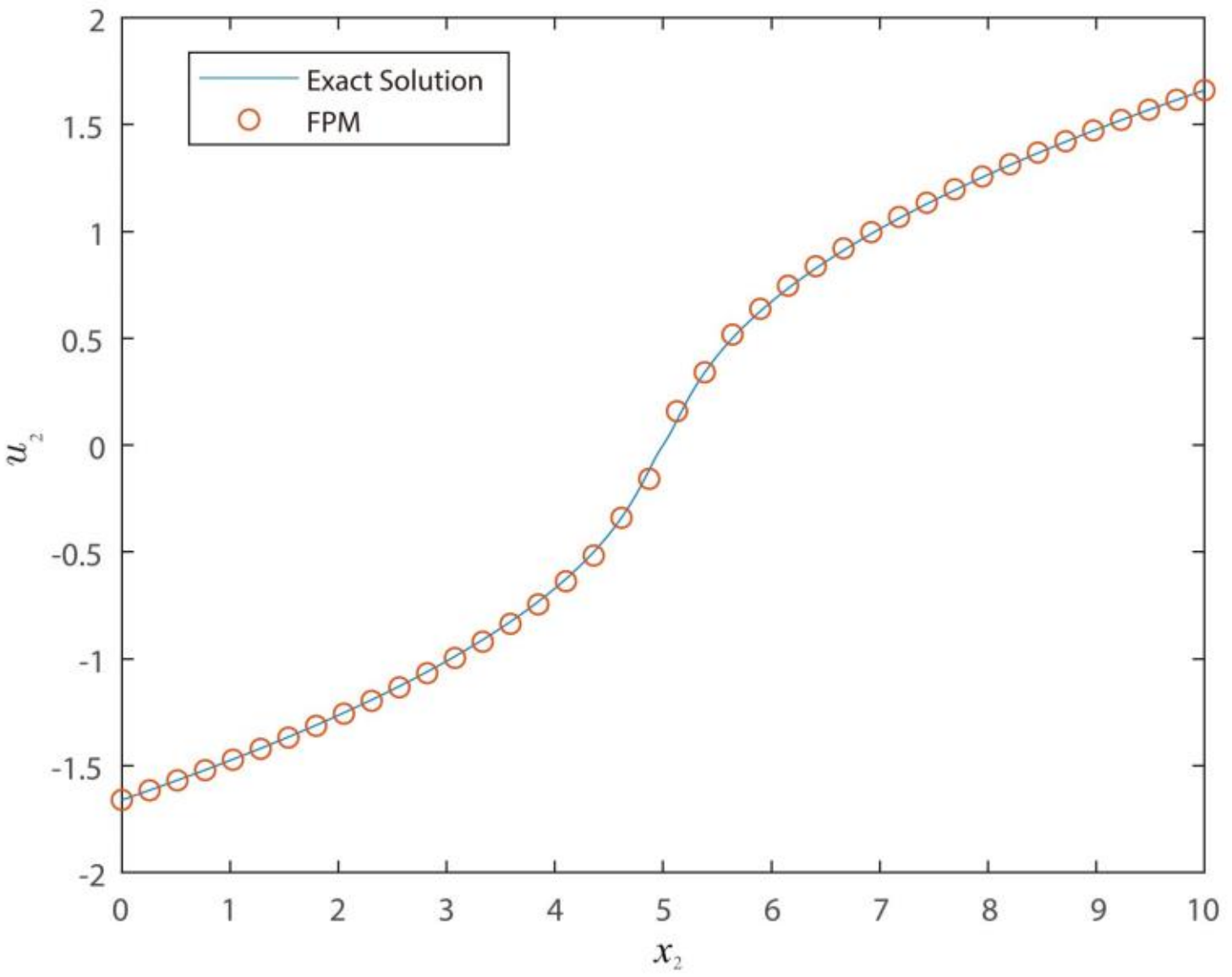

(a)

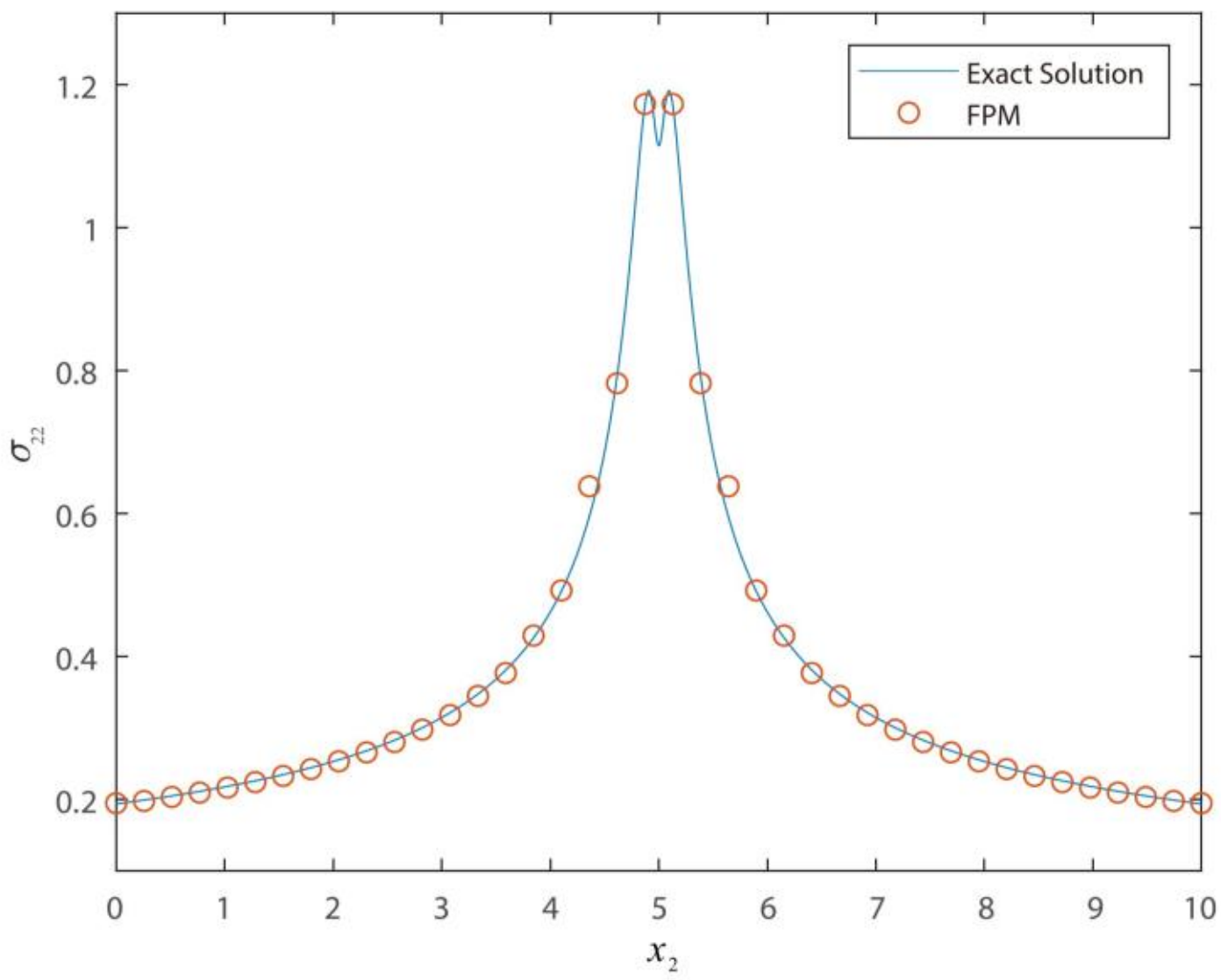

(b)

Figure 27(a). The numerical solution of $u_2$ along $x_1 = 5 + 0.5h$

(b). The numerical solution of $\sigma_{22}$ along $x_1 = 5 + 0.5h$

In addition, we further back-calculate the mode-I stress intensity factor $K_I$ from the numerical solution by calculating the *J*-integral for the computed numerical solution. A 5*h* × 9*h* rectangular contour is used to calculate the *J*-integral. The solution of $K_I$ by computing the *J*-integral is shown in Table VI, which is very accurate as compared with the prescribed exact solution.

Table VI. Stress intensity factor for the pre-existing edge crack problem

| $K_I$ (*J*-integral) | $K_{\text{exact}}$ | *Error* |
|---|---|---|
| 0.9917 | 1 | 0.0083 |

4.7 An Edge-Cracked Plate with Mixed-Mode Loading

In this subsection, an edge-cracked rectangular plate with mixed-mode loading is studied. As shown in Figure 28, the plate with length $L$ = 16, width $W$ = 7, crack length $a$ = 3.5 is fixed at the bottom. Uniformly distributed tangential tractions ($t$ = 1) are applied on its top edge.

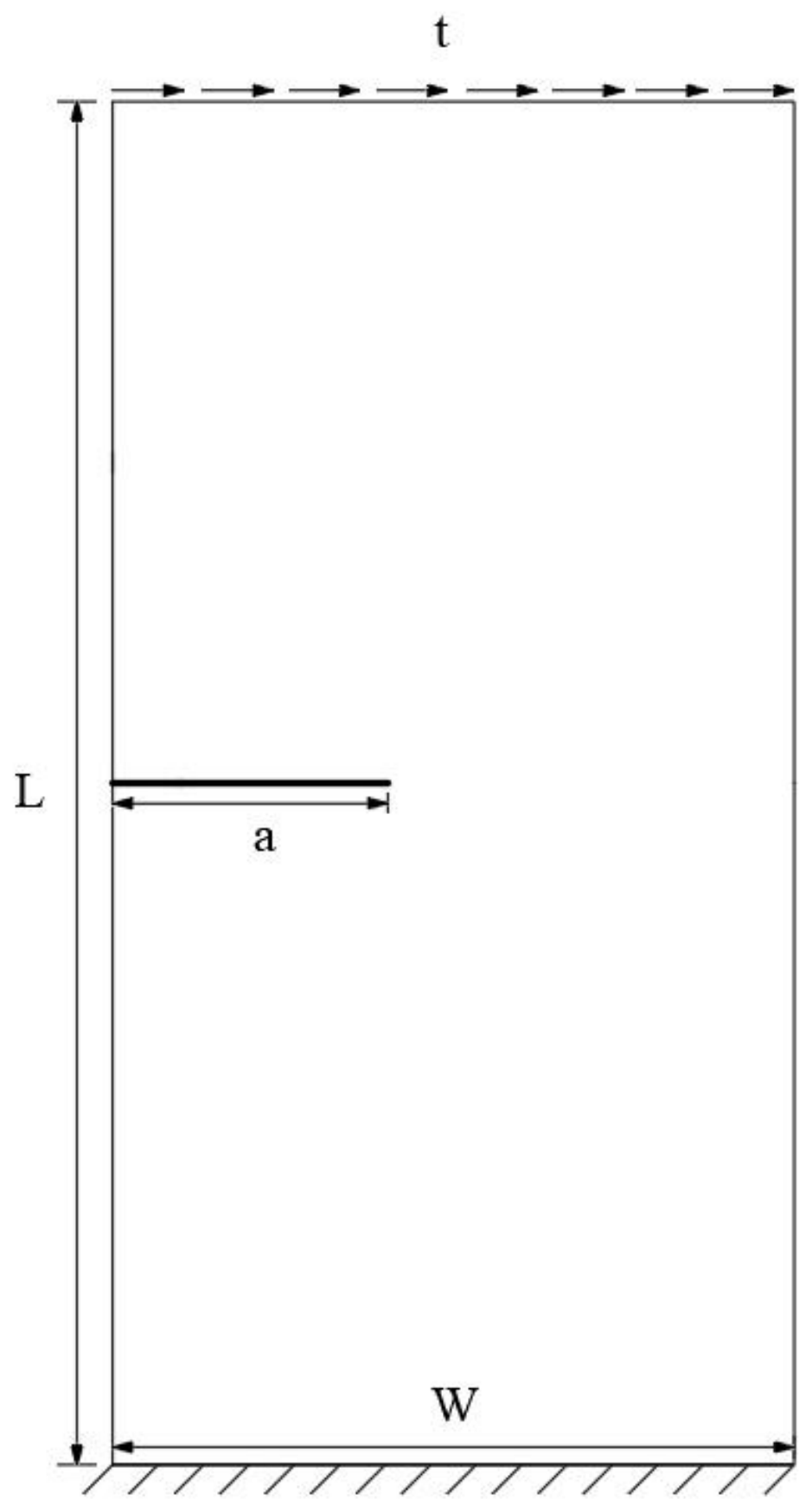

Figure 28. An edge-cracked plate under mixed-mode loading

36×80 Points are scattered regularly in the plate during the simulation. This problem is analyzed considering a Plane Stress condition. We prescribe that the Young's model $E$ = 1, Poisson's ratio $v$ = 0.3 and the penalty coefficient $\eta$ is set to be equal to

$E$.

The mode-I and mode-II stress intensity factors (SIFs) $K_I$ and $K_{II}$ are calculated separately by using the interaction integral method[16]. The interaction integral is defined as

$$M^{(1,2)} = \int_\Gamma W^{(1,2)} n_1 d\Gamma - \int_\Gamma \left[ \sigma_{ij}^{(1)} \frac{\partial u_i^{(2)}}{\partial x_1} + \sigma_{ij}^{(2)} \frac{\partial u_i^{(1)}}{\partial x_1} \right] n_j d\Gamma \tag{4.10}$$

where

$$W^{(1,2)} = \frac{1}{2}\left( \sigma_{ij}^{(1)} \varepsilon_{ij}^{(2)} + \sigma_{ij}^{(2)} \varepsilon_{ij}^{(1)} \right)$$

The superscripts "1" and "2" stand for the yet unsolved real state and an auxiliary state, respectively. $\Gamma$ is the integral contour. According to [16], the following relation between the interaction integral and SIFs exists.

$$M^{(1,2)} = \frac{2}{E}\left( K_I^{(1)} K_I^{(2)} + K_{II}^{(1)} K_{II}^{(2)} \right) \tag{4.11}$$

For the auxiliary state, if $K_I^{(2)} = 1$ and $K_{II}^{(2)} = 0$, then $M^{(1,2)} = \frac{2}{E} K_I^{(1)}$. Similarly, when $K_I^{(2)} = 0$ and $K_{II}^{(2)} = 1$, then $M^{(1,2)} = \frac{2}{E} K_{II}^{(1)}$. Therefore, SIFs can be computed separately by using the interaction integral.

The contour for the interaction integral is defined as a rectangle with length *l* and width *w*. According to [17], reference values are set as $K_I = 34.0$ and $K_{II} = 4.55$. Numerical solutions calculated by the FPM shown in the Table VII agree well with the reference values.

Table VII. Numerical solutions of SIFs by the FPM

| $2l \times 2w$ | $K_I$ | *ratio* | $K_{II}$ | *ratio* |
|---|---|---|---|---|
| $1.2 \times 3.2$ | 33.57 | 0.987 | 4.62 | 1.015 |
| $2.8 \times 6.4$ | 33.80 | 0.994 | 4.58 | 1.007 |
| $6.4 \times 6.4$ | 33.83 | 0.995 | 4.55 | 1.000 |

In order to study the influence of the mesh on computed SIFs, the same problem is solved by FPM using 36 × 80 Points and 70 × 160 Points, respectively ($2l \times 2w = 1.6 \times 3.2$). Results are shown in Table VIII. Both solutions agree well with the reference SIFs.

Table VIII. Numerical solutions of SIFs by the FPM with different meshes

| Number of Points | $K_I$ | *ratio* | $K_{II}$ | *ratio* |
|---|---|---|---|---|
| 36×80 | 34.10 | 1.003 | 4.58 | 1.007 |
| 70×160 | 34.08 | 1.002 | 4.56 | 1.002 |

4.8 Realistic Structures with Complex Geometries

In this subsection, we employ the FPM to simulate realistic structures with complex geometries. An interface program has been designed to connect the preprocessing module of ABAQUS with the FPM solver which is written by MATLAB. Elements meshed in ABAQUS will be converted into subdomains for FPM through the interface program (shown in Figure 2). For every subdomain located inside the problem domain, its interpolation Point is at the center of mass of the subdomain. As for the subdomain on the boundary, its interpolation Point sits on the edge of the subdomain, which coincides with the boundary.

4.8.1 Stress analysis of a wrench

For the model of a wrench shown in the Figure 29, its jaw is set to be fixed and uniformly distributed tractions are imposed on the end of the wrench. This structure is simulated by FPM, and for comparison purposes it is also modeled in ABAQUS by FEM. The mesh for FEM in ABAQUS is illustrated in Figure 30. The number of total elements is equal to 3001, and the element type is set to be CPS4R (Plane Stress, rectangle, reduced integration). The elements employed by ABAQUS are converted to subdomains for the FPM by the interface program.

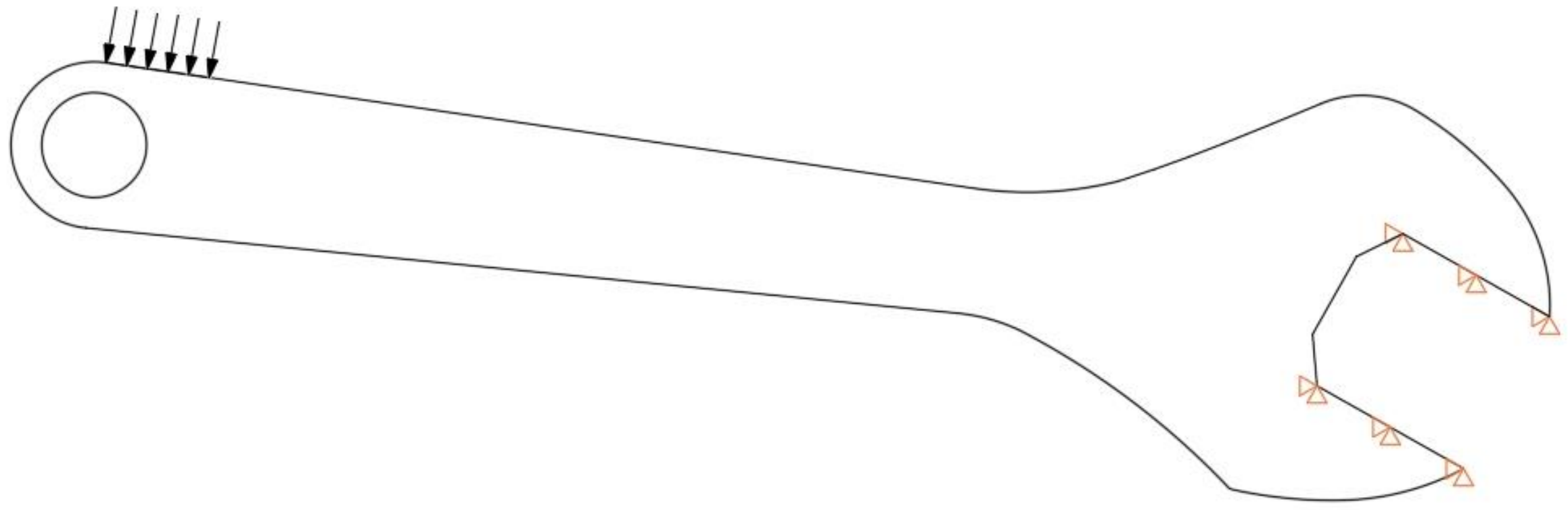

Figure 29. The wrench model

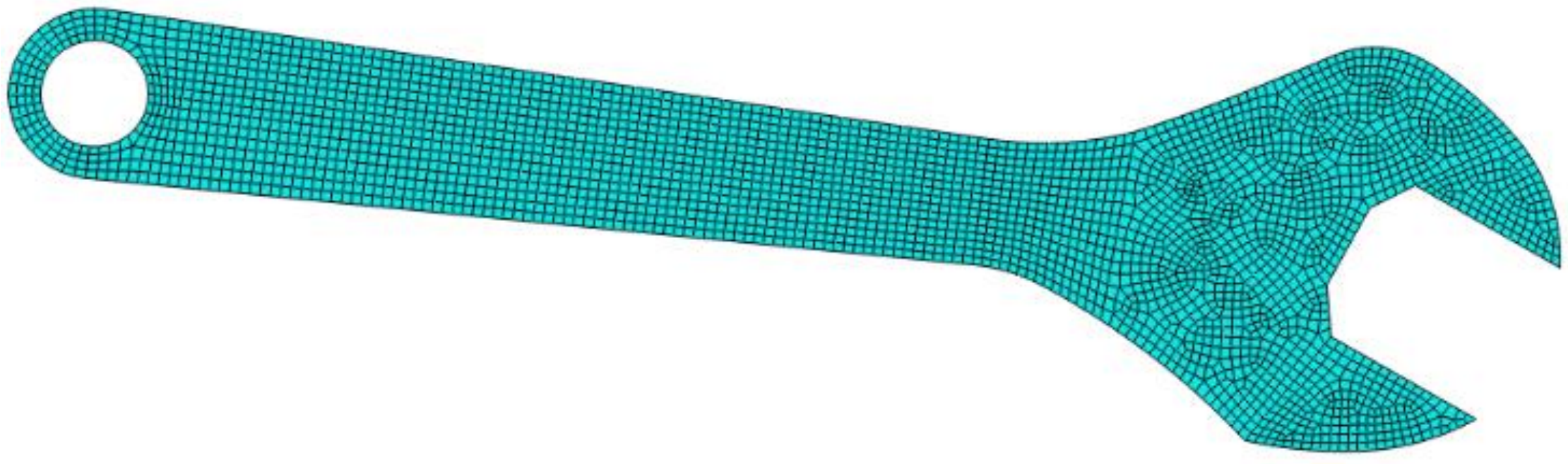

Figure 30. The mesh in ABAQUS for the wrench

The relative difference $r_E$ according to Eq. (4.1) between the results obtained by the two methods is about 7%. Specifically, numerical solutions of $\sigma_{11}$ by the FPM and the FEM are illustrated in Figure 31(a) and (b), respectively.

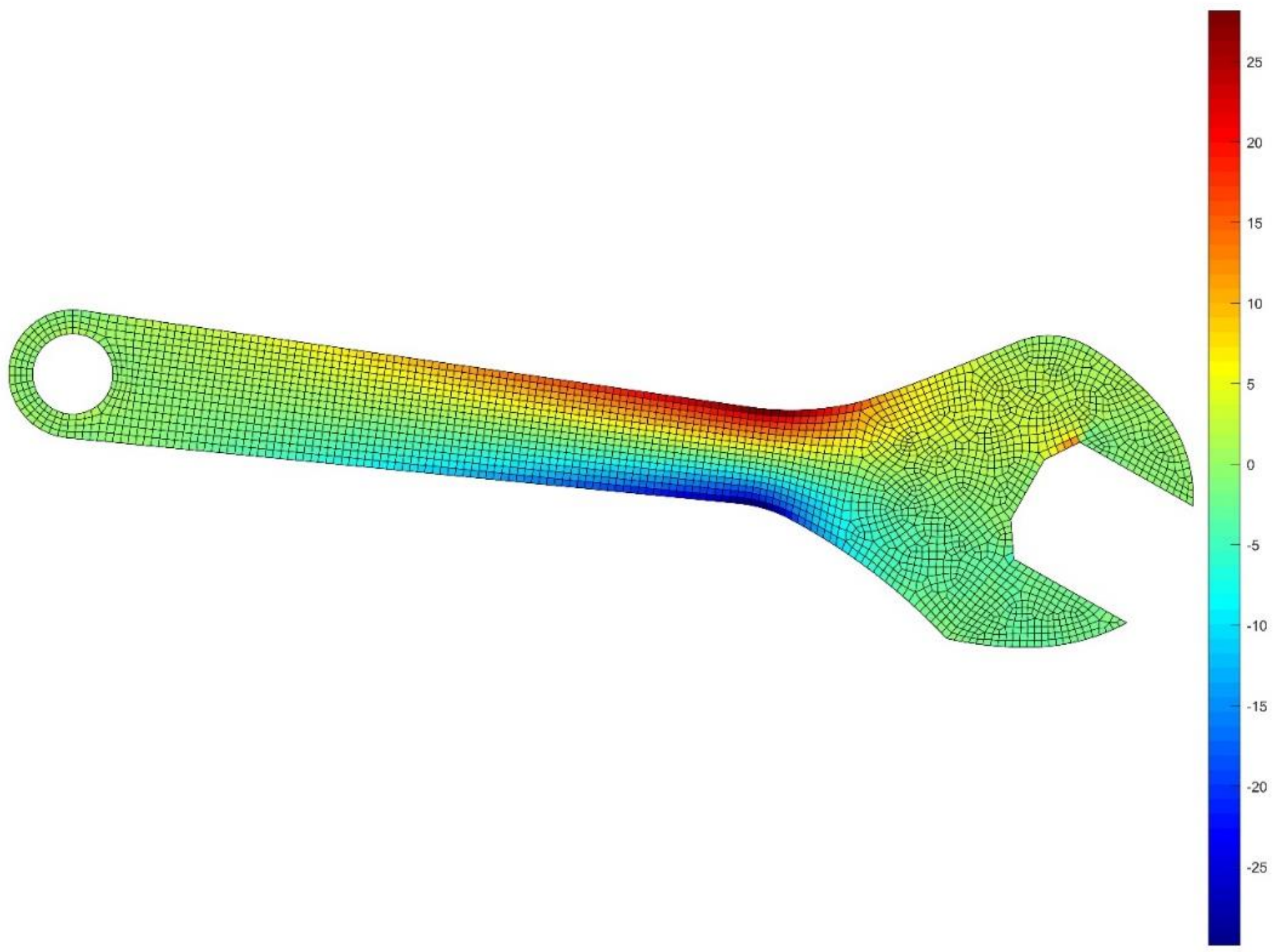

(a)

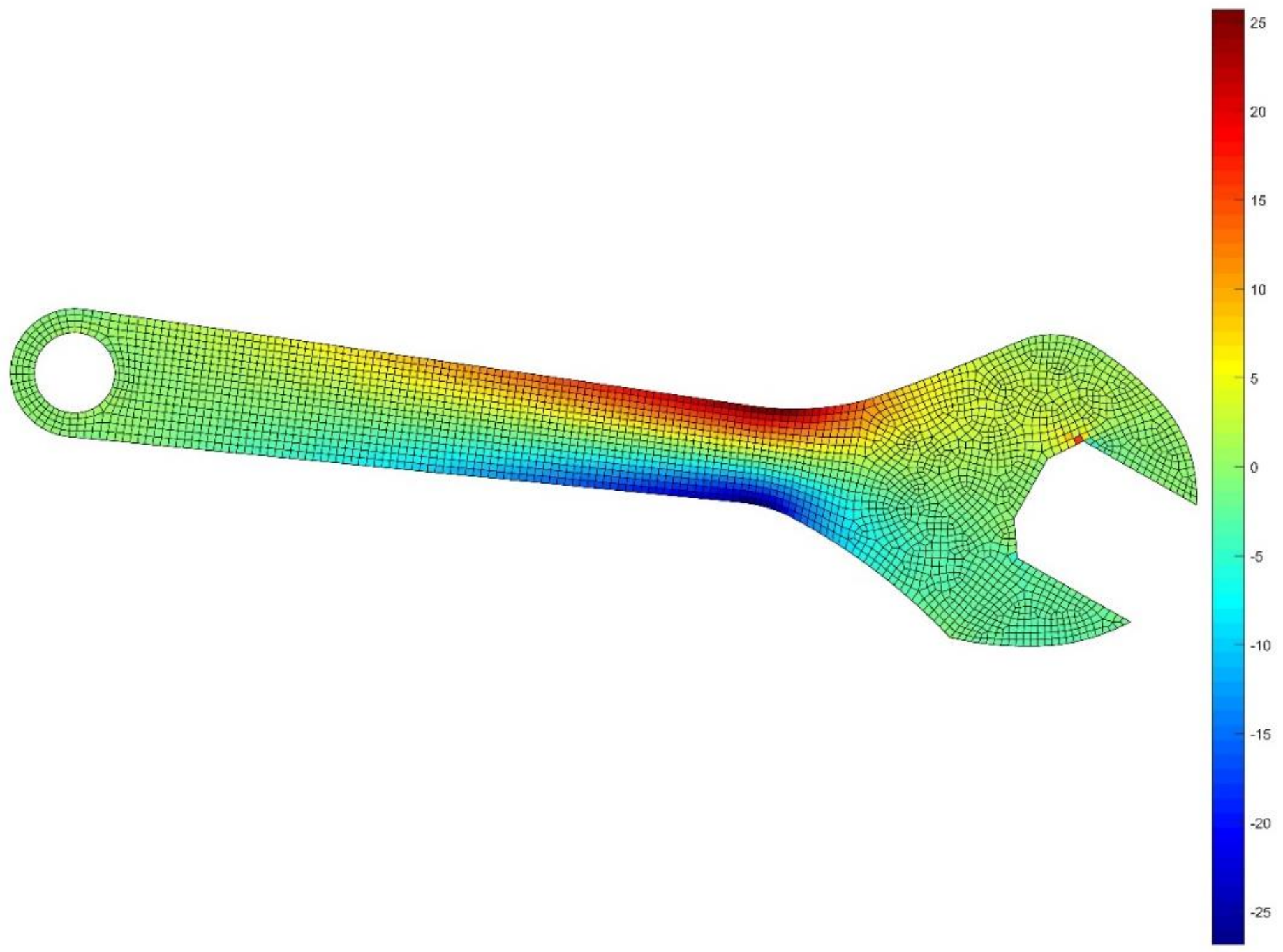

(b)

Figure 31(a). The numerical solution of $\sigma_{11}$ by using the FPM

(b). The numerical solution of $\sigma_{11}$ by the FEM using ABAQUS

4.8.2 Stress analysis of a connecting rod

In this subsection, we employ both FPM and FEM to simulate a connecting rod which is shown in the Figure 32. The internal surface of the smaller hole in the rod is fixed and half of the internal surface of the larger hole is subjected to uniform tractions. The mesh for the FEM in ABAQUS is illustrated in Figure 33. The number of total elements is 2095, and the element type in ABAQUS is also set to be CPS4R. The mesh is alternatively converted into subdomains suitable for FPM by the currently developed interface program.

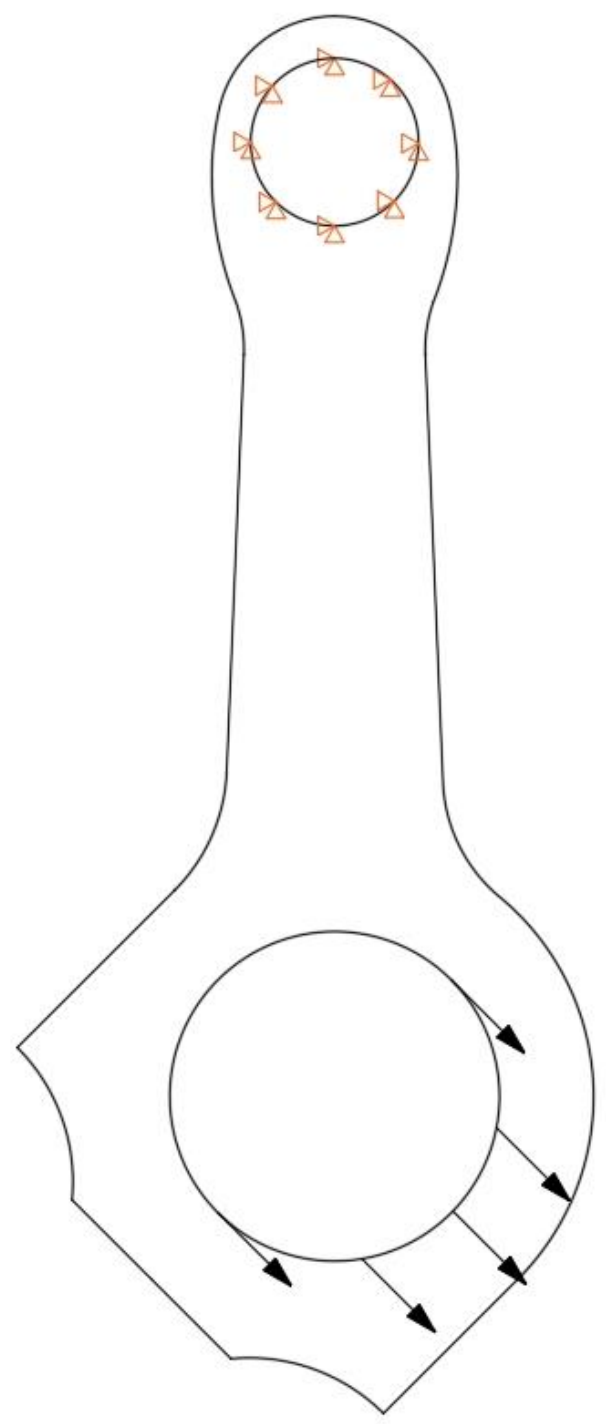

Figure 32. The connecting rod model

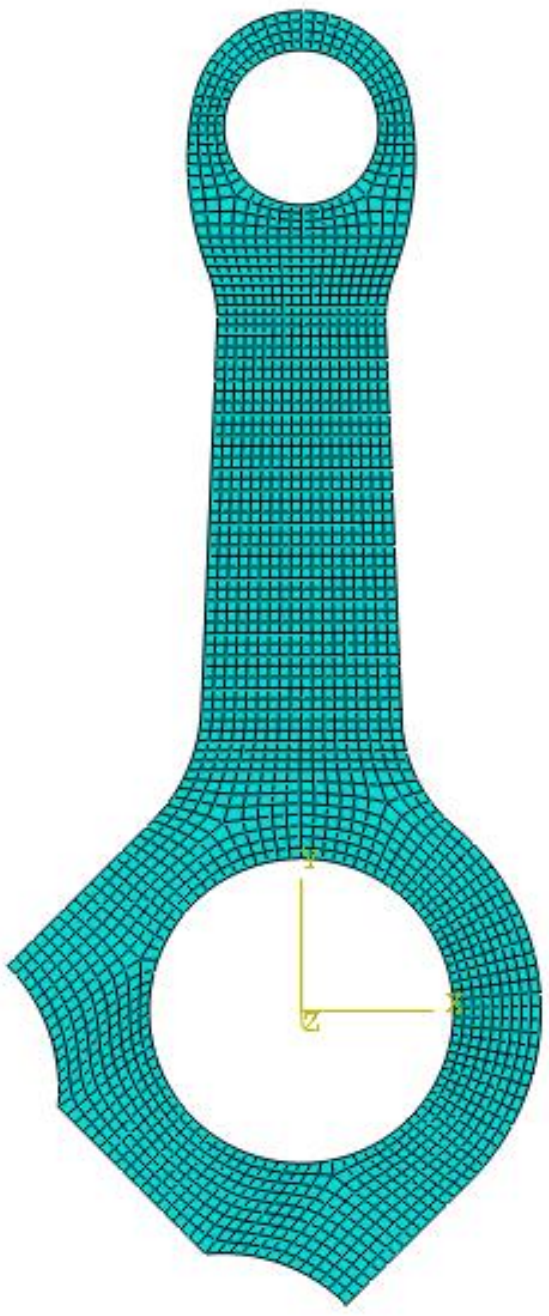

Figure 33. The mesh in ABAQUS software for the connecting rod

The relative difference $r_E$ according to Eq. (4.1) between results obtained by the two methods is about 6%, and numerical solutions of $\sigma_{22}$ by the FPM and the FEM are illustrated in Figure 34(a) and Figure 34(b), respectively.

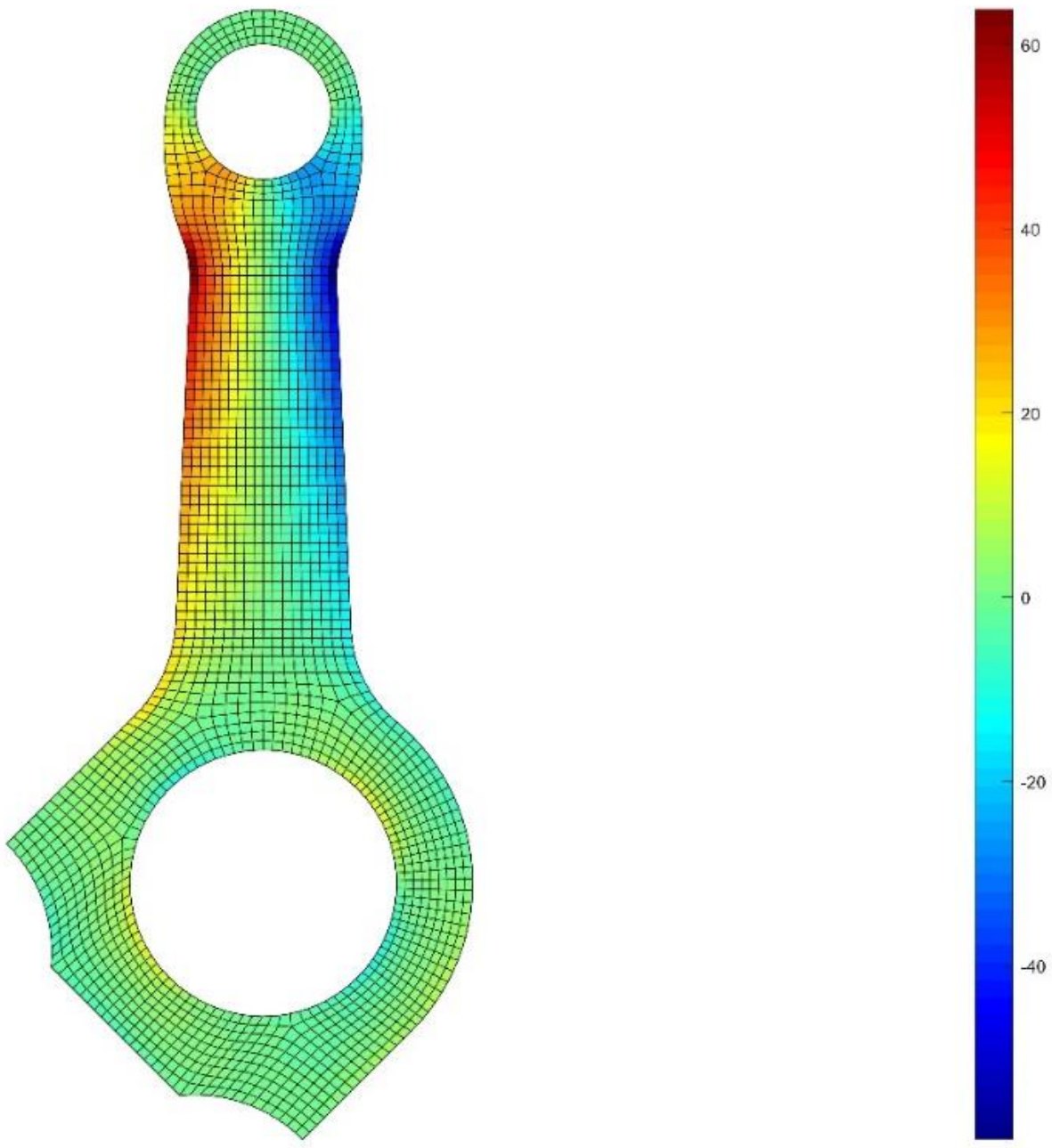

(a)

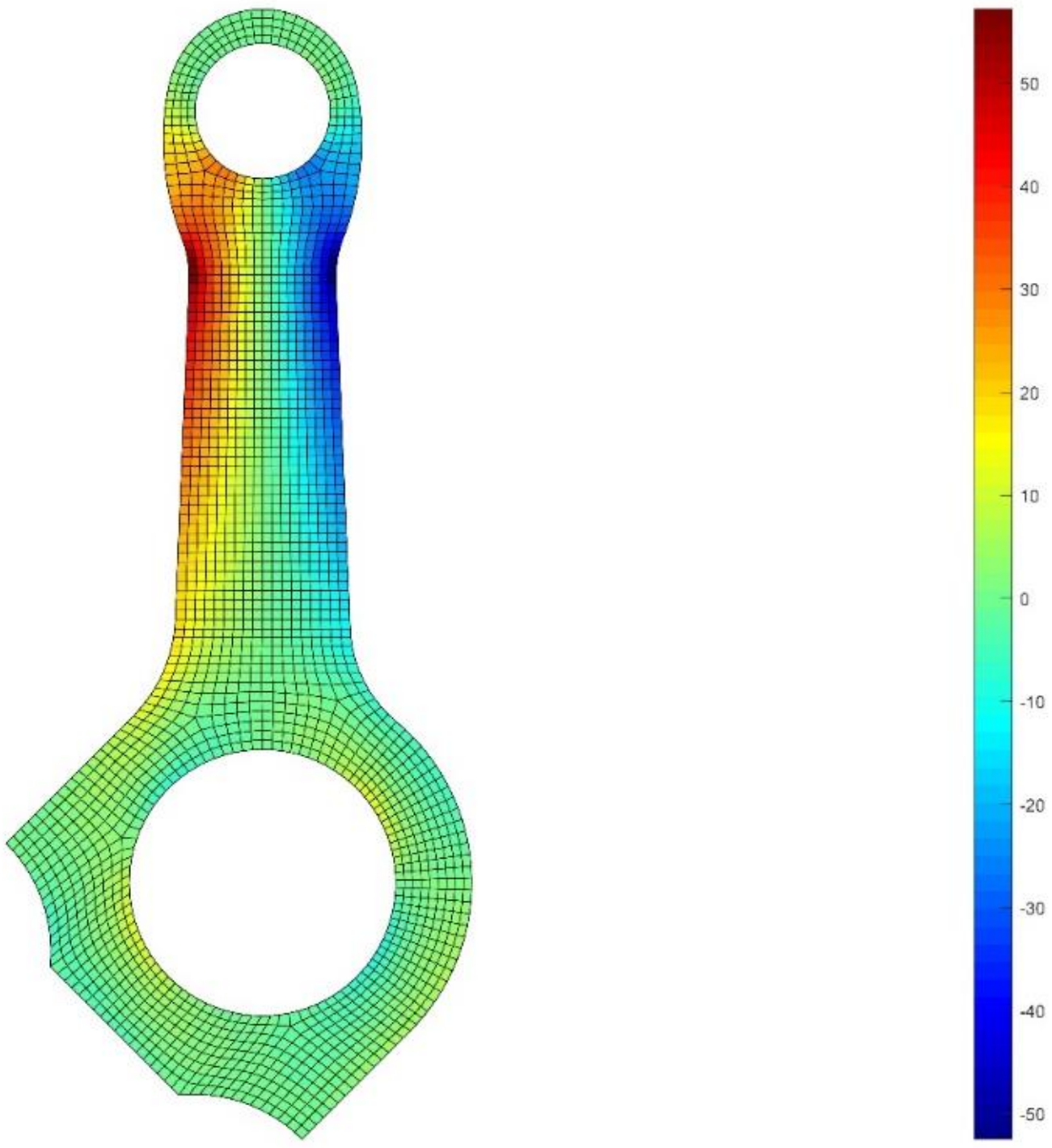

(b)

Figure 34(a). The numerical solution of $\sigma_{22}$ by the FPM

(b). The numerical solution of $\sigma_{22}$ by the FEM in ABAQUS

It takes the FPM solver, written by MATLAB, 9.57s and 7.06s to model the wrench and rod, respectively. In comparison, it takes the ABAQUS software about 15s to solve each problem. Since polynomial trial and test functions are employed in both the FPM and FEM, leading to simple Gauss quadrature for the evaluation of stiffness matrices. Then it is expected that the computational time spent by FPM and FEM with the same nodal distributions is in the same order.

### 4.9 Simulations of Crack Propagation Paths

In this subsection, the FPM is employed to simulate crack propagation paths. All the simulations are based on the Linear Elastic Fracture Mechanics, and the Maximum Hoop Stress criterion proposed by Erdogan and Sih [18] is used to predict the crack propagation paths.

#### 4.9.1 A Plate with A Pre-Existing Oblique Crack

Mageed and Pandey [19] conducted a series of uniaxial tension experiments on 2024-T3 Aluminum alloy sheets with centrally-located oblique cracks. As shown in Figure 35, the initial crack length $a$ = 20*mm*, and the length ($L$) and width ($W$) of specimens are equal to 220*mm* and 110*mm*, respectively. The crack angle $\beta$ is set to be $15°$ or $60°$. Young's modulus $E$ = 71*GPa*, and Poisson's Ratio $v$ = 0.33.

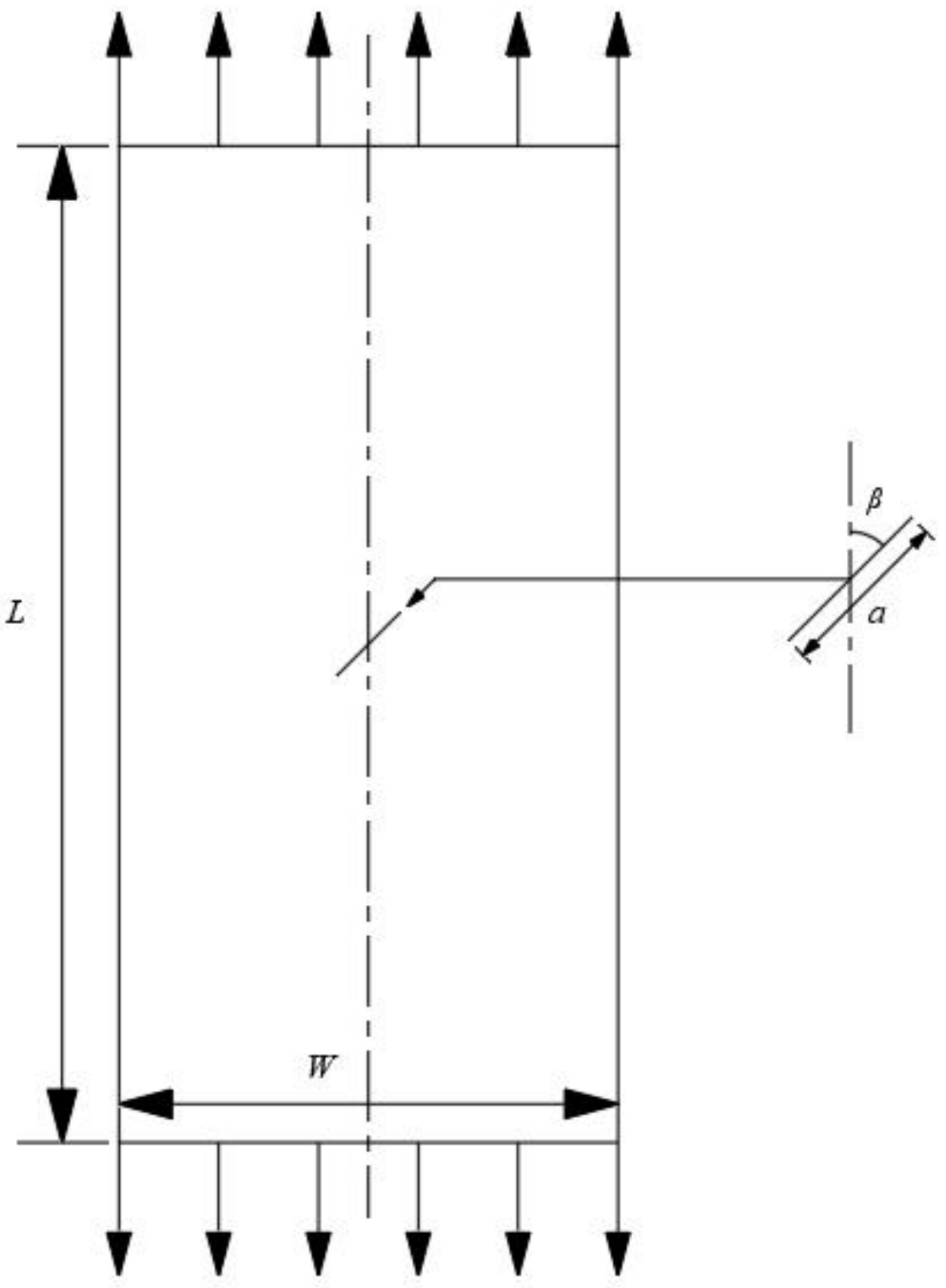

Figure 35. Specimens with pre-existing cracks

About 20,000 Points are distributed irregularly in the FPM model. As we mentioned before, the Maximum Hoop Stress criterion is employed at here to predict the crack growth paths. Because of the discontinuity of computed stresses across the internal boundary, the average of the stresses calculated in the two adjacent subdomains is set as the stress value at the internal boundary. Therefore, in each analysis step, for all the internal boundaries connected to the current crack tip, the one with the Maximum Hoop Stress will be cracked. The crack propagation paths simulated by the FPM as compared with the experiment data when $\beta$=15° and 60° are illustrated in Figure 36.

We can see that the predicted crack paths by the FPM are in good agreement with the experiment results.

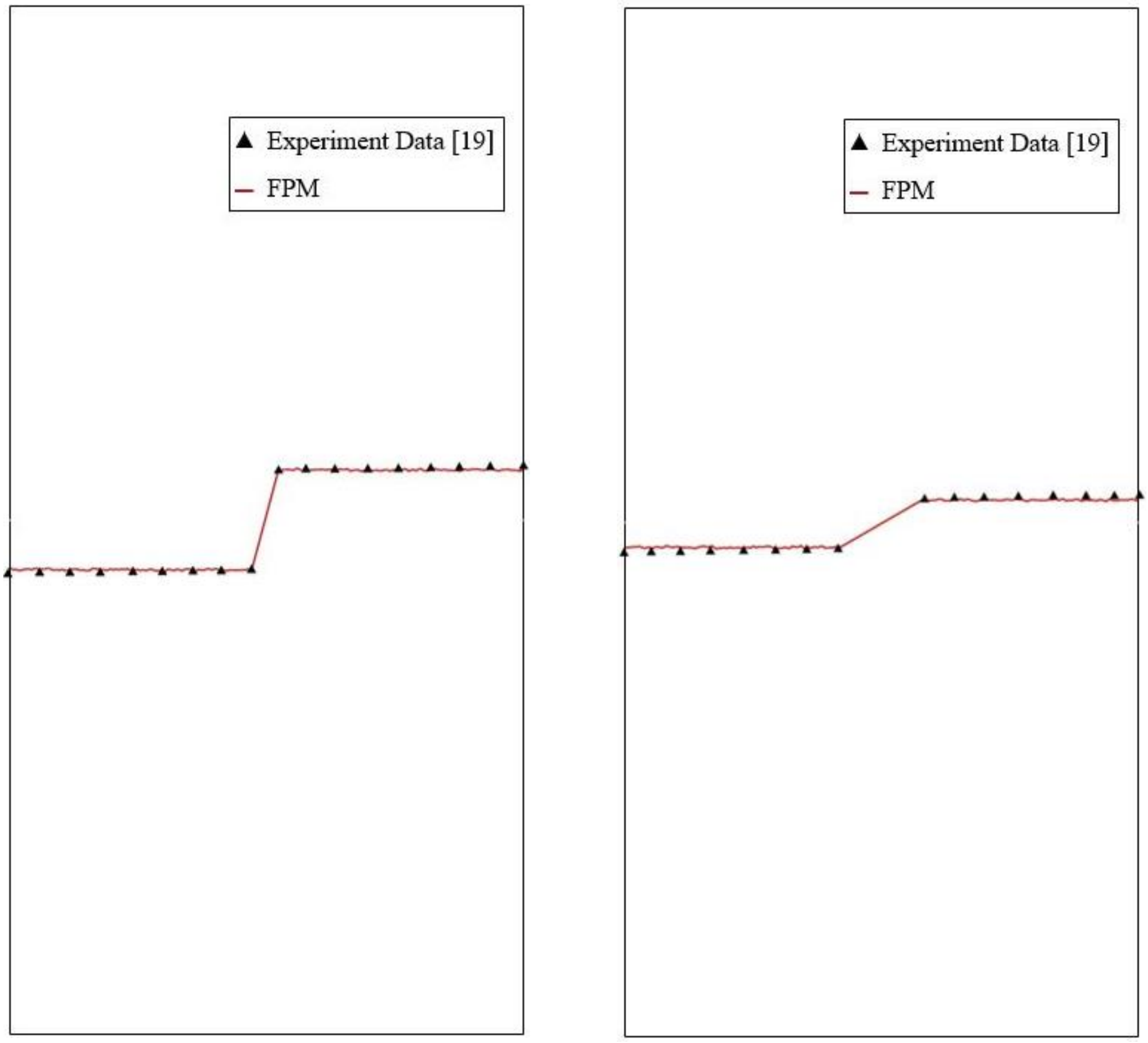

(a) (b)

Figure 36. Crack propagation paths simulated by the FPM

(a) $\beta$=15° (b) $\beta$=60°

For the case where $\beta$=15°, the crack propagation paths simulated by FPM with 7200, 12800, 20000 randomly distributed Points are shown in Figure 37. It can be seen that the crack propagation path simulated by FPM is not sensitive to Point distributions for this case.

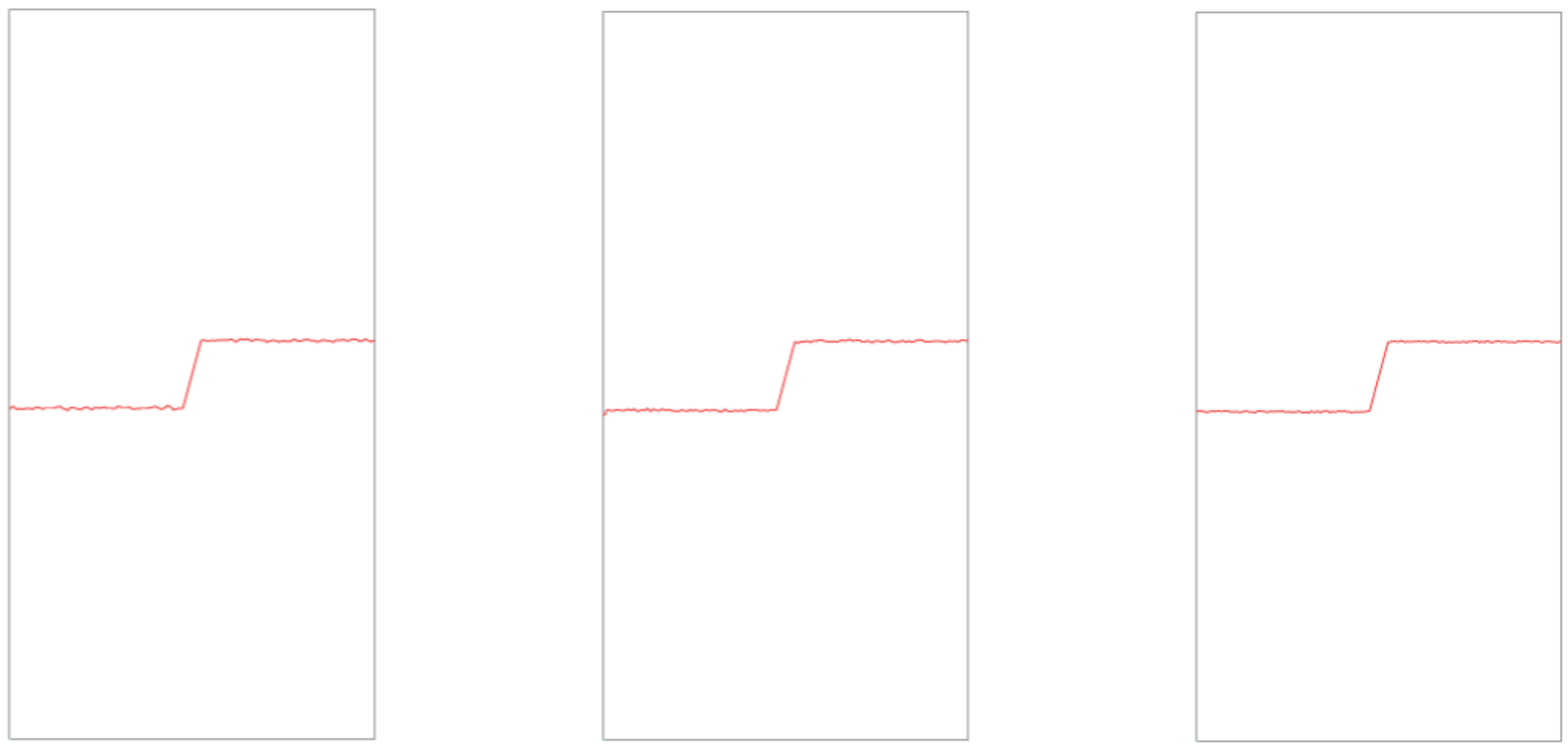

(a) (b) (c)

Figure 37. Crack propagation paths simulated by the FPM with

(a) 7200 random Points (b) 12800 random Points (c) 20000 random Points

4.9.2 A Disk with A Pre-Existing Oblique Crack

Pre-cracked disk specimens of rock-like material were experimentally tested under compressive line loading by Haeri [20]. As shown in Figure 38, The diameter of the disk $D$ = 100$mm$, and the initial crack length $a$ = 30$mm$, and the initial crack angle $\beta$=45°. Young's modulus of the rock-like material: $E$ = 15$GPa$ and Poisson's ratio $v$ = 0.21.

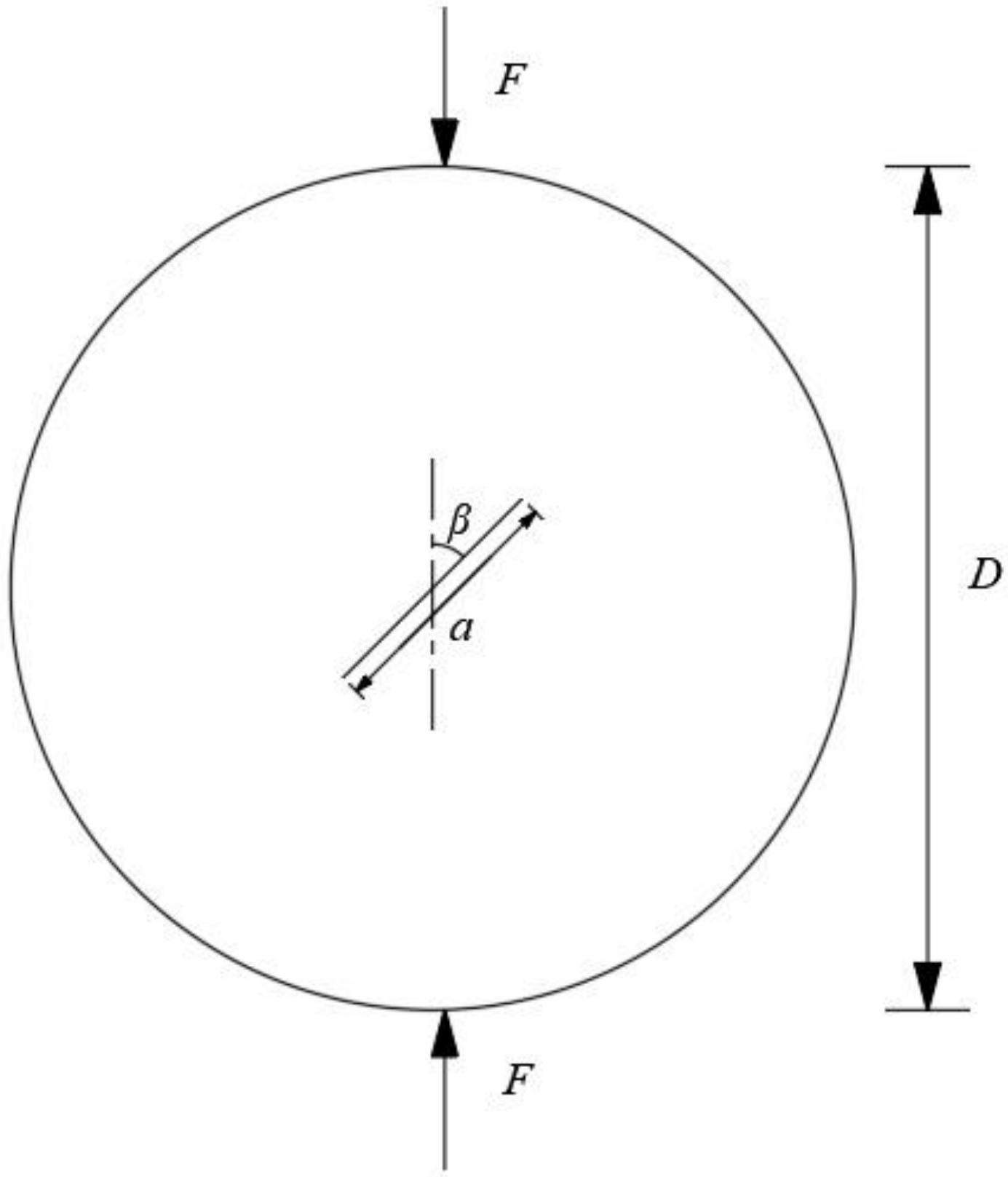

Figure 38. The pre-cracked disk specimen

The experimental results for crack propagation paths [20] are illustrated in Figure 39(a), and the numerical simulation results obtained by Haeri using the Boundary Element Method (BEM) and the Maximum Hoop Stress criterion are also given in Figure 39(b). As for the FPM, 11000 Points are randomly scattered in the structure, and in each analysis step, for all the internal boundaries connected to the current crack tip, the one with the maximum hoop stress will be cracked. The FPM simulations of the crack propagation paths are given in Figure 39(c). We can see that the crack propagation paths predicted by the FPM are in good agreement with the experiment results as well as the BEM simulations.

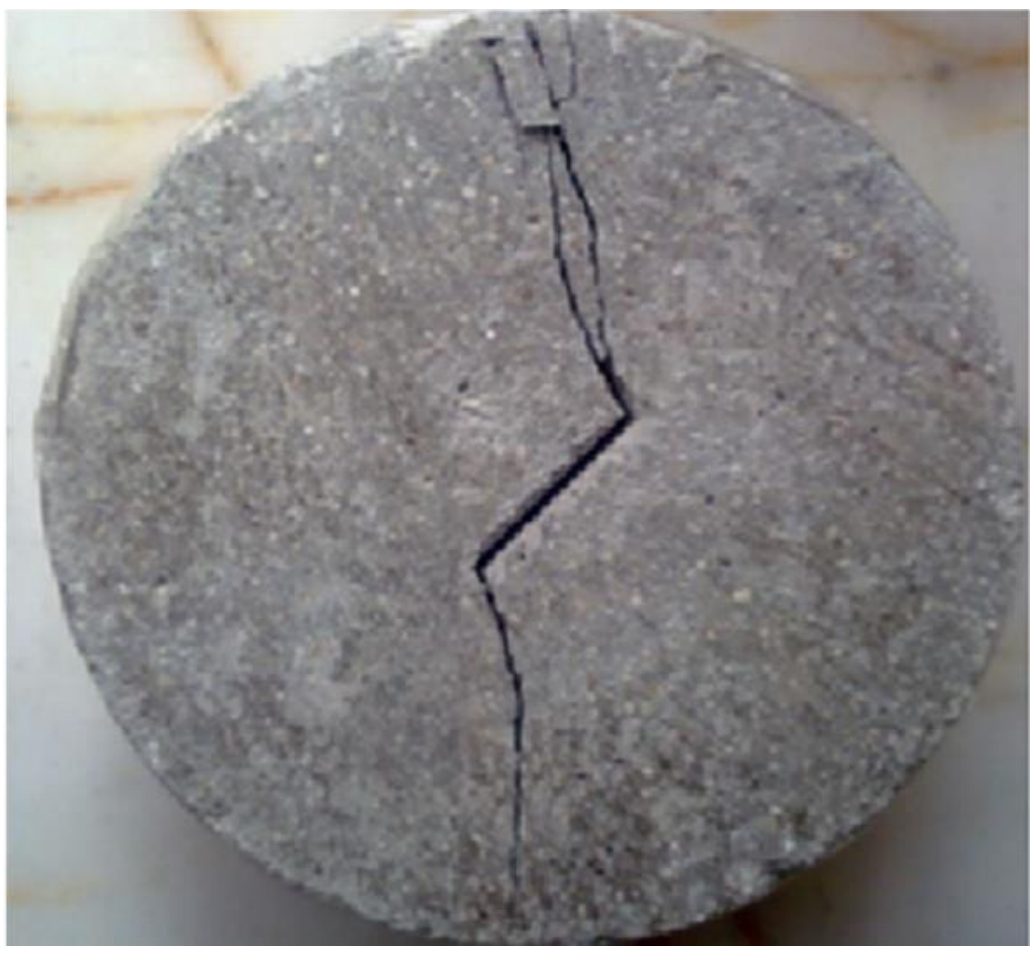

(a)

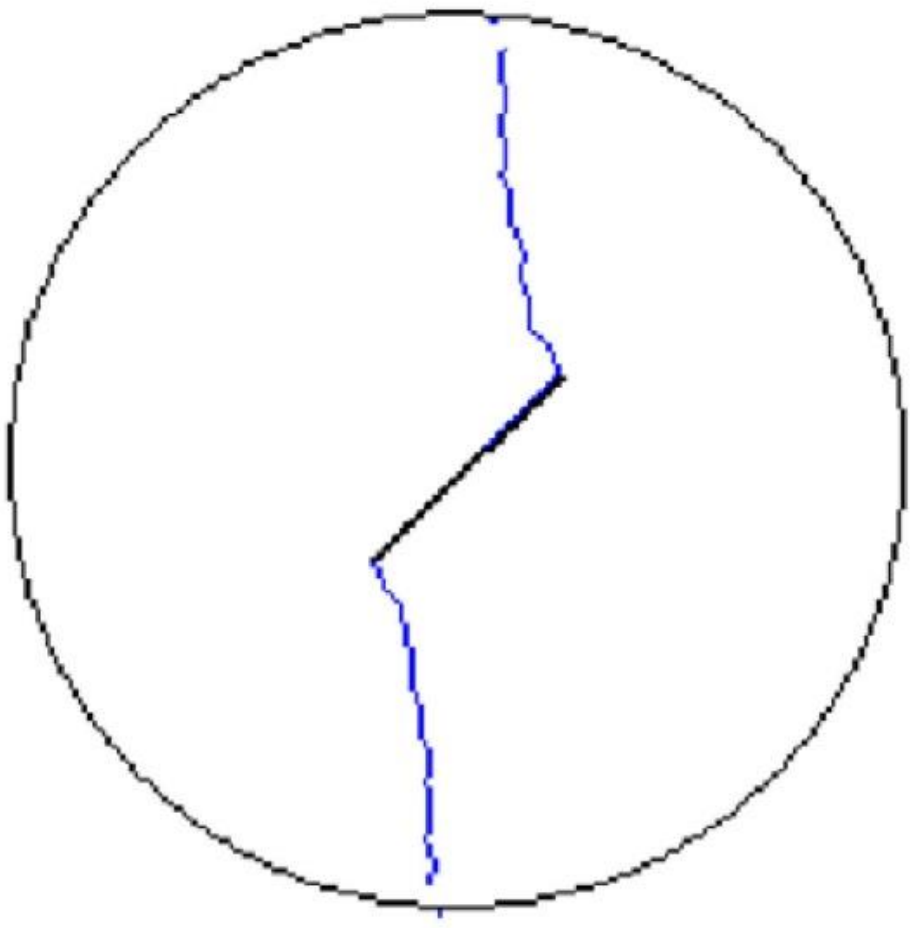

(b)

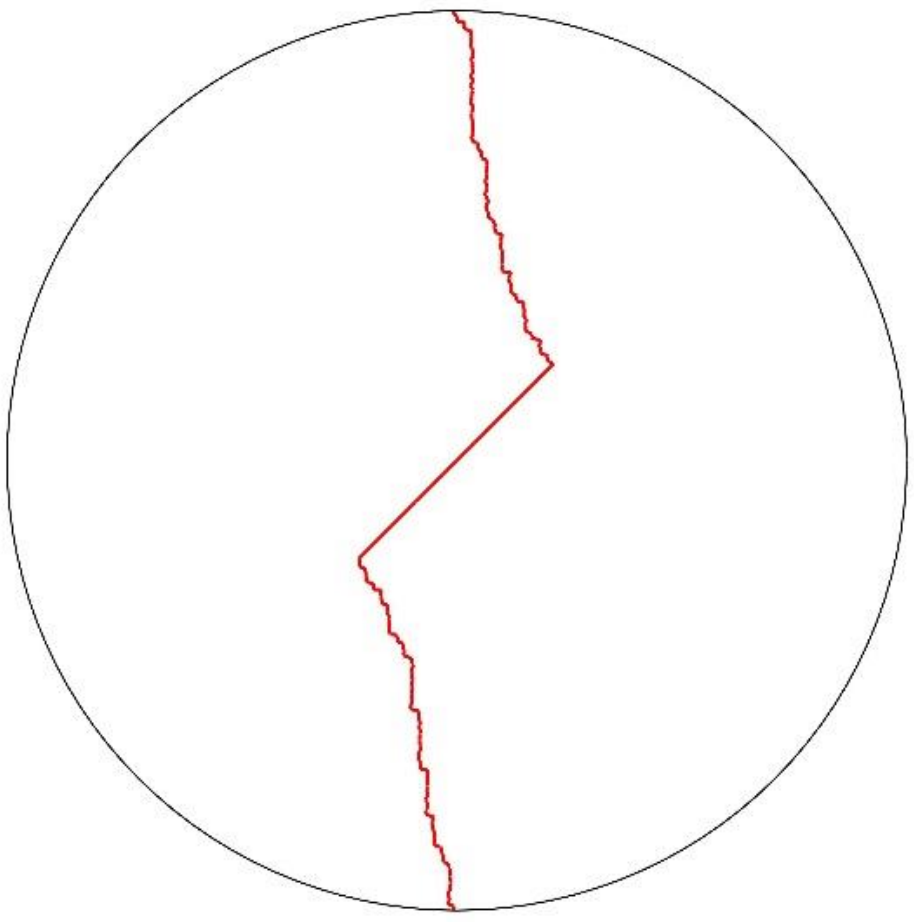

(c)

Figure 39. Crack propagation paths of the pre-cracked disk specimen

(a) Experiment results (b) Boundary Element Method simulations

(c) Fragile Point Method simulations

As discussed in the previous subsection, the computational time and the accuracy of stress analysis by FPM and FEM are close to each other, when the same nodal distribution is used. However, for problems involving crack propagation, FEM necessitates either remeshing the cracked structure in each crack-propagation step (at least separate neighboring elements and insert additional nodes) using software such as Zencrack and Fracn2D, or augmenting the trial functions along the crack path using extended or generalized FEM[21]. In either case, the number of DoFs and the size of the global stiffness matrix will be changed during crack developments. However, in the currently proposed FPM, we just delete the terms related to the IP numerical fluxes of an internal boundary when it is cracked, and adjust the support domain of Points near the crack. In this way, there is no need of separating adjacent elements or adding extra nodes to model crack developments. If the point resolution in the initial subdomain partition is sufficient, crack simulations can keep going with no remeshing. If a very coarse initial mesh is used instead, an adaptive refinement strategy to add Points near the crack tip will be useful to enhance the accuracy of the simulation of near-crack-tip stress fields.

### 4.10 Simulations of the Crack Initiation Process

In this subsection, the FPM is employed to simulate the process of crack initiation and development.

As discussed in Section 3.3, an appropriate traditional-continuum-physics-based criterion is crucially desirable for the practical simulation of crack initiation. However, detailed discussion and judgment of various criteria is a fundamentally important subject itself, which is out of the scope of the present study. In this study, a hoop-stress-based criterion is firstly used because of its simplicity. After that, a new inter-subdomain-boundary bonding-energy-rate based criterion is proposed. And the corresponding numerical simulations of crack initiation and its further developments are demonstrated.

#### 4.10.1 A Hoop-Stress-Based Criterion

To demonstrate the power of FPM for the simulation of crack initiation, square plates with a square hole or a circular hole under biaxial loads is considered, as shown in Figure 40. The width of the plate is equal to 4, the diameter of the circular hole is equal to 2, and the width of the square hole is equal to $\sqrt{2}$. Young's modulus $E$ and the penalty factor $\eta$ are set to be 1, and the Poisson's ratio $v$ is equal to 0.3. A state of Plane Stress is considered. In Figure 40(a)(b), quasi-static biaxial displacements $u_1$ and $u_2$ are gradually applied in the horizontal and vertical directions, respectively. In Figure 40(c)(d), quasi-static biaxial tensile tractions $t_1$ and $t_2$ are applied. The ratio $t_1 : t_2$ or $u_1 : u_2$ is set to be a fixed value.

The entire load history is divided into steps. And in each analysis step, a hoop-stress-based criterion for crack initiation is considered: if the normal traction on a specific internal boundary between two subdomains exceeds a prescribed critical value, this internal boundary will be cracked. In this example, the critical hoop stress value is set to be 1. The crack development results simulated by FPM with simple hoop-stress-based criterion are shown in Figure 41 when $u_1 : u_2 = 0:1,\ 1:1$, and in Figure 42 when $t_1 : t_2 = 0:1,\ 1:1$, respectively.

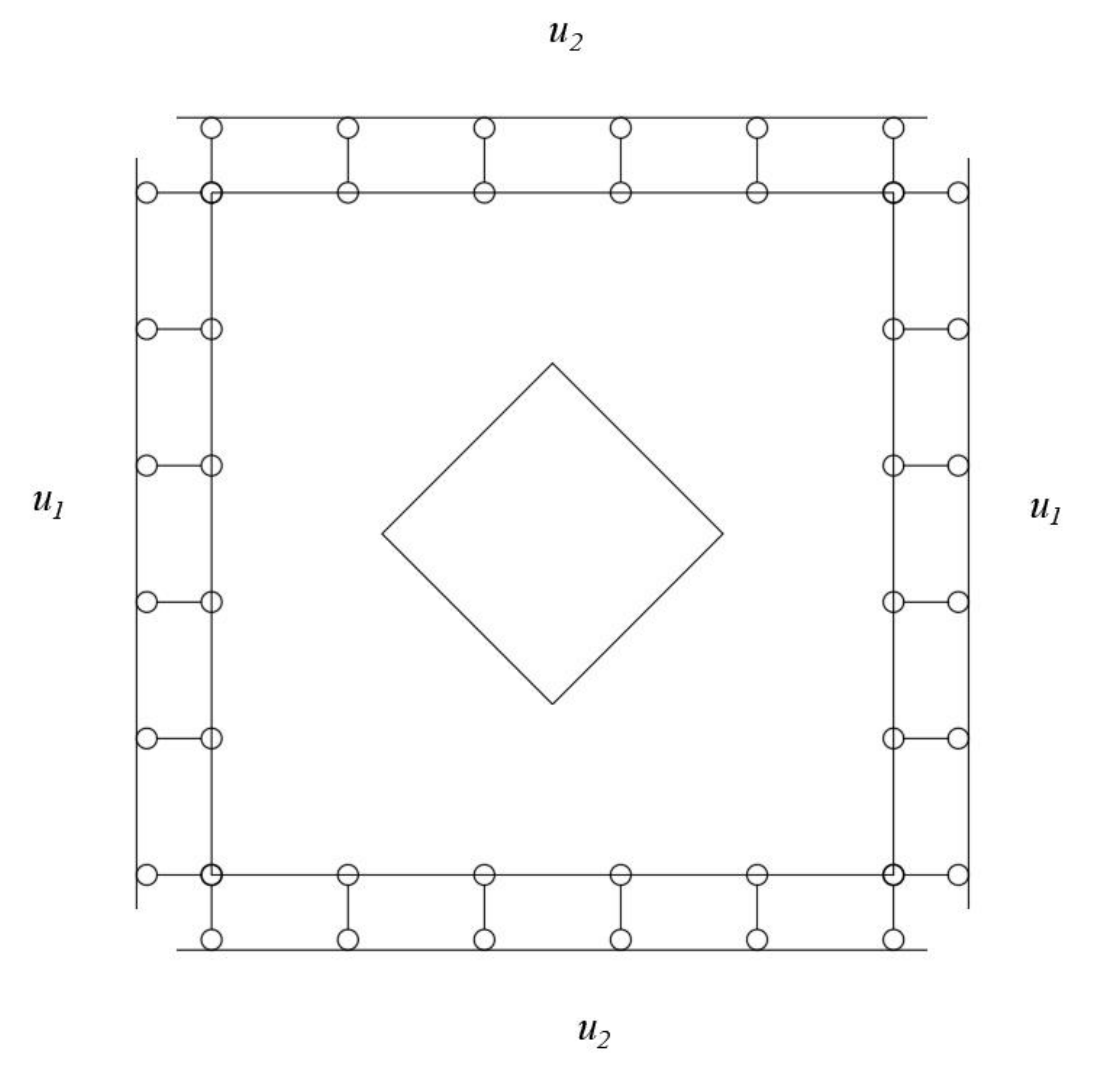

(a)

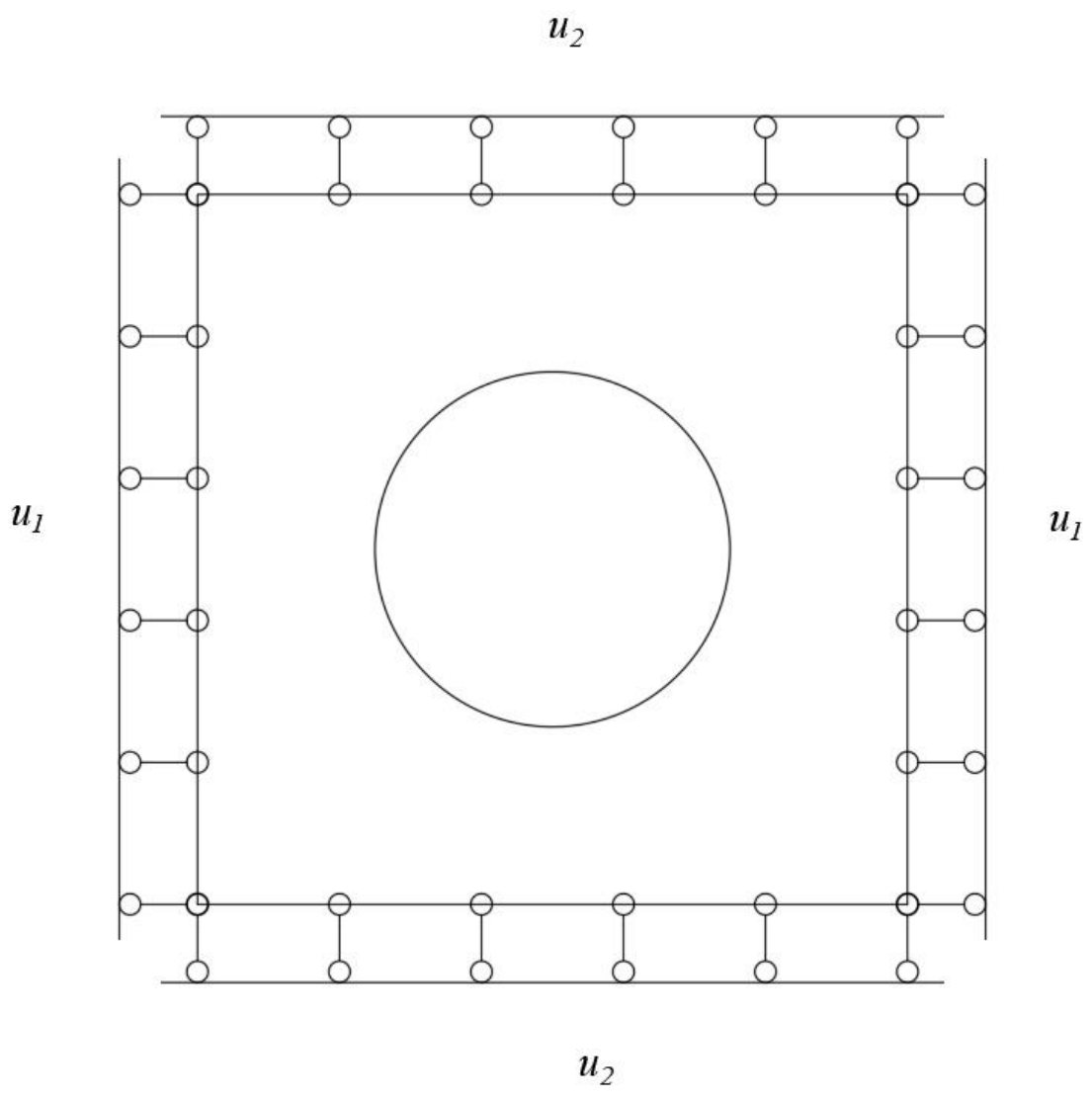

(b)

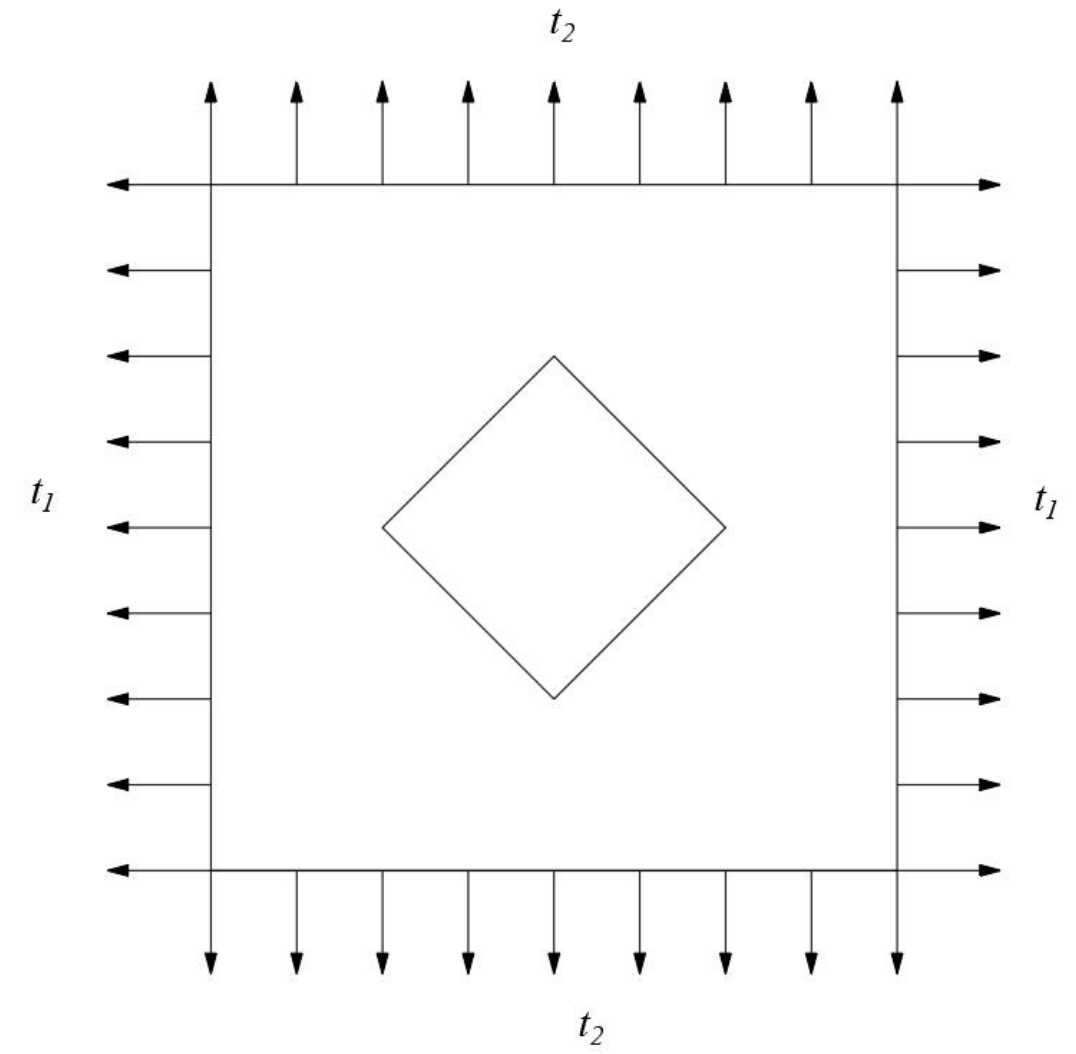

(c)

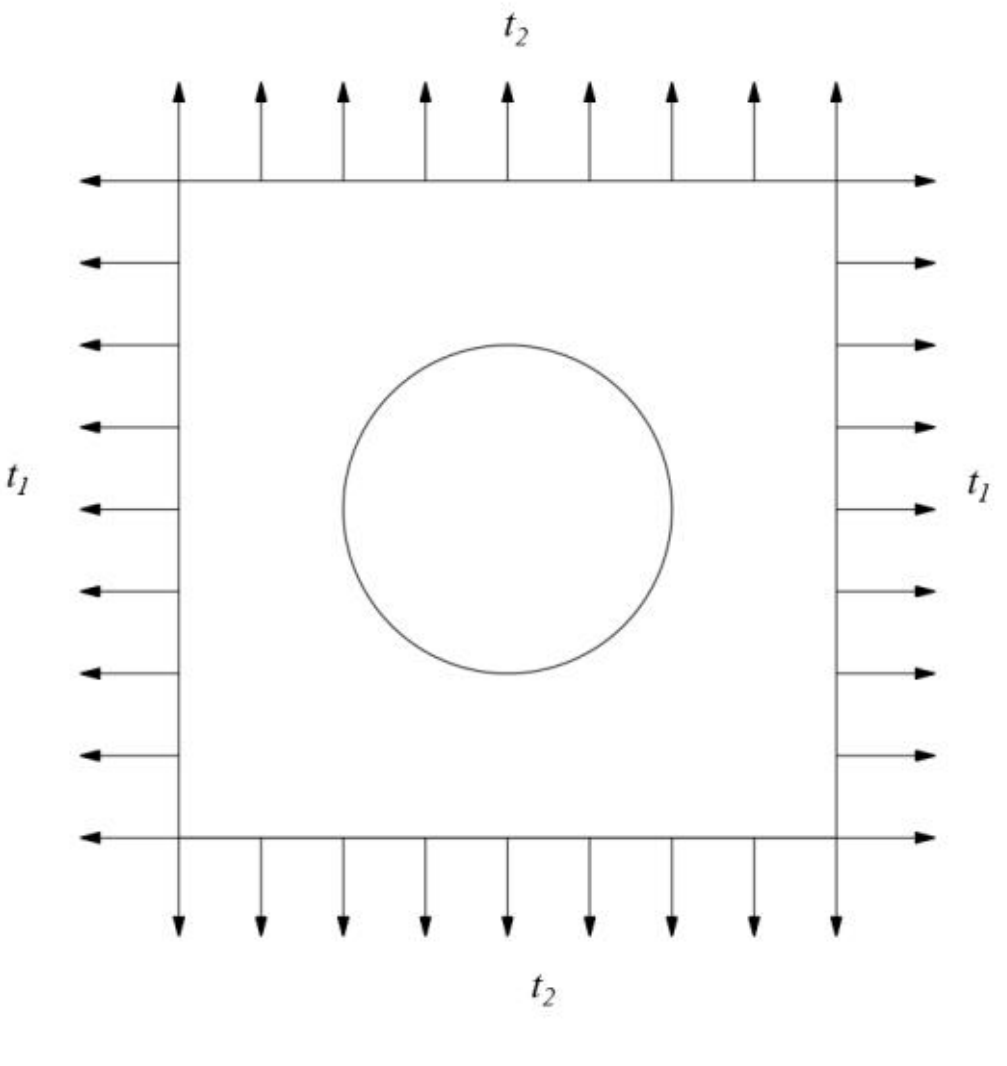

(d)

Figure 40. Square plates loaded

(a) with a square hole, loaded by biaxial displacements

(b) with a circular hole, loaded by biaxial displacements

(c) with a square hole, loaded by biaxial tensile tractions

(d) with a circular hole, loaded by biaxial tensile tractions

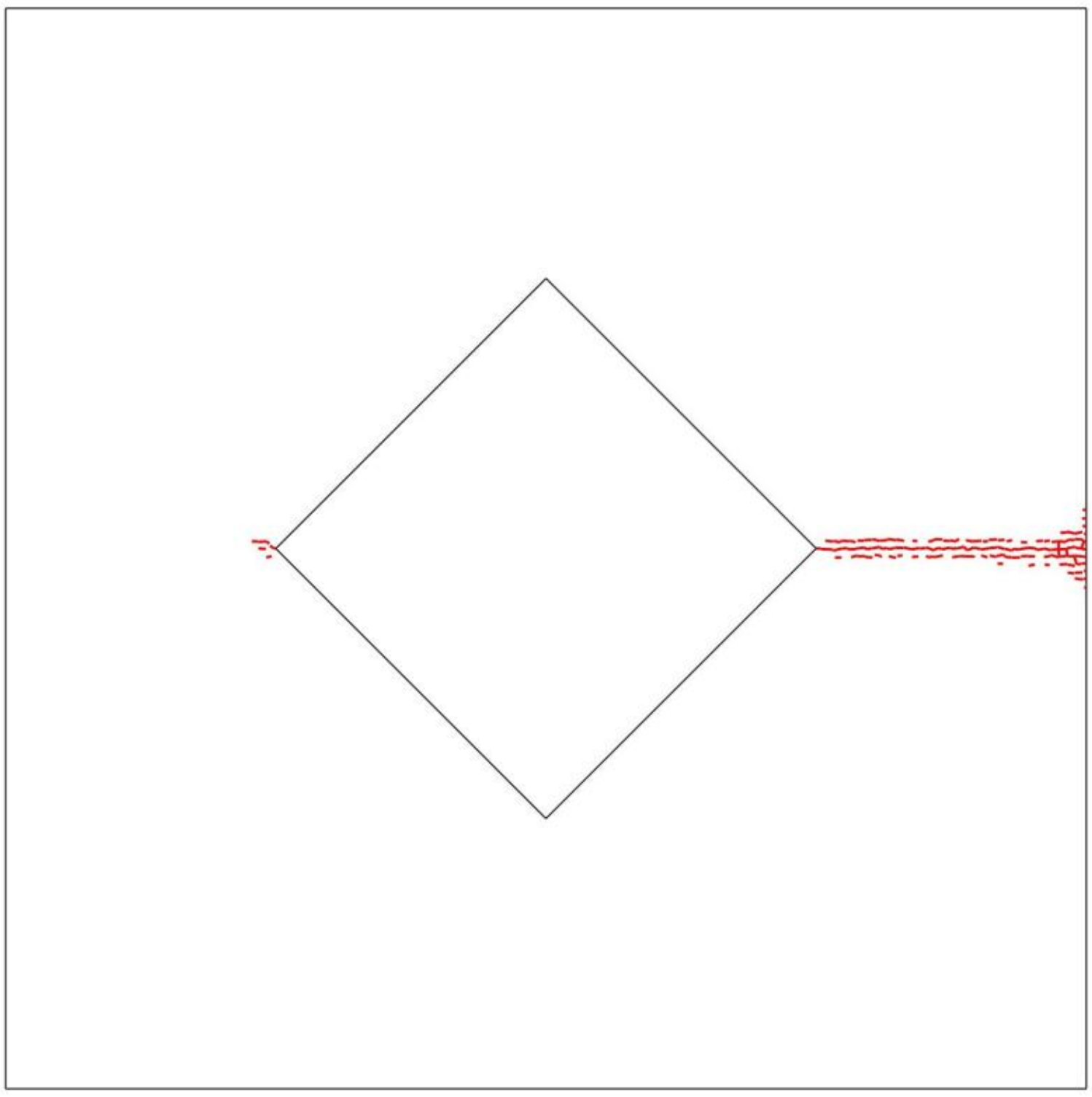

(a)

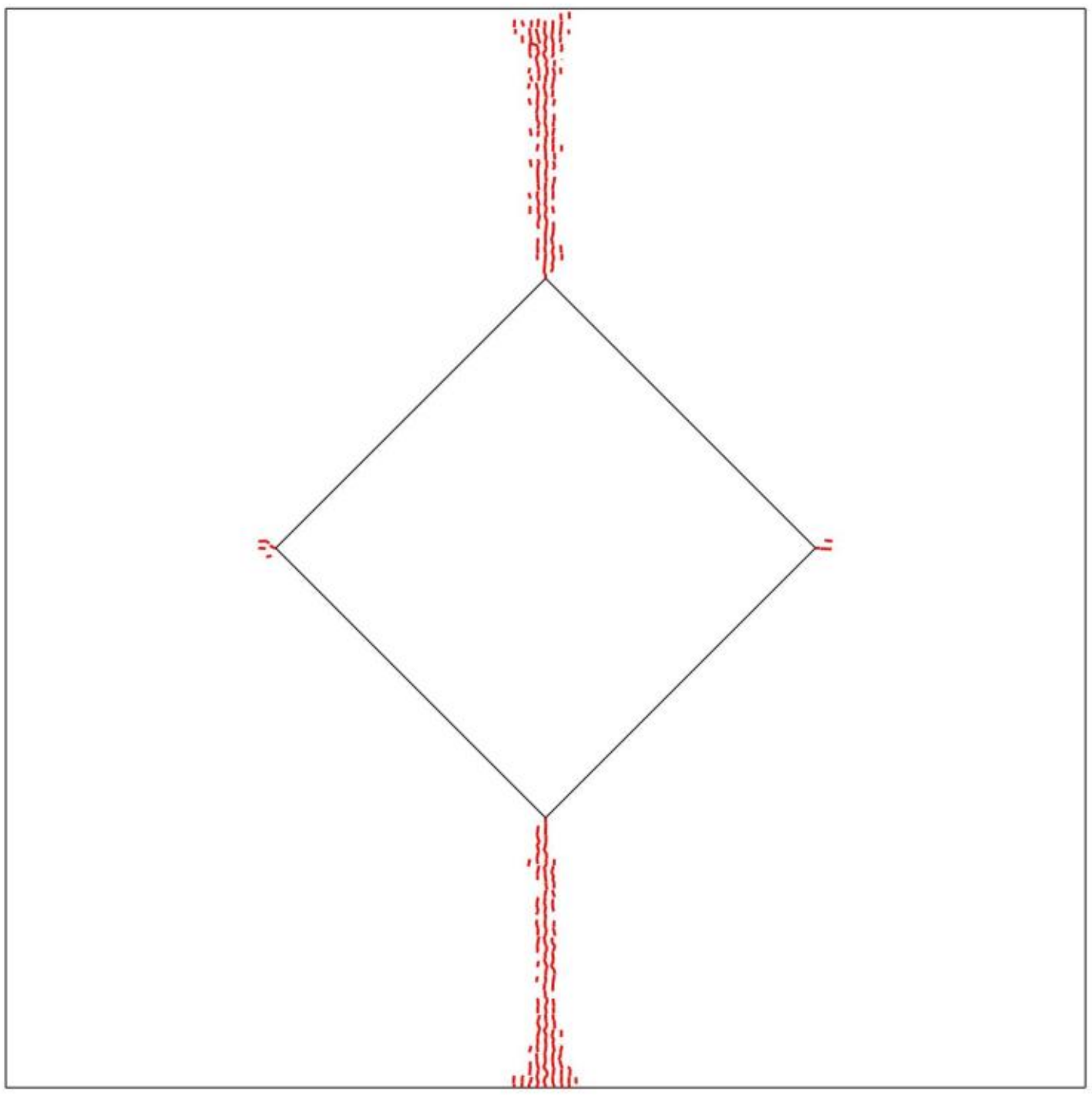

(b)

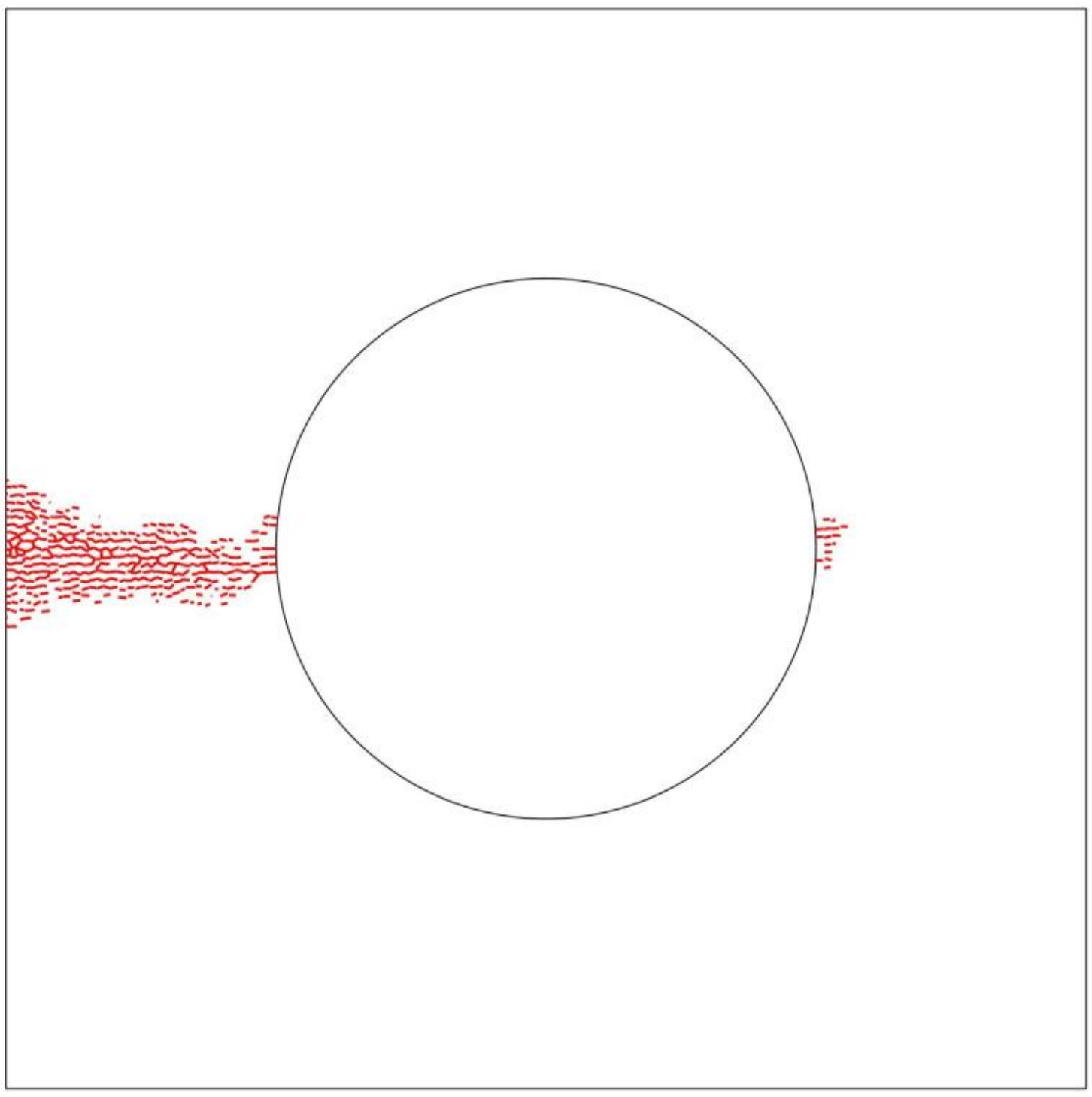

(c)

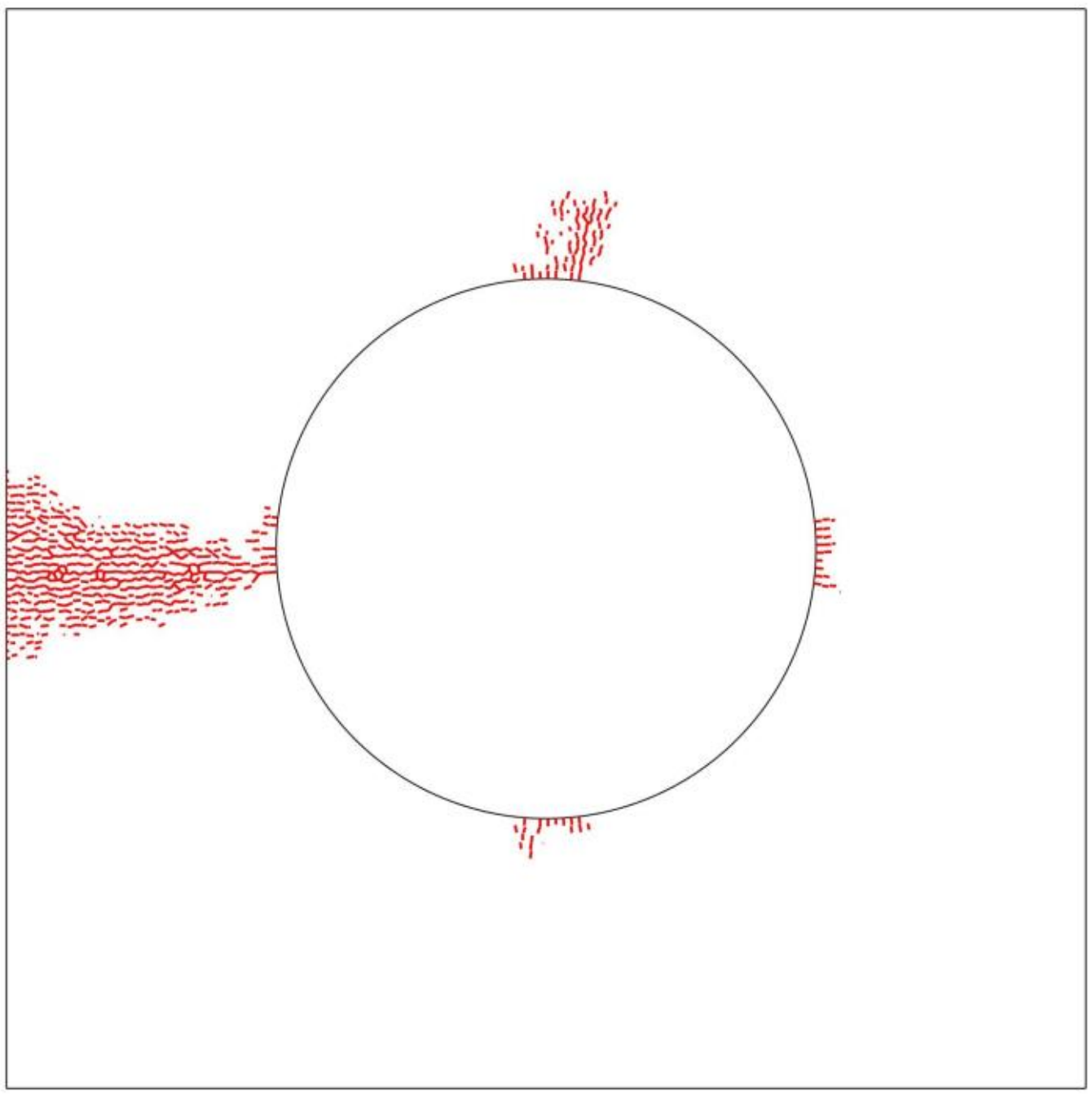

(d)

Figure 41. The crack development results simulated by FPM with simple hoop-stress-based criterion

(a) square hole, $u_1 : u_2 = 0:1$ (b) square hole, $u_1 : u_2 = 1:1$

(c) circular hole, $u_1 : u_2 = 0:1$ (d) circular hole, $u_1 : u_2 = 1:1$

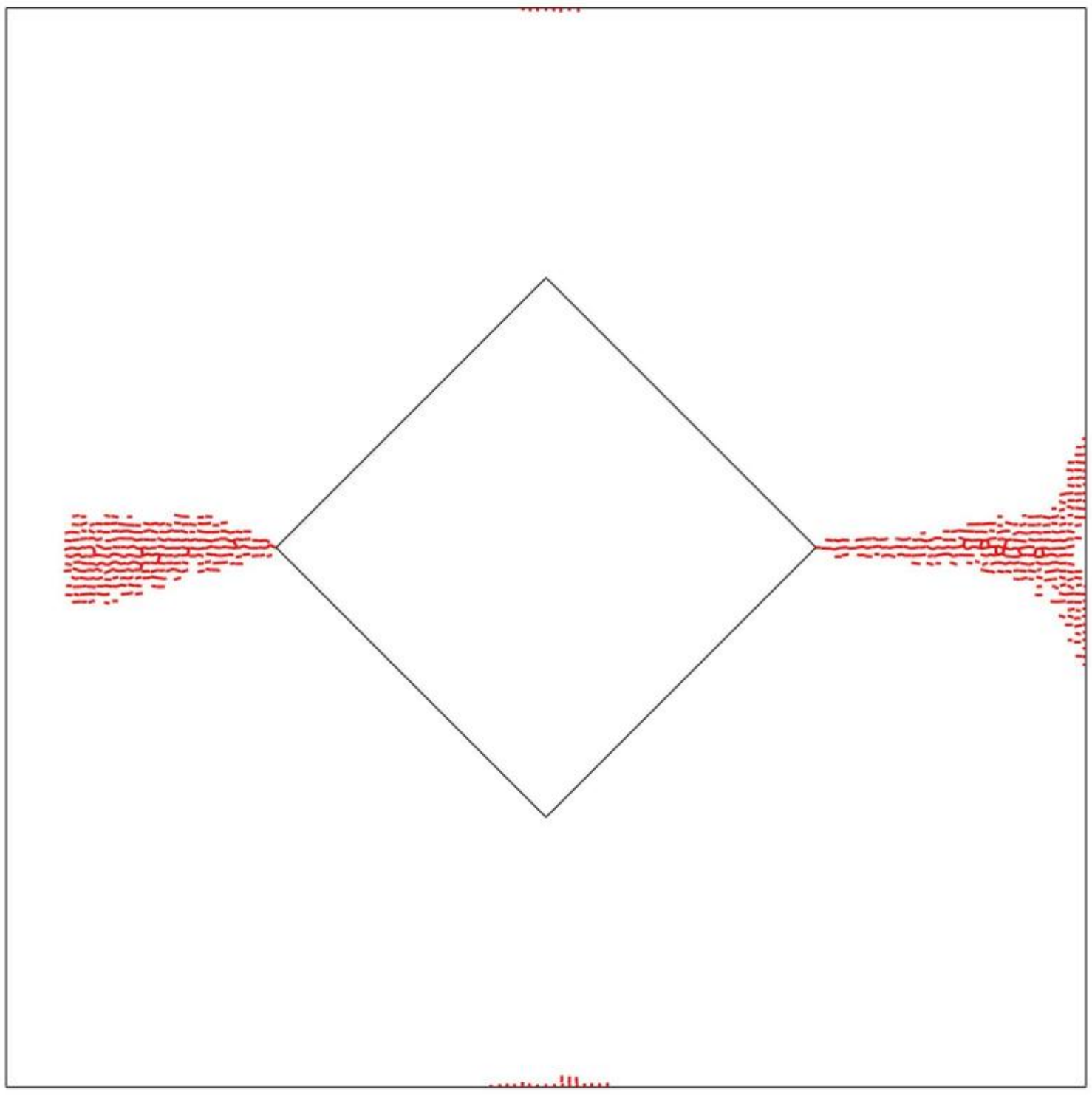

(a)

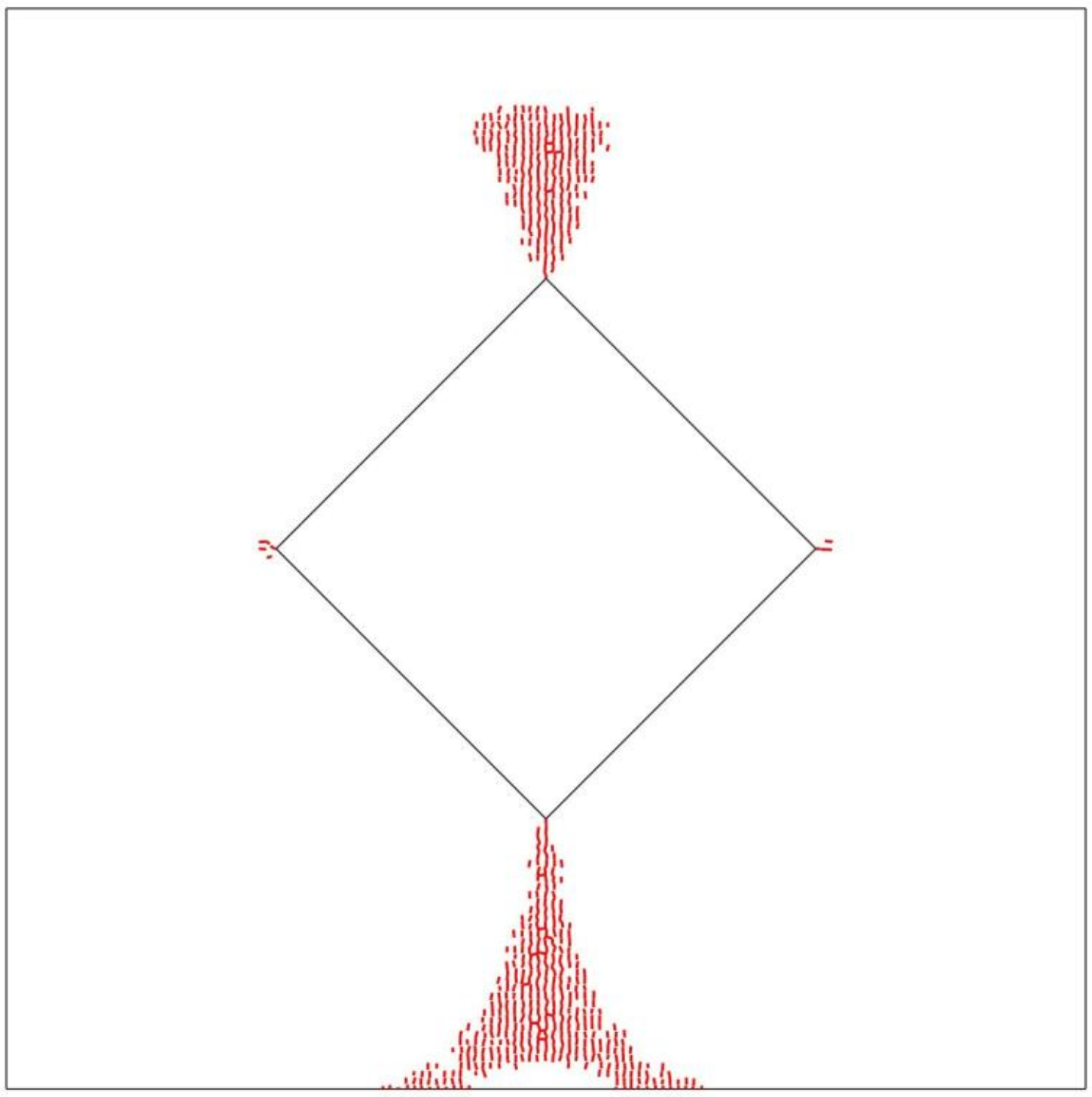

(b)

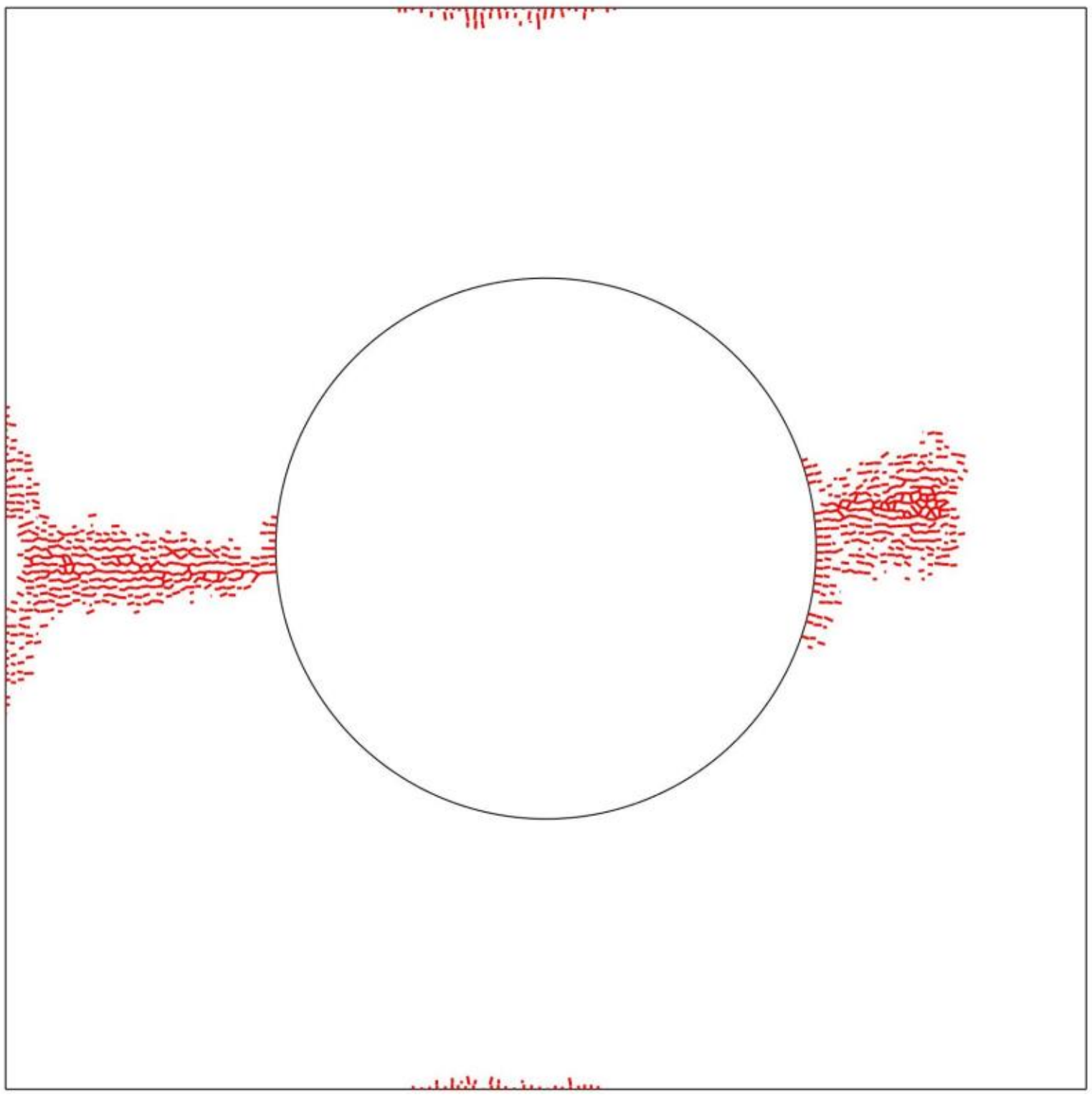

(c)

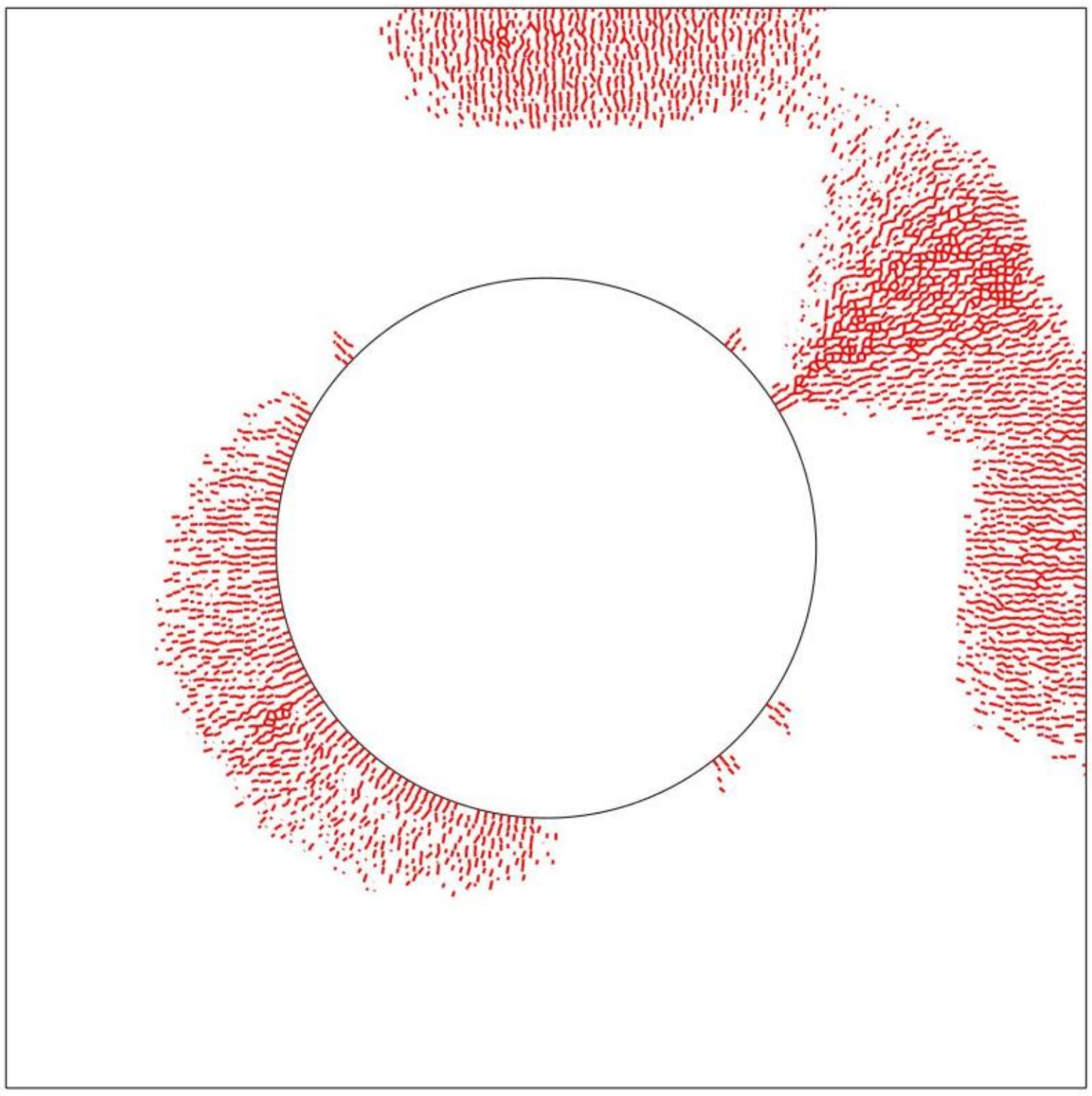

(d)

Figure 42. The crack development results simulated by FPM with the hoop-stress-based criterion

(a) square hole, $t_1 : t_2 = 0:1$ (b) square hole, $t_1 : t_2 = 1:1$

(c) circular hole, $t_1 : t_2 = 0:1$ (d) circular hole, $t_1 : t_2 = 1:1$

4.10.2 A New Inter-Subdomain-Boundary Bonding-Energy-Rate-Based Criterion

A new inter-subdomain-boundary bonding-energy-rate based criterion for crack development is also proposed in this study. Supposing that the body force is negligible, and considering quasi-static loading, we start with the well-known *J*-integral, the definition of which is:

$$J = \int_{\partial\Omega^{\text{int}}} W n_1 d\Gamma - \int_{\partial\Omega^{\text{int}}} \sigma_{ij} n_j \frac{\partial u_i}{\partial \hat{x}_1} d\Gamma \tag{4.12}$$

The Eq. (4.12) is established on a local coordinate system where the $\hat{x}_1$-axis is aligned with the crack. The domain $\Omega^{\text{int}}$ enveloped by the integral contour can be divided into several subdomains $\Omega_I^{\text{int}}$, and $\Gamma_h^{\text{int}}$ is the set of internal boundaries within $\Omega^{\text{int}}$. Then, equivalently (when *W* is a single-valued function of the displacement gradients), we have:

$$\begin{aligned} J &= \sum_I \int_{\partial\Omega_I^{\text{int}}} W n_1 d\Gamma - \int_{\partial\Omega^{\text{int}}} \sigma_{ij} n_j \frac{\partial u_i}{\partial \hat{x}_1} d\Gamma - \sum \int_{\Gamma_h^{\text{int}}} \left[ W n_1^e \right] d\Gamma \\ &= \sum_I \int_{\Omega_I^{\text{int}}} \frac{\partial u_{i,j}}{\partial \hat{x}_1} \sigma_{ij} d\Omega - \int_{\partial\Omega^{\text{int}}} \sigma_{ij} n_j \frac{\partial u_i}{\partial \hat{x}_1} d\Gamma - \sum \int_{\Gamma_h^{\text{int}}} \left[ W n_1^e \right] d\Gamma \end{aligned} \tag{4.13}$$

For the exact solution, Eq. (3.8) should be satisfied over $\Omega^{\text{int}}$ for an arbitrary test function $v_i$. And for the approximate solution $u_i$ given by FPM, by letting $v_i = \frac{\partial u_i}{\partial \hat{x}_1}$, we may argue that Eq. (3.8) is satisfied approximately. Therefore, combining Eq. (3.8) and Eq. (4.13), and considering $v_i = \frac{\partial u_i}{\partial \hat{x}_1}$, we have:

$$\begin{aligned} & J \approx \sum_{e \in \Gamma_h^{\text{int}}} BER; \\ & BER = -\int_e \frac{\eta}{h_e} \left[ \frac{\partial u_j}{\partial \hat{x}_1} \right] \left[ u_j \right] d\Gamma + \int_e \left\{ n_i^e \sigma_{ij} \right\} \left[ \frac{\partial u_j}{\partial \hat{x}_1} \right] d\Gamma \\ & + \int_e \left\{ n_i^e \frac{\partial \sigma_{ij}}{\partial \hat{x}_1} \right\} \left[ u_j \right] d\Gamma - \int_e \left[ W n_1^e \right] d\Gamma. \end{aligned} \tag{4.14}$$

From Eq. (4.14), we can see that the *J*-integral is approximately equal to the summation of integrals over internal boundaries, which we define as *BER*. Thus, we may postulate that *BER* has the physical meaning of the bonding energy rate in the set-up which is very specific to the currently developed algorithm of FPM. Moreover, if we consider one specific internal boundary, and define the local coordinate system where the $\hat{x}_1$-axis is aligned with this internal boundary segment, then we have $n_1^e = 0$. Besides, when a linear trial function is employed, $\frac{\partial \sigma_{ij}}{\partial \hat{x}_1} = 0$ in each subdomain. As a result, the third and fourth terms of Eq. (4.14) vanish. Therefore, the formula for *BER* of this specific internal boundary-segment is simplified as:

$$BER = -\int_e \frac{\eta}{h_e} \left[ \frac{\partial u_j}{\partial \hat{x}_1} \right] \left[ u_j \right] d\Gamma + \int_e \left\{ n_i^e \sigma_{ij} \right\} \left[ \frac{\partial u_j}{\partial \hat{x}_1} \right] d\Gamma \tag{4.15}$$

In this study, we propose that *BER* can be used as one of the possible energy-based criteria to simulate the crack initiation. The BER defined in Eq. (4.15) may ab initio be also used as a criterion for dynamic problems and for arbitrary material behavior. The same plate with a circular or square hole is considered, which is loaded by biaxial displacements or tensile tractions. And in each load step, if the *BER* on a specific internal boundary between two subdomains exceeds a prescribed critical value, this internal boundary will be cracked. In this paper, the critical *BER* value is set to be 1. The results simulated by FPM with this *BER*-based criterion for the plated loaded by biaxial displacements and tensile tractions are shown in Figure 43 and Figure 44, respectively.

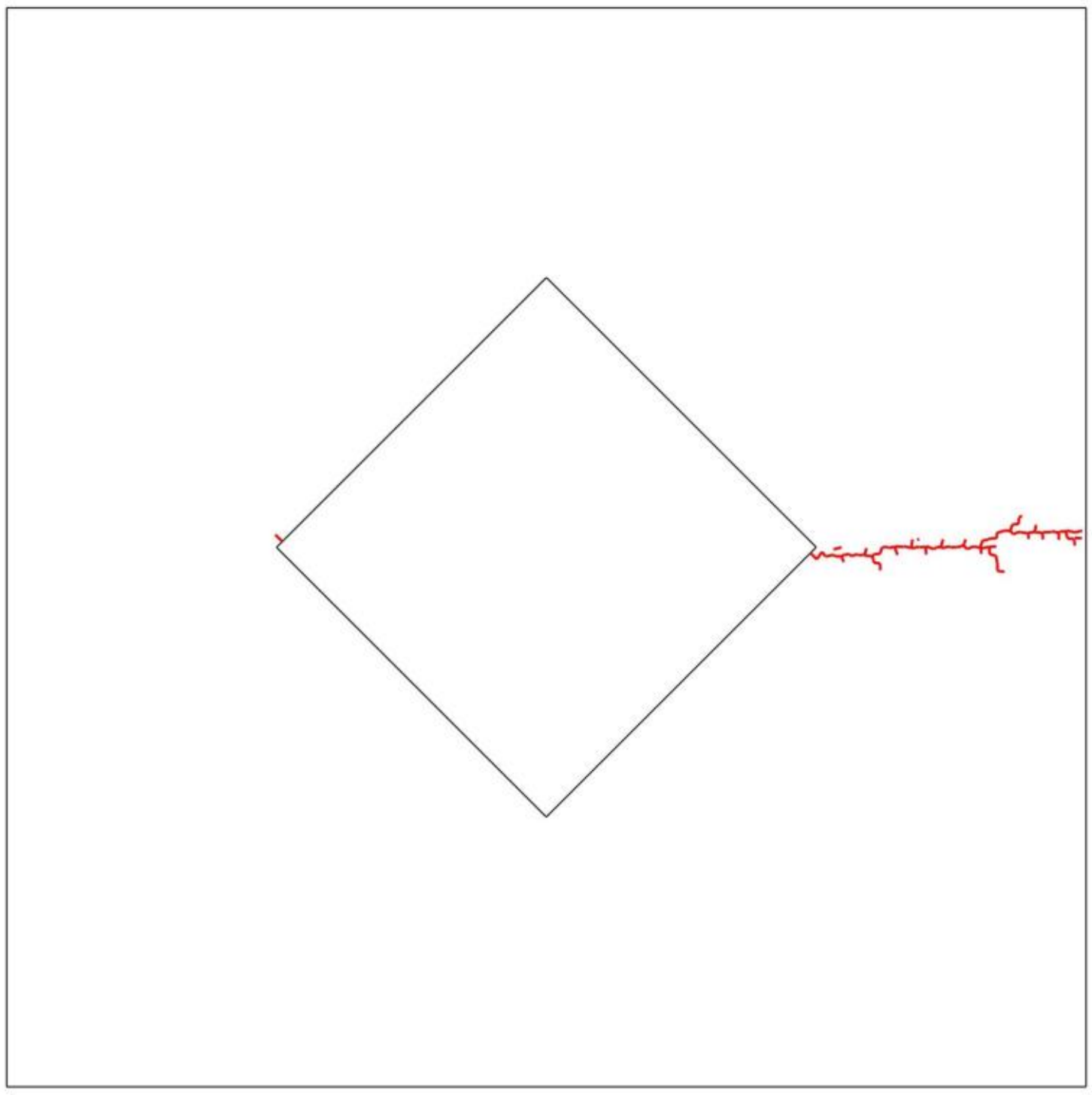

(a)

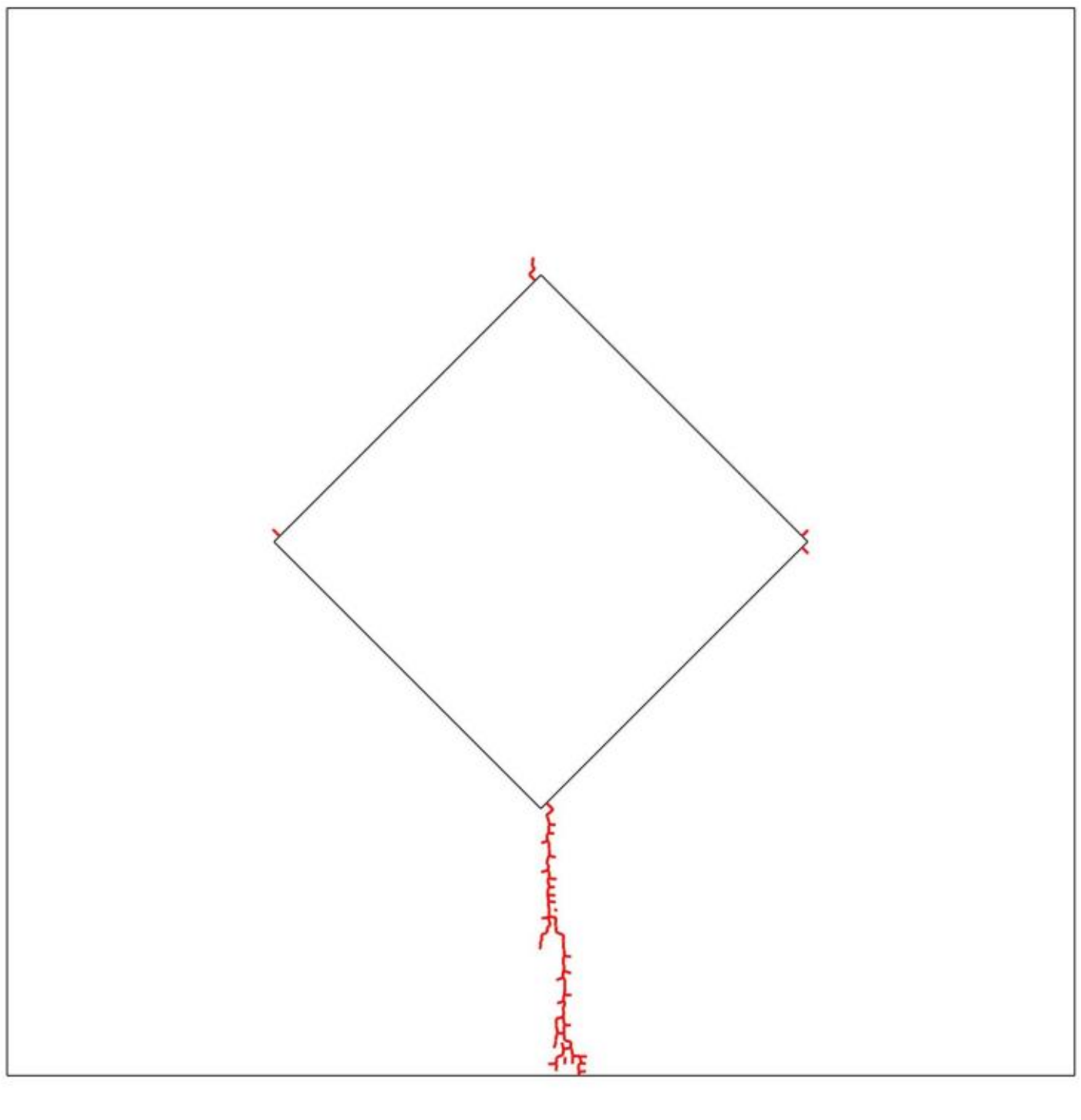

(b)

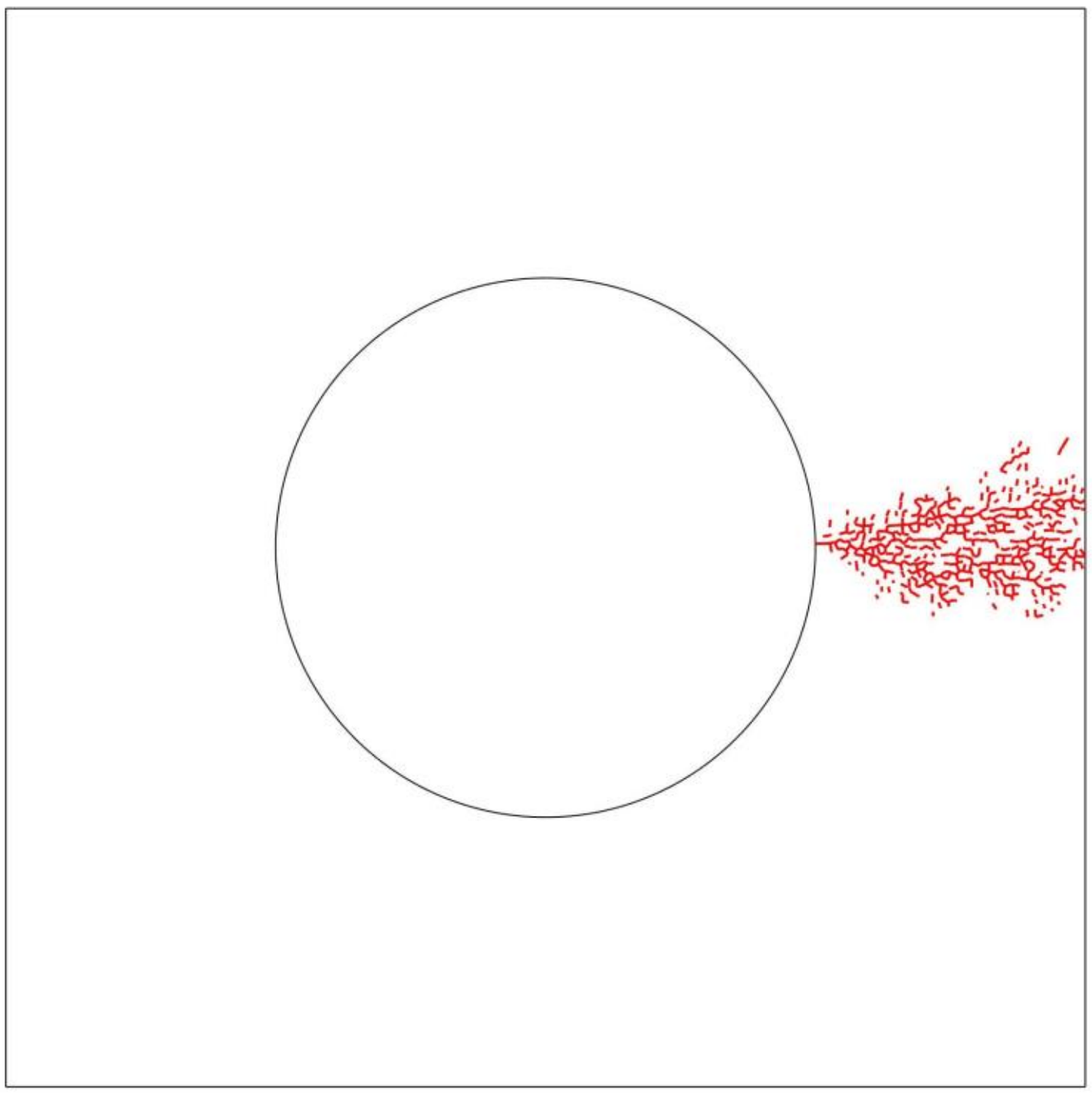

(c)

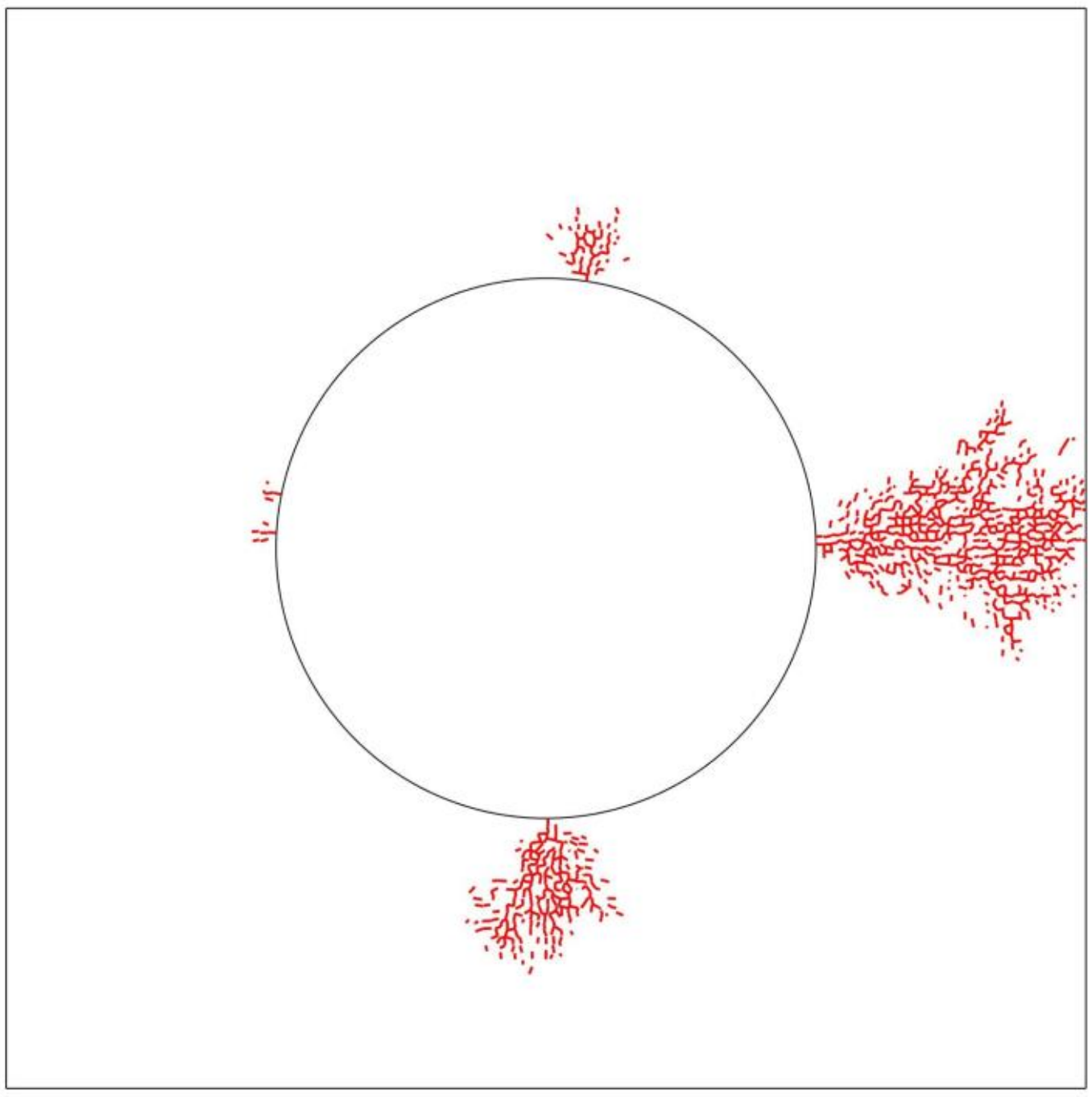

(d)

Figure 43. The crack development results simulated by FPM with the *BER*-based criterion

(a) square hole, $u_1 : u_2 = 0:1$ (b) square hole, $u_1 : u_2 = 1:1$

(c) circular hole, $u_1 : u_2 = 0:1$ (d) circular hole, $u_1 : u_2 = 1:1$

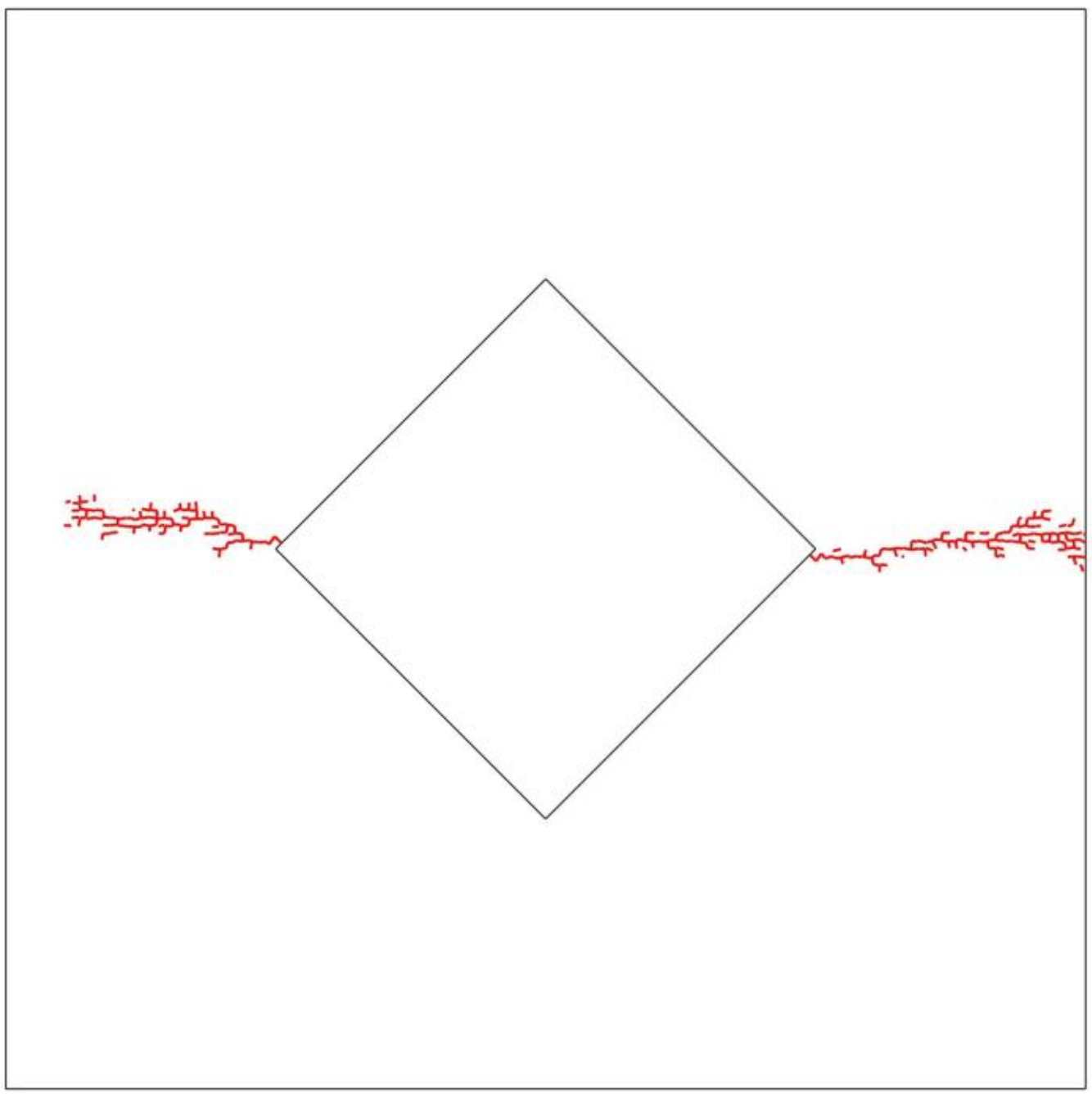

(a)

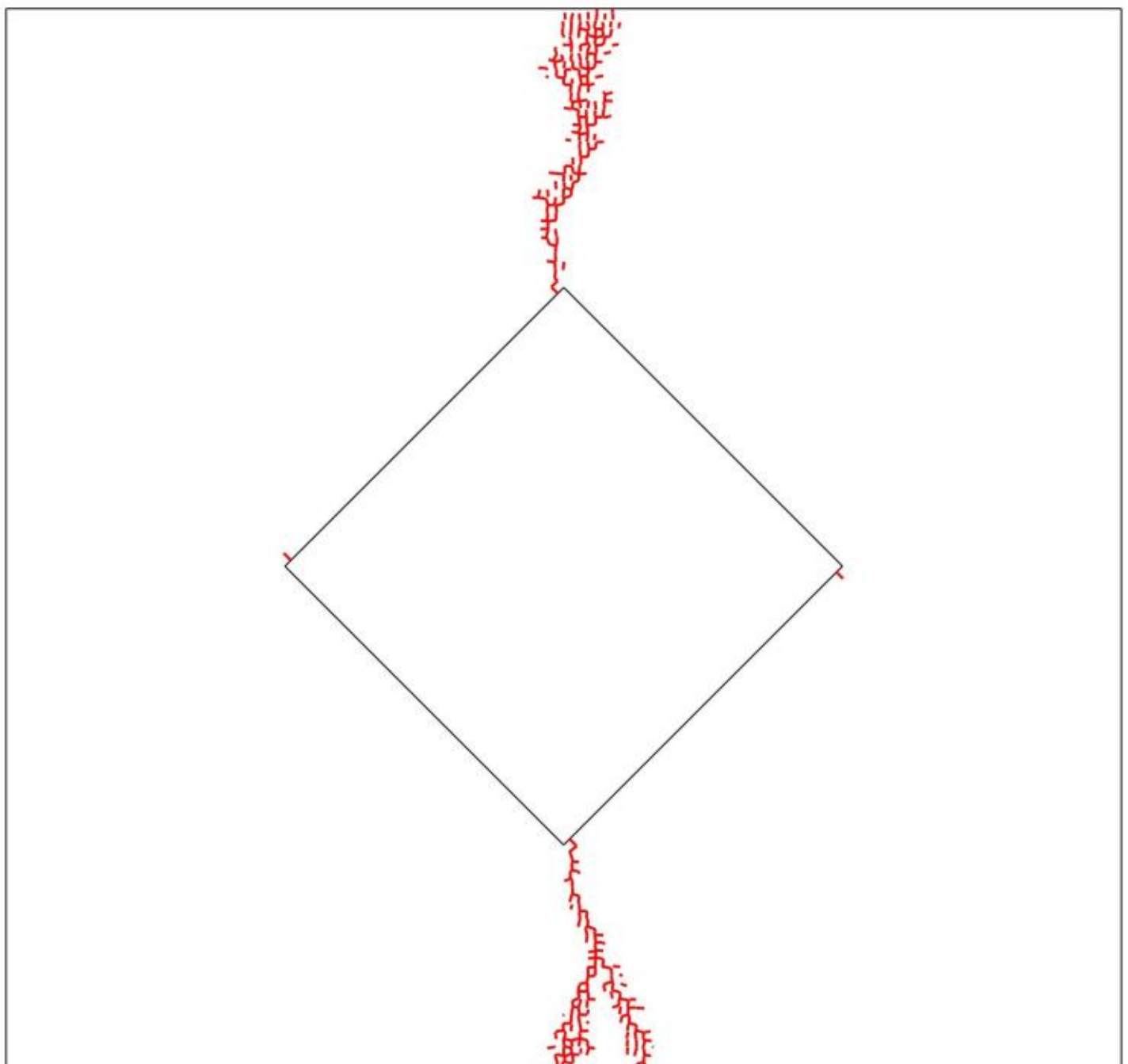

(b)

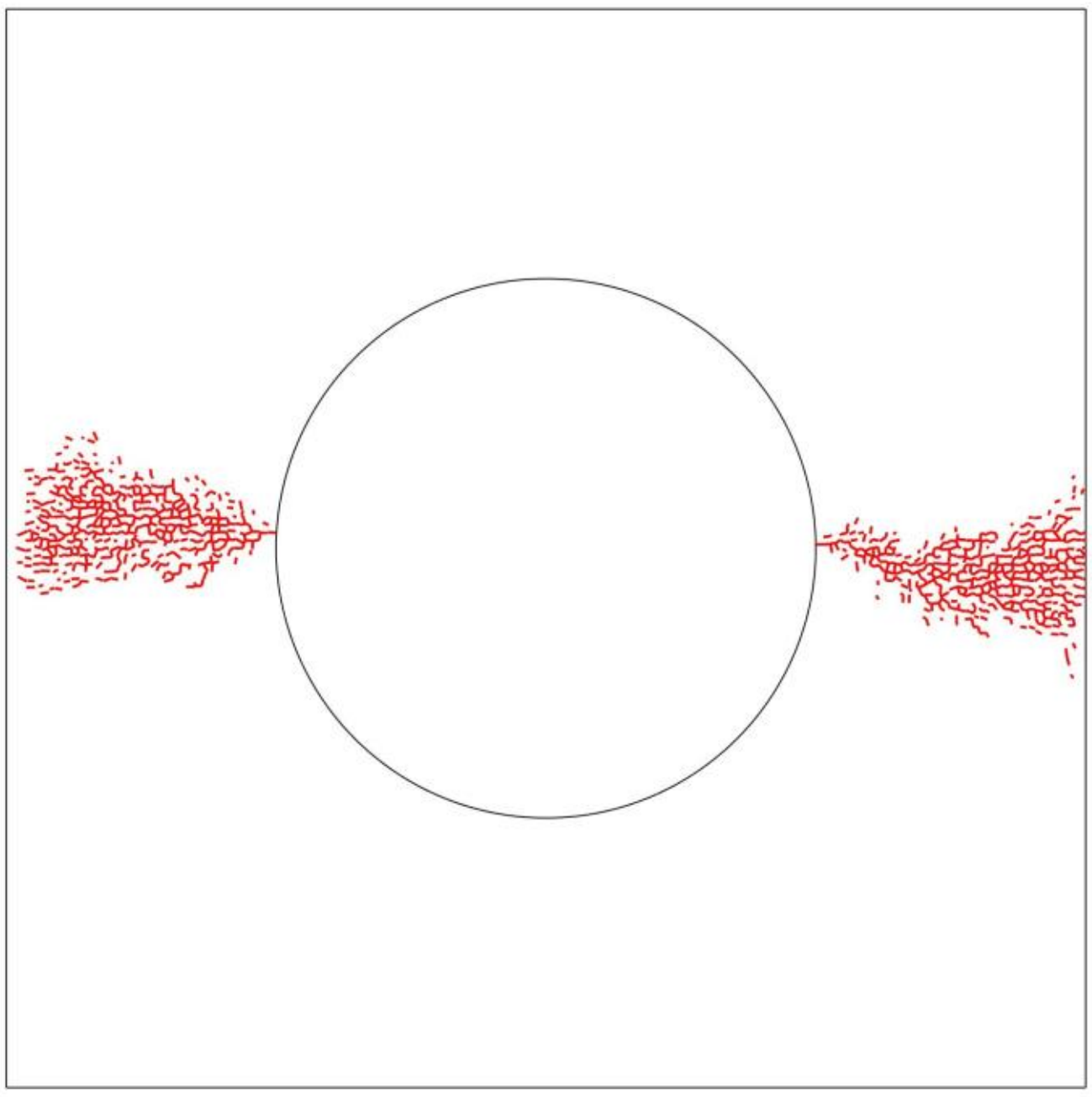

(c)

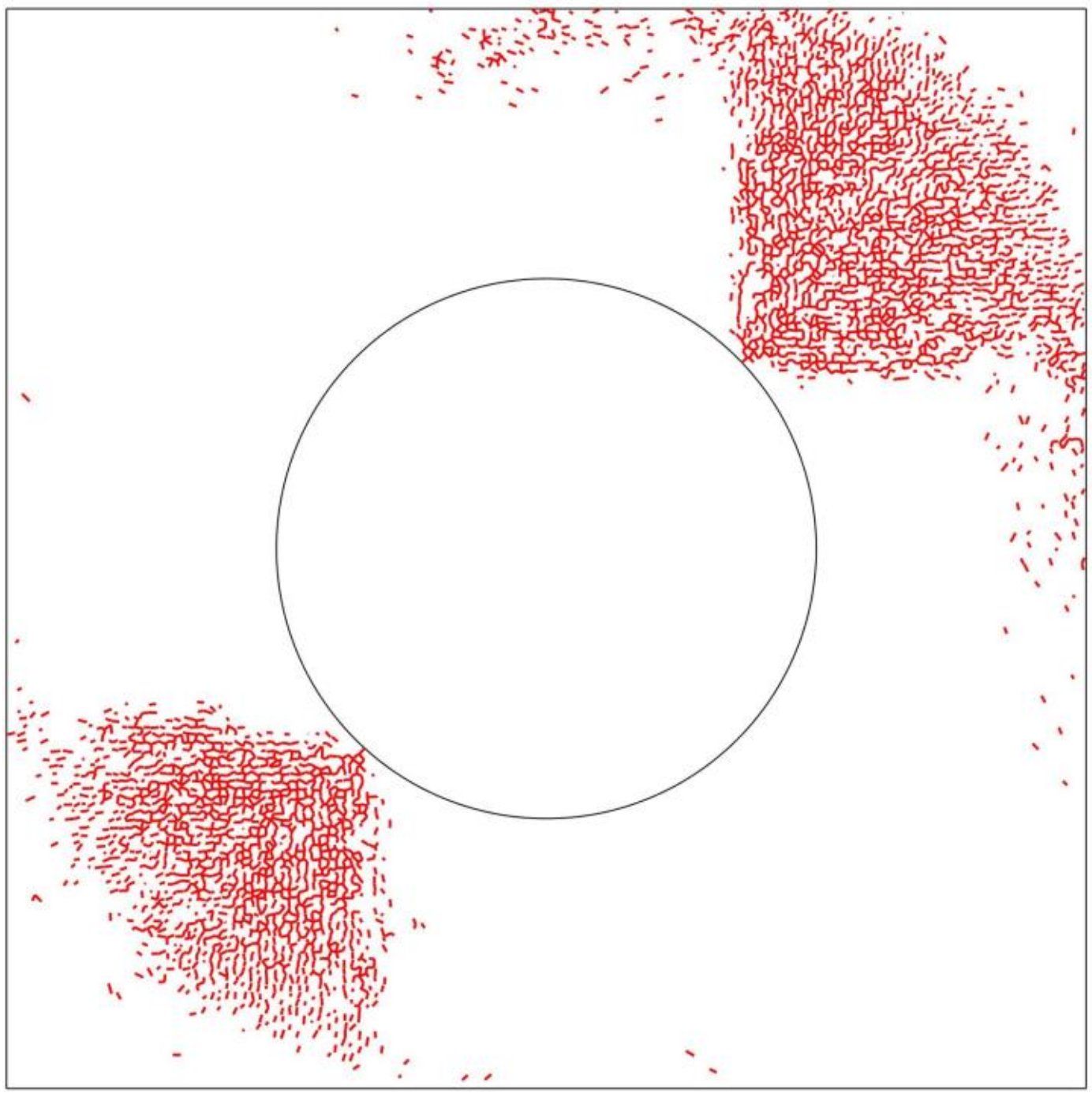

(d)

Figure 44. The crack development results simulated by FPM with the *BER*-based criterion

(a) square hole, $t_1 : t_2 = 0:1$ (b) square hole, $t_1 : t_2 = 1:1$

(c) circular hole, $t_1 : t_2 = 0:1$ (d) circular hole, $t_1 : t_2 = 1:1$

From the simulated results by FPM with the hoop-stress-based criterion and the *BER*-based criterion, we can see that the distribution of cracks is more diffuse with load control. This is expected as crack initiation will increase stress and energy concentration, therefore once the first crack-segment is initiated, more inter-subdomain-boundaries will have local fields exceeding the critical stresses and the critical *BER*. Thus, inter-subdomain-boundaries will be cracked in a sequential way under load control, as soon as the first crack is initiated, which eventually will form a diffuse pattern of cracks. Moreover, it can also be seen that, rupture develops from corners of square plates in a pattern similar to clear-cut single cracks, while diffuse cracks evolve from circular holes in contrast. This is expected with the current hoop-stress-based criterion and the *BER*-based criterion, because the stress and energy are more concentrated in the plate with a square hole as compared to the plate with a circular hole.

It should be noted that, the critical magnitude of *BER* may be measured through a "Hybrid Experimental-Numerical Approach"[22]. An experiment on a center-cracked specimen or a double-cantilever specimen can be conducted to measure load and point-displacement versus crack-growth data. The experiment data can be simulated in an FPM simulation so as to compute the *BER*. In this way, the critical magnitude of *BER* can be measured indirectly by the "Hybrid Experimental-Numerical Approach", which can be then used to model and design of complex engineering structures. Relevant studies will be conducted in our future work.

## 5. Conclusions

In this paper, we have formulated the algorithmic framework and details of the "Fragile Points Method" (FPM) for linear elasticity. Discontinuous, piecewise-polynomial and Point-based trial and test functions are constructed in the FPM and Numerical Flux Corrections are introduced to resolve the inconsistency caused by the discontinuity. In FPM, a sparse, symmetric and positive definitive global stiffness matrix can be obtained by efficiently computing and assembling Point Stiffness and Boundary Stiffness Matrices. The convergence, robustness, consistency and high accuracy of the FPM in elasticity problems are demonstrated by benchmark as well as realist problems. Moreover, we successfully used the FPM for simulations of crack propagation and initiation with different criteria.

In conclusion, the salient features of this new method are:

1. In FPM, one can randomly sprinkle the Fragile Points in the problem domain and use Voronoi Cell partitioning around these Points to fill up the problem domain. Alternatively, the Points are taken to be the centroids of the arbitrarily meshed "Elements" generated through, for example, ABAQUS or ANSYS. However, we do not use the usual Finite Element (FEM) methodology of ABAQUS, ANSYS etc. The trial functions within each of the Voronoi Cells, or alternatively the "Elements", are derived from the Taylor Series expansions of these functions around each Point in the Cell or "subdomain". The derivatives in the Taylor series at any Point are obtained using a variety of ways such as the method of Generalized Finite Differences, Compactly Supported Radial Basis Functions, Differential Quadrature, etc., in terms of the values of the functions only at a finite set of neighboring FPM Points in the problem domain;

2. Depending on the order of the Taylor series expansions, we can generate constant, linear and higher order strains in in each Voronoi Cell or "Element". Thus FPM is a meshless or element-free Galerkin Method, in its simplest form;

3. These Point-based trial and test functions are, however, locally discontinuous at the boundaries of Voronoi Cells or Elements, and these discontinuities are remedied by using Interior Penalty Numerical Fluxes;

4. Unlike the Element Stiffness Matrices in FEM, the FPM leads naturally to Point Stiffness Matrices, and their Corrections due to the Numerical Flux. The bandwidth of each Point Stiffness Matrix depends on the number of surrounding Points which influence the derivatives at each Fragile Point;

5.The integration of the Galerkin weak-form in each Voronoi Cell or "Element" is performed exactly by using only a one-point integration, when the Taylor series expansion contains only first-order derivatives. When second order derivatives are used in the Taylor series, only a simple 2 by 2 integration in each “Element” or triangulation-based three-point Hammer integration in each polygonal Cell is sufficient;

6. There is no incompressibility locking or shear locking in FPM, and thus there is no hourglass control or artificial viscosity;

7. Since there are no element-based interpolations as in FEM, FPM is not sensitive to mesh distortion as in FEM;

8. Crack and rupture initiation and propagation involve neither any remeshing if there are sufficient Points initially, nor trial function enrichment as the crack propagates, as in Generalized FEM, XFEM, Zencrack, etc.;

9. The entire algorithm of the Fragile Points Method is based only on simple continuum physics/mechanics without invoking any nonlocal theories, etc. Crack/rupture initiation/propagation is modeled by simply releasing the segment of the boundary that connects two Voronoi Cells or Elements surrounding two neighboring Points, when certain fracture criteria are met. This leads to the name " Fragile Points Method" due to its ability to simply model the phenomena of material "Fragility";

10. Problems with first and higher order strain gradient effects can be handled in the same way with simplicity, with only displacements as the degrees of freedom at the Points. Likewise, extreme problems of large deformation, plasticity, impact, penetration, and fragmentation can also be handled easily with FPM. The analysis of such problems will be presented in our forthcoming papers.

## Acknowledgements

The first two authors thankfully acknowledge the support from the National Key Research and Development Program of China (No. 2017YFA0207800) and Beihang Advanced Discipline Center for Unmanned Aircraft System (ADBUAS-2019-SP-05). The authors benefited from the constructive criticisms of anonymous reviewers.